\documentclass[11pt,a4paper]{article}
\usepackage{amsfonts}
\usepackage{amsmath}
\usepackage{amsbsy}
\usepackage{theorem}
\usepackage{graphicx}

\usepackage{mathtools}
\usepackage{multirow,array}
\usepackage{graphicx}
\usepackage{subfigure}
\usepackage{epstopdf}
\usepackage{float} 
\usepackage{color, xcolor}
\newtheorem{definition}{Definition}
\newtheorem{theorem}{Theorem}

\newtheorem{corollary}{Corollary}
\newtheorem{remark}{Remark}
\allowdisplaybreaks

\begin{document}

\sf
\title{\bf Global Fr\'echet regression from time correlated bivariate curve data in manifolds
}
\date{}

 \maketitle


\author{A. Torres--Signes$^{1},$ M.P. Fr\'{\i}as$^{2}$ and M.D. Ruiz-Medina$^{3}$}

\noindent{\small $^{1}$\ Department of Statistics and Operation Research, University of M\'alaga,\\
$^{2}$  Department of Statistics and Operation Research, University of Ja\'en\\
$^{3}$ Department of Statistics and Operation Research, University of Granada}

\begin{abstract}

Global Fr\'echet regression is addressed   from  the observation  of a strictly stationary   bivariate curve process, evaluated in a finite--dimensional compact differentiable Riemannian manifold, with bounded positive smooth sectional curvature. The involved univariate curve processes respectively define the  functional response and regressor, having the same Fr\'echet functional    mean. The supports  of the marginal probability measures of the regressor  and response processes are assumed to be contained in a ball, whose radius ensures the injectivity of the exponential map. This map  has time--varying origin  at  the common   marginal  Fr\'echet functional   mean.
A weighted Fr\'echet mean approach is  adopted in the definition of the theoretical loss function. The regularized  Fr\'echet weights are computed,  in the time--varying tangent space from the log--mapped regressors.   Under these assumptions,  and some  Lipschitz regularity  sample path conditions, when a unique minimizer exists, the uniform weak--consistency of the empirical Fr\'echet curve predictor is obtained, under mean--square  ergodicity of the log--mapped regressor process in the first two moments.
  A simulated example in the sphere  illustrates the finite sample size  performance  of the proposed  Fr\'echet  predictor.  Predictions in time of the spherical coordinates of the magnetic field vector are obtained from the time--varying geocentric  latitude and longitude of the satellite NASA's MAGSAT spacecraft   in the real--data example analyzed.

\end{abstract}

\noindent \emph{Keywords}   Fr\'echet  functional  regression, Riemannian manifold, time correlated  manifold--valued  bivariate curve data,   weak--consistency.




\section{Introduction} \label{sec:1}

\setcounter{section}{1} \setcounter{equation}{0} 

Nonparametric regression techniques  have been widely  applied   in the last few decades to solve prediction problems from data lying on a Riemannian manifold (see \cite{Bhattacharya.12}).  Special attention has been paid to  kernel estimation, and local polynomial regression  exploiting the local character of the exponential map.
In the context of manifold--valued response and Euclidean predictors,  external local regression embeds  the manifold where the response lies onto  a higher dimensional Euclidean space, obtaining a local regression estimate in that space, and  back to the manifold via a homeomorphism relating the corresponding tangent spaces (see, e.g.,   \cite{Linthomaszhu17}).   In the same context of manifold--valued response and Euclidean regressors, in \cite{Zhu09},  an intrinsic regression model for the analysis of positive definite matrices is proposed for medical imaging processing.  It is well--known that positive definite matrices do not form a vector space. For these random elements  in Riemannian manifolds, a semiparametric regression model is proposed, considering a link function to map from the Euclidean space of covariates
to the Riemannian manifold of positive definite matrices
(see  also  approaches introduced in  \cite{Bhattacharya.12}, \cite{Marzio.14}, \cite{Khardani22},  \cite{KimPark13},    \cite{Patrangenaru.16}, \cite{Pelletier05}, \cite{Pelletier06},  mainly from a   nonparametric statistical framework).

Different versions of Functional Principal Component Analysis  on manifolds have been derived motivated by several fields of applications.
An intrinsic principal component analysis of Riemannian manifold-valued functional data, the so--called  Riemannian Functional Principal Component Analysis (RFPCA),  is obtained  in \cite{Dai.18}.
We also refer to the nonlinear manifold representation of  $L^{2}$ random functions themselves, lying in a low-dimensional but unknown manifold, or to the consideration of functional predictors lying on a smooth low-dimensional manifold (see \cite{Chen.12}, \cite{Dimeglio14},  \cite{Lazar17}, \cite{Lin.17}).

Among other fields of application, the presented Fr\'echet functional  regression   approach  in finite--dimensional compact Riemannian manifolds   is  of interest, for example,  in the analysis of bivariate plain trajectories data   (see \cite{Dai.18}  where RFPC is applied);  prediction from movement data in a pandemic (see \cite{Torres.21} for the Euclidean setting); prediction of
wind directions in environmental risk assessment  (see Section 9.2 in  \cite{Marzio.14});
 prediction of  uniform deviations   of comet orbits (see, e.g.,  \cite{Jupp03})  (see also \cite{Pigoli16} on  kriging  techniques for manifold--valued random fields).
   Another  field of application where the statistical analysis of  manifold data,  beyond the sphere, plays a crucial role is brain imaging data.
 For example, in the detection of brain  diseases or disorders,  the statistical analysis of covariates providing information about  the brain   structural and  morphological characteristics  of individuals  plays a crucial role in the construction of  RAVENS  maps reflecting    local volumetric    group differences  (see, e.g.,  \cite{Zhu13};  \cite{Zhu09}). In this paper, in the real--data example analyzed in Section \ref{rde}, we illustrate the performance of the proposed Fr\'echet functional regression approach,    in the prediction of  the time--varying  spherical coordinates of the magnetic field vector   from the geocentric  latitude and longitude of the satellite NASA's MAGSAT spacecraft, extending the purely spatial analysis  of \cite{Marzio.14}.

  To the best of our knowledge,  no  systematic approach or theoretical analysis has been developed, in the context of global intrinsic  regression
  in manifolds from time correlated bivariate curve data. This issue has been  addressed in the present paper adopting the theoretical framework of bivariate curve processes with values in a compact finite--dimensional  Riemannian manifold,    displaying suitable geometrical characteristics (see conditions (i)--(ii) in Section \ref{s31} below).
   Our proposal is based on the weighted Fr\'echet mean approach formulated in \cite{Petersen.19} for Euclidean regressors, and response evaluated in a metric space. Specifically, we extend this formulation  to the context of infinite--dimensional response and regressors, correlated in time, and evaluated
    in a finite--dimensional compact Riemannian manifold. One of the  main difficulties arising  in addressing this extended formulation is related to the uniform continuity of the theoretical  predictor, as well as of its empirical version in probability, both characteristics ensuring uniform weak--consistency. In particular,  mean--square first--order and second--order ergodicity of the log--mapped regressor in the time--varying tangent space is assumed in Section \ref{s31} for proving the weak--consistency of the empirical loss function.
     Fr\'echet weights are computed in   a regularized version of the ambient  Hilbert space $\mathbb{H}$ in the time--varying tangent space,   in terms of the semi--inner product of the Reproducing Kernel Hilbert Space (RKHS) of the log--mapped regressor process (see, e.g., \cite{Galeano.15};  \cite{RMAL19}; \cite{Kuelbs.70}). The involved regularizer matrix operator is defined from the pure point spectral properties of the matrix autocovariance operator of the log--mapped regressors in Section \ref{s32}.

    Among some additional geometrical and sample path  Lipschitz  regularity conditions, another challenging topic to be  addressed in our work   is the  existence and uniqueness of the Fr\'echet functional predictor, that requires some probabilistic structural restrictions on the underlying bivariate curve process. Specifically, strictly stationarity of this process is assumed. The inclusion on a closed ball of the support of the marginal infinite--dimensional probability measures of the response and regressor processes is also assumed, ensuring injectivity of the exponential maps, with time--varying origin at the  Fr\'echet functional mean, that minimizes the  quadratic mean geodesic dispersion of the curve values  of the regressor and response processes.
 An extensive literature exists on  the existence and uniqueness of the global $L^{p}$ center of mass (minimizer of the $L^{p}$--energy function)
of an arbitrary probability measure on a  manifold. This topic has been addressed, in particular,  for complete and connected  Riemannian manifolds, under some compactness assumptions on the support of the underlying probability measure   (see, e.g., Theorem 2.1 in \cite{Afsari11} and  \cite{Buss01}).  Specifically, if the support of the probability measure is contained in a ball, whose radius  is bounded  by a function of
 $p,$ the injectivity radius of the manifold,  and an upper bound on the manifold sectional curvatures, existence and uniqueness hold. Convexity   is a notion that  plays a crucial role in this problem (see, e.g., the overview presented in Section 1.1  in  \cite{Afsari11}).
  In \cite{Afsari13},   the convergence of a  constant step--size gradient descent algorithm is also  investigated  for solving this problem, in a general manifold context, analyzing the effect of the curvature and  topology of the manifold on the behavior of the algorithm.
  Beyond the above geometrical restrictions on the Riemannian manifold, we refer to the reader  to Theorem 1 in
     \cite{LeBarden17}.      As commented, in this paper, we pay attention to the $L^{2}$ center of mass of a probability measure  defining
   its Fr\'echet mean. The values of our   Fr\'echet  curve  predictor can be identified with the respective $L^{2}$ centers of mass of a family of probability measures indexed by the curve arguments of this predictor.

   The  performance of the proposed  Fr\'echet functional regression methodology is illustrated through simulations in a numerical example in the sphere in $\mathbb{R}^{3}.$ The implemented simulation algorithm is based on the generation of the curve values of the regressor  process by subordination,  via the poinwise application of  the inverse
    von Mises-Fisher distribution transform, to a family of correlated vector diffusion processes, whose sample paths are mapped into the unit ball of the three  dimensional Euclidean space.  The curve response process is then generated in the time--varying tangent space, by applying a bounded correlation operator to the log--mapped regressor curve values, and  adding a strong Gaussian white noise process, generated  from the Karhunen--Lo\'eve expansion. These simulations show a good  finite sample size performance of our  Fr\'echet functional regression predictor.

    In the real--data example analyzed in the context of world magnetic models, $5$--fold cross validation is implemented for testing the quality of our Fr\'echet functional   predictions  of the time--varying  spherical coordinates of the magnetic field vector,   from the geocentric  latitude and longitude of the satellite NASA's MAGSAT spacecraft. Data  have been obtained from NASA's National Space Science Data Center,  in the period 02/11/1979--06/05/1980. They have been  recorded every half second, and correspond to the first satellite NASA's MAGSAT spacecraft, which orbited the earth  every 88 minutes during seven months at around 400 km altitude.

  The outline of the paper is the following. Some preliminaries on Riemannian  manifolds  and Fr\'echet regression are given in Section \ref{bga}.
  Section \ref{rmv} introduces the assumed  conditions, and the proposed regression methodology from manifold--valued bivariate curve data  correlated in time. Under suitable conditions, the uniform weak--consistency of the empirical  Fr\'echet functional regression predictor is derived in Theorem \ref{th1} of Section \ref{wc}.  The proof of this result is given in the Appendix \ref{app1}. Simulations  in Section \ref{se} illustrate the finite sample  size performance of the proposed manifold  Fr\'echet functional  regression predictor.  The performance of this predictor is also analyzed from a real data example  in Section  \ref{rde}.   Appendix \ref{app2} complements the sample information and cross validation results displayed in Section  \ref{rde}.   Final comments are summarized in Section \ref{corl}.
 \section{Preliminaries}
\label{bga}

Let $\mathcal{M}$ be a smooth manifold with topological dimension $d$ in a
Euclidean space $\mathbb{R}^{d_{0}},$ $d\leq d_{0}$. Denote  by $\{ \mathcal{T}_{p}\mathcal{M},\ p\in \mathcal{M} \}$  the tangent spaces at the points of $\mathcal{M}.$ A
\emph{Riemannian metric} on $\mathcal{M}$ is a family of inner products $\mathcal{G}(p):\mathcal{T}_{p}\mathcal{M}\times \mathcal{T}_{p}\mathcal{M}\longrightarrow \mathbb{R}$ that smoothly varies over $p\in \mathcal{M}.$ Endowed with this Riemannian metric, $\left(\mathcal{M},\mathcal{G}\right)$  is a Riemannian  manifold. The metric on $\mathcal{M}$ induced by $\mathcal{G}$ is the geodesic distance $d_{\mathcal{M}}.$

 The \emph{exponential map} $\exp_{p}(v)$ is defined for $v\in \mathcal{T}_{p}\mathcal{M},$ and for each
$p\in \mathcal{M},$
  in terms of a locally length minimizing curve  $\gamma_{v}=\left\{\exp_{p}(tv), \ t\in [0,1]\right\},$ called geodesic, such that
for every  $v\in \mathcal{T}_{p}\mathcal{M},$   $\exp_{p}(v)=\gamma_{v}(1),$ where  $v\in \mathcal{T}_{p}\mathcal{M}.$ That is,  $\gamma_{v}$ is the unique geodesic with initial location $\gamma_{v}(0)=p,$ and velocity $\gamma_{v}^{\prime }(0)=v.$
The inverse of the exponential map is called the \emph{logarithm map}, and is denoted by $\log_{p},$ $p\in \mathcal{M}.$
The \textit{radius of injectivity} $\mbox{inj}_{p}$ at a point $p$ of the manifold
  is the radius of the largest ball about the origin of the tangent space $\mathcal{T}_{p}\mathcal{M}$  on which $\exp_{p}$ is a diffeomorphism, for each $p\in \mathcal{M}$. If $(\mathcal{M}, d_{\mathcal{M}})$ is a complete metric space, then $\exp_{p}$ is defined on the entire tangent space, and  $\exp_{p}$ is  a diffeomorphism at a neighborhood of
the origin of $\mathcal{T}_{p}\mathcal{M}$ (see, e.g.,  \cite{Chavel06}).

 We now briefly introduce the basic probabilistic  and function elements involved in the formulation of our approach (see Section \ref{rmv} on the conditions assumed).
Denote by $(\Lambda ,\mathcal{A},P)$ the basic probability space.
Consider the  space $\left(\mathcal{C}_{\mathcal{M}}(\mathcal{T}),d_{\mathcal{C}_{\mathcal{M}}(\mathcal{T})}\right)=\left\{x:\mathcal{T}\to \mathcal{M}:\ x\in\mathcal{C}(\mathcal{T})\right\},$
constituted by  $\mathcal{M}$-valued continuous functions on a compact interval $ \mathcal{T} $ with the supremum geodesic distance $$d_{\mathcal{C}_{\mathcal{M}}(\mathcal{T})}(x(\cdot),y(\cdot))=\sup_{t\in \mathcal{T}}d_{\mathcal{M}}\left( x(t),y(t)\right), \quad \forall x(t),y(t)\in \left(\mathcal{C}_{\mathcal{M}}(\mathcal{T}), d_{\mathcal{C}_{\mathcal{M}}(\mathcal{T})}\right).$$

Let $Z=\{Z_{s},\ s\in \mathbb{Z}\}$   be a family  of random elements  in  $\left(\mathcal{C}_{\mathcal{M}}(\mathcal{T}),d_{\mathcal{C}_{\mathcal{M}}(\mathcal{T})}\right) $   indexed by $\mathbb{Z}.$
    Specifically, $Z:\mathbb{Z}\times (\Lambda,\mathcal{A},\mathcal{P})\to \mathcal{C}_{\mathcal{M}}(\mathcal{T}),$   and \linebreak
 $\mathcal{P}\left(\omega\in \Lambda;\ Z_{s}(\cdot ,\omega)\in \left(\mathcal{C}_{\mathcal{M}}(\mathcal{T}),d_{\mathcal{C}_{\mathcal{M}}(\mathcal{T})}\right)\right)=1,$ for every $s\in  \mathbb{Z}.$
Here, $Z_{s}(t)$  denotes the pointwise value at $t\in \mathcal{T}$ of the   curve $Z_{s}$ in $\mathcal{M},$ for each  $s\in \mathbb{Z}.$ Consider the    ambient Hilbert space $\mathbb{H}$  of square--integrable vector functions  given by   \begin{equation}\mathbb{H}=\left\{h(\cdot)=(h_{1}(\cdot), \dots, h_{d_{0}}(\cdot))^{T}:\ \mathcal{T}\to \mathbb{R}^{d_{0}}:\ \int_{\mathcal{T}}h(t)^{T}h(t)dt<\infty\right\}.\label{ahs}
 \end{equation}
\noindent This space is  equipped with the inner product $\left\langle h,f\right\rangle_{\mathbb{H}}=\int_{\mathcal{T}}h(t)^{T}f(t)dt,$ and norm   $\|h\|_{\mathbb{H}}=\left[\left\langle h,h\right\rangle_{\mathbb{H}}\right]^{1/2},$ for any $h, f\in \mathbb{H}.$

For each  $s\in \mathbb{Z},$ define
the functional (curve)  Fr\'echet mean $\mu_{Z_{s},\mathcal{M}} $  as
\begin{eqnarray}
 \mu_{Z_{s},\mathcal{M}}(t)&=&\mbox{arg min}_{p\in \mathcal{M}}E\left([d_{\mathcal{M}}\left(Z_{s}(t),p\right)]^{2}\right)\nonumber\\
 &=&\mbox{arg min}_{p\in \mathcal{M}}\int_{\mathcal{C}_{\mathcal{M}}(\mathcal{T})} [d_{\mathcal{M}}\left(z_{s}(t),p\right)]^{2}dP_{Z_{s}(t)}(z_{s}(t)),\ t\in \mathcal{T},\nonumber\\
 \label{FM}
\end{eqnarray}
\noindent where  $dP_{Z_{s}(t)}$ denotes the  probability measure induced by $Z_{s}(t),$ the \linebreak $t$--projection  of the infinite--dimensional marginal probability measure $dP_{Z_{s}}$ of $Z_{s},$  for each $s\in \mathbb{Z}.$
Thus,   $\mu_{Z_{s},\mathcal{M}}$ is the curve in $\mathcal{M}$  providing the best  pointwise approximation of $Z_{s}$  in the mean quadratic geodesic distance sense. The sample path continuity of  $Z_{s}\in  \left(\mathcal{C}_{\mathcal{M}}(\mathcal{T}),d_{\mathcal{C}_{\mathcal{M}}(\mathcal{T})}\right)$ allows the following equivalent definition of the continuous function $\mu_{Z_{s},\mathcal{M}}(t):$

\begin{eqnarray}\mu_{Z_{s},\mathcal{M}}(\cdot )&=&\mbox{arg min}_{ z(\cdot )\in \mathcal{C}_{\mathcal{M}}(\mathcal{T})}
E\left(\int_{\mathcal{T}}[d_{\mathcal{M}}\left(Z_{s}(t), z(t)\right)]^{2}dt\right)\nonumber\\
&=&\mbox{arg min}_{ z(\cdot )\in \mathcal{C}_{\mathcal{M}}(\mathcal{T})}\int_{\mathcal{C}_{\mathcal{M}}(\mathcal{T})}\int_{\mathcal{T}}[d_{\mathcal{M}}\left(z_{s}(t), z(t)\right)]^{2}dt
dP_{Z_{s}}(z_{s}).\nonumber\\
\label{FFM}
\end{eqnarray}

   The response  $Y=\{ Y_{s},\ s\in \mathbb{Z}\}$ and regressor  $X=\{ X_{s},\ s\in \mathbb{Z}\}$  curve processes are introduced in  Section \ref{rmv},  under the above scenario, that is,  $X,$ and $Y$ satisfy the conditions of process $Z.$
\subsection{Fr\'echet regression in metric spaces with finite-dimensional Euclidean regressors}\label{Frims}

This section reviews some material from \cite{Petersen.19} to motivate the  regression approach presented.
Let $ (\Omega, d_{\Omega }) $ be a metric space.  Consider a random vector $(X, Y) \sim F ,$ where $dF(x,y)$ is the  probability measure induced by $(X, Y),$ with $ X $ and $ Y $ respectively taking values in $ \mathbb{R}^p $ and $ \Omega .$   The   marginal distributions of $ X $ and $ Y$ are denoted as  $ F_X $ and $ F_Y,$ respectively.  Assume that the mean vector  $ \mu_{X} = E(X),$ and the variance--covariance matrix  $ \Sigma  $ of the  regressors $X,$ as well as the conditional distributions $ F_{X \mid Y} $ and $ F_{Y \mid X} $ exist (see, e.g., Chapter V in \cite{Parthasarathy}).

\begin{definition}

The Fr\'echet regression predictor $m_{\oplus}(x)$  of $ Y $ given the observed value, $ X = x \in  \mathbb{R}^p,$ is defined as follows:
\begin{eqnarray}&&m_{\oplus}(x) = \arg\min_{\omega \in \Omega} M_{\oplus}( \omega, x)=
\arg\min_{\omega \in \Omega}
 E(d_{\Omega }^2(Y, \omega)\mid X = x)\nonumber\\
&&=\arg\min_{\omega \in \Omega}\int_{\Omega } d_{\Omega }^2(y, \omega) dF_{Y \mid X}(x, y).\label{fcm}
\end{eqnarray}

\end{definition}

Let us consider the case of real--valued response, i.e.,  $ \Omega = \mathbb{R},$ and   $d_{\Omega }=d_{E},$ where  $d_E(y_{1},y_{2})=|y_{1}-y_{2}|.$ It is well-known that  the least--squares linear global regression  predictor  $m_{L}(x)$ is computed in a parametric framework from the  minimizer
\begin{equation}
(\beta_0^*,\ \beta_1^* ) = \arg\min_{\beta_0 \in \mathbb{R},\  \beta_1 \in \mathbb{R}^p} \int \left[ \int y \ d F_{Y \mid X} (x,y) - (\beta_0 + \beta_1^T(x-\mu_{X}))\right]^2 \ d F_X (x),
\label{regrectamc}
\end{equation}
\noindent \noindent with  $E[X]=\mu_{X} \in \mathbb{R}^{p},$ where now $d F_{Y \mid X} (x,y) $ denotes the conditional probability measure on $\mathbb{R}$ induced by $Y$ given $X_{p\times 1}=x_{p\times 1}\in \mathbb{R}^{p},$ and, as before,
 $d F_X (x)$ denotes the marginal probability measure on $\mathbb{R}^{p}$ induced by $X.$  It is well--known that the scalar intercept $ \beta_0^* =E[Y],$ and the slope vector $ \beta_1^*=\Sigma^{-1}\sigma_{YX} ^{T}$, where $ [\sigma_{YX}] _{1\times p}= E[(Y-\mu_{Y})_{1\times 1}(X - \mu_{X})^{T}_{1\times p}],$ with $\mu_{Y}=E[Y],$
 and $\Sigma_{p\times p} =E\left[(X- \mu_{X})_{p\times 1}(X - \mu_{X})^T_{1\times p}\right]$ denoting  the matrix of variances and covariances of the regressor  vector  $X.$
Hence, for every $x\in \mathbb{R}^{p},$ one can write
\begin{equation}\label{Pet2.3}
m_{L}(x)=  \beta_0^* + (\beta_1^*)^T(x-\mu_{X}),
\end{equation}
\noindent  which can be equivalently expressed as
\begin{eqnarray}\label{Pet2.5}
m_{L}(x) &=& E(Y) + \sigma_{YX}\Sigma^{-1} (x - \mu_{X})\nonumber\\
&=&\int_{\mathbb{R}\times \mathbb{R}^{p}} y [1+(z - \mu_{X})^T \Sigma^{-1} (x - \mu_{X}) ]\ d F(z,y)\nonumber\\
&=&\int_{ \mathbb{R}\times \mathbb{R}^{p}} y s(z,x)\ d F(z,y).
\end{eqnarray}

\noindent Since the weight function $s(z,x) = 1 + (z-\mu_{X})^T \Sigma^{-1}(x-\mu_{X})$ satisfies $ \int_{\mathbb{R}^{p}\times \mathbb{R}}  s(z,x) d F(z,y) = 1,$  restricting our  attention to the unit ball of the RKHS of $X,$ one can consider the family of   bivariate probability measures
$\left\{P_{x}(dz,dy),\ x \in \mathbb{R}^{p}\right\},$ given by
  $P_{x}(dz,dy)=s(z,x) d F(z,y),$ $z\in \mathbb{R}^{p},$ $y\in \mathbb{R},$ for each $x\in \mathbb{R}^{p}.$ Consider the marginal  $P_{x}(dy)=\int_{\mathbb{R}^{p}}s(z,x) d F(z,y).$ Equation (\ref{Pet2.5}) can be symbolically rewritten as
\begin{eqnarray}
m_{L}(x)&=&\int_{\mathbb{R}} y P_{x}(dy) =\arg\min_{\omega  \in \mathbb{R}} \int_{\mathbb{R}} d_E^2 (y,\omega ) P_{x}(dy)\nonumber\\
&=&\arg\min_{\omega  \in \mathbb{R}} \int_{\mathbb{R}} d_E^2 (y,\omega )\int_{\mathbb{R}^{p}}s(z,x) d F(z,y)\nonumber\\
&=&\arg\min_{\omega  \in \mathbb{R}}E\left[ s(X,x) d_E^2 (Y, \omega  ) \right].
\label{eqrep}
\end{eqnarray}
  \noindent  The weighted Fr\'echet mean approach proposed in  \cite{Petersen.19}, under independent data,  consists of  replacing   the Euclidean distance $d_E $ by the distance $d_{\Omega }$ in an arbitrary metric space $(\Omega, d_{\Omega }),$ to cover the case where  $Y$  is  evaluated in such a  metric space$(\Omega, d_{\Omega }).$ That is,
  \begin{equation}
m_{L}(x)=\arg\min_{\omega  \in \Omega } E\left[ s(X,x) d_{\Omega }^{2} (Y, \omega  ) \right]
\label{Pet2.6inb}
\end{equation}
 \noindent (see \cite{Petersen.19}).  Our proposal is formulated under dependent curve data, and corresponds to the  case of  $(\Omega , d_{\Omega })= \left(\mathcal{C}_{ \mathcal{M}}(\mathcal{T}), d_{\mathcal{C}_{\mathcal{M}}(\mathcal{T})}\right),$ and $X\in \left(\mathcal{C}_{ \mathcal{M}}(\mathcal{T}), d_{\mathcal{C}_{\mathcal{M}}(\mathcal{T})}\right),$
 extending the above Euclidean regressor framework to the $\mathcal{M}$--valued  infinite--dimensional case.

 \section{Fr\'echet regression under dependent curve data in Riemannian manifolds}
\label{rmv}
In Section \ref{s31},  conditions (i)-(v)  are formulated   to provide  a suitable  geometrical and probabilistic scenario, allowing  the definition of our weighted Fr\'echet mean based  theoretical and empirical loss functions, and  ensuring the existence and uniqueness of   Fr\'echet curve means. Under  these conditions, the curve regressor process is mapped into  the time--varying tangent space.   Fr\'echet weights  are then computed in a regularized version of the RKHS of the log--mapped  regressor process, obtained by applying  a  smoother matrix operator satisfying    conditions (a)--(b)   in Section \ref{s32}. Finally, the theoretical and empirical loss functions are introduced  in Section \ref{s33}  under the  conditions established in Sections  \ref{s31}  and  \ref{s32}.

\subsection{Model assumptions}
\label{s31}
 The following conditions  of geometrical nature will be considered:
\begin{itemize}
\item[(i)] $\mathcal{M}$ is a $d$--dimensional compact  and connected Riemannian submanifold  of a Euclidean space  $\mathbb{R}^{d_{0}},$ $d\leq d_{0},$ with geodesic distance $d_{\mathcal{M}}$  induced by
the Euclidean metric.
\item[(ii)] The  sectional curvature of manifold  $\mathcal{M}$ is bounded, positive, and of smooth variation.
\end{itemize}
\begin{remark}
\label{rem1}
The exponential map is defined on the entire tangent space under (i) (see, e.g.,  \cite{Dai.18}). Under (ii),  the geodesic distance between two points in the manifold   is upper bounded by  the  Euclidean distance of their corresponding tangent vectors (see Assumption A2,  and Proposition 1 in \cite{Dai.18}).
\end{remark}

Let $Y=\{ Y_{s},\ s\in \mathbb{Z}\}$ and $X=\{ X_{s},\ s\in \mathbb{Z}\}$ be the response $Y$ and regressor $X$ curve processes  in the Riemannian  manifold $\mathcal{M}.$   Denote by $\mathcal{Y}_{\mathcal{C}_{\mathcal{M}}(\mathcal{T})}\subseteq \left(\mathcal{C}_{ \mathcal{M}}(\mathcal{T}), d_{\mathcal{C}_{\mathcal{M}}(\mathcal{T})}\right),$ and   $\mathcal{X}_{\mathcal{C}_{\mathcal{M}}(\mathcal{T})}\subseteq \left(\mathcal{C}_{ \mathcal{M}}(\mathcal{T}), d_{\mathcal{C}_{\mathcal{M}}(\mathcal{T})}\right)$ the respective  supports of their marginal probability measures (see conditions (iv)--(v) below).
The following  conditions are assumed on the bivariate curve process $(X,Y):$
\begin{itemize}
\item[(iii)]
For every time $s_{i}\in \mathbb{Z},$ the  random Lipschitz constants $L_{Y}(Y_{s_{i}})$ and  $L_{X}(X_{s_{i}})$ of $Y_{s_{i}}$ and $X_{s_{i}}$
    are almost surely (a.s.) finite. The Lipschitz constants
     $L(\mu_{Y_{s_{i}},\mathcal{M}})$ and  $L(\mu_{X_{s_{i}},\mathcal{M}})$  of the Fr\'echet means  $\mu_{Y_{s_{i}},\mathcal{M}}$ and $\mu_{X_{s_{i}},\mathcal{M}}$   are also finite.  Particularly,  assume that $E\left[\left(L_{X}(X_{s_{i}})\right)^{2}\right]<\infty,$  and \linebreak  $E\left[\left(L_{Y}(Y_{s_{i}})\right)^{2}\right]<\infty,$  for any $s_{i}\in \mathbb{Z}.$
     Note that,  for any curve $z(\cdot),$  $L(z)=\sup_{t\neq s}\frac{d_{\mathcal{M}}(z(t),z(s))}{\left|t-s\right|}.$
\item[(iv)] The $\mathcal{M}$--valued bivariate curve process $\{(Y_{s}, X_{s}),\ s\in \mathbb{Z}\}$ is strictly stationary.
Furthermore,   $\{ \log_{\mu_{X_{0},\mathcal{M}}(t)}\left(X_{s}(t)\right),\ s\in \mathbb{Z}\}$ is mean--square ergodic in the first moment  in the norm of $\mathbb{H},$ and in the second--order moments
in the norm of the space $\mathcal{S}(\mathbb{H})$ of Hilbert--Schmidt operators on $\mathbb{H}.$
\item[(v)]   We assume that
$X=\{ X_{s},\ s\in \mathbb{Z}\}$ and  $Y=\{ Y_{s},\ s\in \mathbb{Z}\}$
       have the same  Fr\'echet  functional mean. The supports of the marginal probability measures $dP_{X_{0}}(\cdot )$ and $dP_{Y_{0}}(\cdot )$ respectively induced by $X_{0}(\cdot )$ and $Y_{0}$ are included in the ball of the space $\left(\mathcal{C}_{ \mathcal{M}}(\mathcal{T}), d_{\mathcal{C}_{\mathcal{M}}(\mathcal{T})}\right),$ centered at the Fr\'echet functional  mean $\mu_{X_{0},\mathcal{M}}=\mu_{Y_{0},\mathcal{M}}$ with radius
$R=\inf_{t\in \mathcal{T}}\ \mbox{inj}_{\mu_{X_{0},\mathcal{M}}(t)}.$
  Here, $\mbox{inj}_{\mu_{X_{0},\mathcal{M}}(t)}$  denotes the injectivity radius of the exponential map whose  origin  is $\mu_{X_{0},\mathcal{M}}(t),$ for each  $t\in \mathcal{T}.$
\end{itemize}

  The global regression methodology, based on the weighted Fr\'echet mean approach, arising from the linear correlation between the response and regressor in the Euclidean setting (see Section \ref{Frims}) has sense, since the log--mapped versions of both curve processes,  $X$ and $Y,$ lie in the same time--varying tangent space, whose origin is the Fr\'echet functional  mean $\mu_{X_{0},\mathcal{M}}=\mu_{Y_{0},\mathcal{M}}$    under conditions (iv)--(v).

\begin{remark}
\label{remcondint}
Condition (v)  ensures that, for every $t\in \mathcal{T},$
 and $s_{i}\in \mathbb{Z},$  the geodesic connecting  $X_{s_{i}}(t)$ and  $\mu_{X_{0},\mathcal{M}}(t)$ is unique, ensuring that the tangent vectors do not switch directions under small perturbations of $\mu_{X_{0},\mathcal{M}}(t).$  The same assertion holds for  $Y_{s_{i}}(t)$  and $\mu_{Y_{0},\mathcal{M}}(t),$ $t\in \mathcal{T},$  $s_{i}\in \mathbb{Z}.$ This  condition is crucial in the implementation of the weighted Fr\'echet mean approach from time correlated bivariate curve data evaluated in $\mathcal{M}$ in a consistent way.

 In practice, a curve  clustering around the   intrinsic Fr\'echet functional mean  $\mu_{X_{0},\mathcal{M}}=\mu_{Y_{0},\mathcal{M}}$
 is observed under conditions (i)--(v). The same feature holds  around the empirical Fr\'echet functional mean $\widehat{\mu}_{X_{0},\mathcal{M}}$   (see, e.g., Figures \ref{fig:1}--\ref{fig:3} in the simulation study in Section \ref{se}). See also Figure \ref{fig:0rda} in Section  \ref{rde}, where the bivariate curve data  (left--hand side),  the  empirical  Fr\'echet   curve mean (center),  and responses  (right--hand side)  are  displayed. Note that this curve clustering  is observed in most of the real--data problems  cited in the Introduction. That is the case of   flight trajectory data set  (see, e.g., Section 5.2 of \cite{Dai.18},  where  Riemannian functional principal component analysis (RFPCA)  is applied).
\end{remark}

\subsection{Regularization of infinite--dimensional Fr\'echet weights}
\label{s32}
As commented, Fr\'echet weights  are computed from the log--mapped regressor process $\left\{\log_{\widehat{\mu}_{X_{0},\mathcal{M}}(\cdot)}(X_{s}(\cdot)), \ s\in \mathbb{Z}\right\}$  in the  time--varying tangent space.  Denote by $\mu (t) =E\left[\log_{\mu_{X_{0},\mathcal{M}}(t)}\left(X_{0}(t)
\right)\right],$  $t\in \mathcal{T},$ and by

\begin{eqnarray}&&\mathcal{R}_{X}=E\left[\left(\log_{\mu_{X_{0},\mathcal{M}}(\cdot)}\left(X_{s}(\cdot)
\right)-\mu (\cdot)\right)\otimes \left(\log_{\mu_{X_{0},\mathcal{M}}(\cdot)}\left(X_{s}(\cdot)
\right)-\mu (\cdot) \right)^{T}\right]\nonumber\\ &&\hspace*{1cm}=E\left[ \left(\log_{\mu_{X_{0},\mathcal{M}}(\cdot)}\left(X_{0}(\cdot)
\right)-\mu (\cdot)\right)\otimes \left(\log_{\mu_{X_{0},\mathcal{M}}(\cdot)}\left(X_{0}(\cdot)
\right)-\mu (\cdot)\right)^{T} \right],\nonumber\\
\label{cots}
\end{eqnarray}
  \noindent the matrix autocovariance operator of $\left\{\log_{\widehat{\mu}_{X_{0},\mathcal{M}}(\cdot)}(X_{s}(\cdot)), \ s\in \mathbb{Z}\right\}.$ Its inverse, defining the  semi--inner  product $\left\langle f,g\right\rangle_{\widetilde{H}}=\left\langle \mathcal{R}_{X}^{-1}(f),g\right\rangle_{\mathbb{H}},$ of the  RKHS  $\widetilde{H}=\mathcal{R}_{X}^{1/2}(\mathbb{H})$ of $\left\{\log_{\widehat{\mu}_{X_{0},\mathcal{M}}(\cdot)}(X_{s}(\cdot)), \ s\in \mathbb{Z}\right\}$  is not bounded in the ambient Hilbert space $\mathbb{H}.$  Fr\'echet weights are then computed in a regularized version $\boldsymbol{\mathcal{K}}(\mathbb{H})\subset \widetilde{H}$ of
  $\mathbb{H}$ (see condition (b) below),  obtained from  smoother matrix operator  $\boldsymbol{\mathcal{K}}$  given by
   \begin{eqnarray}&\sqrt{\boldsymbol{\mathcal{K}}}(\psi)&=\left[\begin{array}{cccc}
 \sqrt{\mathcal{K}}& 0 &\dots & 0\\
0 &\sqrt{\mathcal{K}}& \dots & 0\\
\dots &\dots &\dots & \sqrt{\mathcal{K}}\\
 \end{array}\right]\left[\begin{array}{c} \psi_{1}\\ \vdots\\ \psi_{d_{0}}\\  \end{array}\right]\nonumber\\
 &&=
 \left[\sqrt{\mathcal{K}}(\psi_{1}), \dots , \sqrt{\mathcal{K}}(\psi_{d_{0}})\right]^{T},\quad \forall \psi =(\psi_{1},\dots,\psi_{d_{0}} )^{T} \in \mathbb{H}.
 \nonumber
 \label{ro} \end{eqnarray}
 \noindent  Thus, $\sqrt{\boldsymbol{\mathcal{K}}}$ is a diagonal operator with     constant
 functional entries equal to $\sqrt{\mathcal{K}},$ where $\mathcal{K}$ is a  trace self--adjoint integral operator on $\mathbb{H},$  satisfying
 $$\mathcal{K}(\phi_{n})(t)=\int_{\mathcal{T}}k(t,s)\phi_{n}(s)ds=\gamma_{n}\phi_{n}(t),\ t\in \mathcal{T},\ \gamma_{n}>0,\ n\geq 1,$$
 \noindent with $\sum_{n=1}^{\infty}\gamma_{n}=1.$   The norm induced by $\boldsymbol{\mathcal{K}}$ can be expressed  in terms of the  eigenvalue sequence
 $\left\{ \gamma_{n},\ n\in \mathbb{N}_{0}\right\}$ as
 \begin{eqnarray}
 &&\|\psi\|_{\boldsymbol{\mathcal{K}}}^{2}=\left\langle \psi, \psi\right\rangle_{\boldsymbol{\mathcal{K}}}=\sum_{i=1}^{d_{0}}\sum_{k\geq 1}\gamma_{k}|\phi_{k}(\psi_{i})|^{2}=\|\psi\|_{\widetilde{H}_{W}}^{2},\label{eqhw}
\end{eqnarray}
\noindent where we have denoted by $\widetilde{H}_{W}$ the separable Hilbert space of vector  functions with finite $\|\cdot\|_{\boldsymbol{\mathcal{K}}}$--norm induced by the inner product  $\left\langle \psi, \varphi\right\rangle_{\widetilde{H}_{W}},$ given by, for
  $\psi =(\psi_{1},\dots,\psi_{d_{0}})^{T},  \varphi=(\varphi_{1},\dots,\varphi_{d_{0}})^{T}\in \mathbb{H}=H^{d_{0}},$
\begin{eqnarray}
&\left\langle \psi, \varphi\right\rangle_{\widetilde{H}_{W}}&=\left\langle \sqrt{\boldsymbol{\mathcal{K}}}(\psi),\sqrt{\boldsymbol{\mathcal{K}}}(\varphi )\right\rangle_{\mathbb{H}}  =
 \int_{\mathcal{T}} \left[\sqrt{\boldsymbol{\mathcal{K}}}(\psi)(t)\right]^{T}\sqrt{\boldsymbol{\mathcal{K}}}(\varphi )(t)dt
\nonumber\\
&&
=\sum_{i=1}^{d_{0}}\int_{\mathcal{T}}\int_{\mathcal{T}} k(t,s)\psi_{i}(t) \varphi_{i}(s) dtds\nonumber\\
&&=
\sum_{i=1}^{d_{0}}\sum_{k\geq 1}\gamma_{k}\phi_{k}(\psi_{i})
\phi_{k}(\varphi_{i})=\sum_{i=1}^{d_{0}}\sum_{k\geq 1}\gamma_{k}\left\langle \phi_{k},\psi_{i}\right\rangle_{H}
\left\langle \phi_{k},\varphi_{i}\right\rangle_{H}.\nonumber\\
\label{hwip}
\end{eqnarray}

 Assume that $\left\{ \gamma_{n},\ n\geq 1\right\}$  satisfy
 \begin{itemize}
\item[(a)]For every $f\in \mathbb{H},$
 $\|f\|_{\mathbb{H}}\geq \|f\|_{\widetilde{H}_{W}}.$ Hence,   $\mathbb{H}\subset \widetilde{H}_{W}$ in a continuous way.
 \item[(b)] For every $\psi \in \sqrt{\boldsymbol{\mathcal{K}}}(\mathbb{H}),$
$\| \sqrt{\boldsymbol{\mathcal{K}}}^{-1}(\psi) \|_{\mathbb{H}}\geq \|\psi \|_{\widetilde{H}}.$ Therefore,
 $\left(\sqrt{\mathcal{K}}(\mathbb{H}), \left\langle\cdot,\cdot\right\rangle_{\boldsymbol{\mathcal{K}}^{-1}} \right)$ $\subset \left(\widetilde{H},
 \left\langle\cdot,\cdot\right\rangle_{\widetilde{H}}\right)$ in a continuous way.
\end{itemize}

Hence, the following continuous inclusions (embeddings, denoted as $\hookrightarrow $)  hold:
\begin{equation}\left(\sqrt{\mathcal{K}}(\mathbb{H}), \left\langle\cdot,\cdot\right\rangle_{\boldsymbol{\mathcal{K}}^{-1}} \right)\hookrightarrow \left(\widetilde{H},
 \left\langle\cdot,\cdot\right\rangle_{\widetilde{H}}\right)\hookrightarrow \left(\mathbb{H},\left\langle \cdot, \cdot \right\rangle_{\mathbb{H}}\right) \hookrightarrow\left(\widetilde{H}_{W},
 \left\langle\cdot,\cdot\right\rangle_{\widetilde{H}_{W}}\right),\label{emb}
 \end{equation}
\noindent and the  following regularized version of Fr\'echet weights is computed:\begin{eqnarray}&&
s\left(X_{0}(\cdot ),x(\cdot)\right)=\nonumber\\ &&\hspace*{-0.5cm}=
\left[1+ \left\langle \log_{\mu_{X_{0},\mathcal{M}}(\cdot)}\left(x(\cdot)
\right)-\mu (\cdot) , \sqrt{\mathcal{K}}\mathcal{R}_{X}^{-1}\left( \sqrt{\mathcal{K}}\left(\log_{\mu_{X_{0},\mathcal{M}}(\cdot)}\left(X_{0}(\cdot)
\right)-\mu (\cdot) \right)\right)
\right\rangle_{\mathbb{H}}\right]
\nonumber\\
&&=\left[
1+ \left\langle \sqrt{\mathcal{K}}\left(\log_{\mu_{X_{0},\mathcal{M}}(\cdot)}\left(x(\cdot )\right)-\mu (\cdot)\right), \sqrt{\mathcal{K}}\left(\log_{\mu_{X_{0},\mathcal{M}}(\cdot)}\left(X_{0}(\cdot )\right)-\mu (\cdot)\right)\right\rangle_{\widetilde{H}}
\right].\nonumber\\
\label{tlfbb2}\end{eqnarray}

\begin{remark}
 Equation (\ref{tlfbb2}) restricts  the support of the proposed theoretical Fr\'echet  functional  predictor $\widehat{Y}(x(\cdot)),$ in equation (\ref{eqfp}), and its empirical version  $\widehat{Y}_{n}(x(\cdot)),$ in equation (\ref{predte}),
  to $\exp_{\mu_{X_{0},\mathcal{M}}(\cdot)}\left(\sqrt{\boldsymbol{\mathcal{K}}}(\mathbb{H})\right)\subset \exp_{\mu_{X_{0},\mathcal{M}}(\cdot)}\left(\widetilde{H}\right)\subset \mathcal{X}_{\mathcal{C}_{\mathcal{M}}(\mathcal{T})}
 \subseteq \left(\mathcal{C}_{ \mathcal{M}}(\mathcal{T}), d_{\mathcal{C}_{\mathcal{M}}(\mathcal{T})}\right),$  allowing their computation in  a continuous way. This property  is applied  in
the proof of  Theorem \ref{th1} to derive   uniform weak--consistency.
\end{remark}

\subsection{Formulation of the theoretical and empirical loss functions}
\label{s33}

Under conditions (i)--(v), let us consider, for every $z(\cdot)\in \mathcal{Y}_{\mathcal{C}_{\mathcal{M}}(\mathcal{T})}$ and   $x(\cdot )\in \mathcal{X}_{\mathcal{C}_{\mathcal{M}}(\mathcal{T})},$ the loss function $M\left(z(\cdot), x(\cdot )\right)$ given by:

\begin{eqnarray}&&
M\left(z(\cdot), x(\cdot )\right)=
E\left[s\left(X_{0}(\cdot ),x(\cdot)\right)\int_{\mathcal{T}}\left[d_{\mathcal{M}}\left( Y_{0}(t),z(t)\right)\right]^{2}dt\right],
\label{tlfbb}\end{eqnarray}
\noindent  where   the regularized Fr\'echet weights $s\left(X_{0}(\cdot ),x(\cdot)\right)$ have been introduced in  equation (\ref{tlfbb2}).
The proposed Fr\'echet predictor is  obtained as the solution to the following minimization problem:
\begin{eqnarray}\widehat{Y}(x(\cdot))&=&\mbox{arg min}_{z(\cdot)\in \mathcal{Y}_{\mathcal{C}_{\mathcal{M}}(\mathcal{T})}}M\left(z(\cdot), x(\cdot )\right),\ x(\cdot )\in \mathcal{X}_{\mathcal{C}_{\mathcal{M}}(\mathcal{T})}.\label{eqfp}
\end{eqnarray}

Let
$((X_{1}(\cdot), Y_{1}(\cdot)) ,\dots, (X_{n}(\cdot), Y_{n}(\cdot)))$ be
a bivariate  functional sample of size  $n$  of correlated curves in time of the $\mathcal{M}$--valued response and regressor curve  processes.
 For each $x(\cdot )\in \mathcal{X}_{\mathcal{C}_{\mathcal{M}}(\mathcal{T})},$
and $z(\cdot)\in \mathcal{Y}_{\mathcal{C}_{\mathcal{M}}(\mathcal{T})},$
the empirical version of  (\ref{tlfbb}), based on  $((X_{1}(\cdot), Y_{1}(\cdot) ),\dots, (X_{n}(\cdot), Y_{n}(\cdot))),$  is defined as
\begin{eqnarray}&&
\widehat{M}_{n}(z(\cdot), x(\cdot ))=\frac{1}{n}\sum_{i=1}^{n}s_{n}\left(X_{i}(\cdot ),x(\cdot)\right)\int_{\mathcal{T}}\left[d_{\mathcal{M}}\left( Y_{i}(t),z(t)\right)\right]^{2}dt, \nonumber\\
&&\hspace*{3cm} \forall  z(\cdot)\in \mathcal{Y}_{\mathcal{C}_{\mathcal{M}}(\mathcal{T})},\quad x(\cdot )\in \mathcal{X}_{\mathcal{C}_{\mathcal{M}}(\mathcal{T})},
\label{elfbb}
\end{eqnarray}
    \noindent where
    \begin{eqnarray}&& s_{n}\left(X_{i}(\cdot ),x(\cdot)\right)=\left[1+\left\langle \sqrt{\mathcal{K}}\left(\boldsymbol{\gamma}_{x(\cdot ),\overline{X}_{n}}\right), \widehat{\mathcal{R}}_{X}^{-1}\left(\sqrt{\mathcal{K}}\left(\boldsymbol{\gamma}_{X_{i}(\cdot ),\overline{X}_{n}}\right)\right)\right\rangle_{\mathbb{H}}\right],\nonumber\\
    &&\quad   \quad \quad \forall  x(\cdot )\in \mathcal{X}_{\mathcal{C}_{\mathcal{M}}(\mathcal{T})}, \quad  i=1,\dots,n,\nonumber\\ \label{elfbb2}
\end{eqnarray}
\noindent with
\begin{eqnarray}
\widehat{\mathcal{R}}_{X}&=&\frac{1}{n}\sum_{i=1}^{n}
    \left[\log_{\widehat{\mu}_{X_{0},\mathcal{M}}(\cdot)}\left(X_{i}(\cdot)\right)-\overline{X}_{n}(\cdot )\right]\otimes \left[\log_{\widehat{\mu}_{X_{0},\mathcal{M}}(\cdot)}\left(X_{i}(\cdot)\right)-\overline{X}_{n}(\cdot )\right]^{T}\nonumber\\
\boldsymbol{\gamma}_{x(\cdot ),\overline{X}_{n}}&=&
\log_{\widehat{\mu}_{X_{0},\mathcal{M}}(\cdot)}\left(x(\cdot )\right)-\overline{X}_{n}(\cdot ),\quad x(\cdot )\in \mathcal{X}_{\mathcal{C}_{\mathcal{M}}(\mathcal{T})},\nonumber\\
\boldsymbol{\gamma}_{X_{i}(\cdot ),\overline{X}_{n}}&=&\log_{\widehat{\mu}_{X_{0},\mathcal{M}}(\cdot)}\left(X_{i}(\cdot )\right)-\overline{X}_{n}(\cdot ),\quad i=1,\dots,n,\label{elfbb3}\\
\overline{X}_{n}(\cdot)&=&\overline{\log_{\widehat{\mu}_{X_{0},\mathcal{M}}(\cdot)}(X(\cdot))}
=\frac{1}{n}\sum_{i=1}^{n}\log_{\widehat{\mu}_{X_{0},\mathcal{M}}(\cdot)}(X_{i}(\cdot))\nonumber\\
 \widehat{\mu}_{X_{0},\mathcal{M}}(\cdot)&=&\mbox{arg min}_{\widetilde{x}(\cdot )\in \mathcal{X}_{
\mathcal{C}_{\mathcal{M}}(\mathcal{T})}} \frac{1}{n}   \sum_{i=1}^{n}\int_{\mathcal{T}}\left[d_{\mathcal{M}}\left(X_{s_{i}}(t),\widetilde{x}(t)\right)\right]^{2}dt.
\label{efmvf}
\end{eqnarray}
\noindent Weak--consistency, in the supremum geodesic distance,  of $\widehat{\mu}_{X_{0},\mathcal{M}}(\cdot)$ has been derived in Proposition 2 in \cite{Dai.18} under curve independent data evaluated in a manifold. Under conditions  (i)--(v), this result holds under weak--dependent curve data evaluated in $\mathcal{M}.$

The empirical predictor is  the  solution of the following minimization problem:
\begin{eqnarray}
\widehat{Y}_{n}(x(\cdot))&=&
    \mbox{arg min}_{z(\cdot)\in  \mathcal{Y}_{\mathcal{C}_{\mathcal{M}}(\mathcal{T})}}
        \widehat{M}_{n}\left(z(\cdot), x(\cdot )\right),\quad x(\cdot )\in \mathcal{X}_{\mathcal{C}_{\mathcal{M}}(\mathcal{T})}.\label{predte}\end{eqnarray}

Conditions (i)--(v) ensure the existence and uniqueness of the theoretical $\widehat{Y},$ in (\ref{eqfp}),  and empirical $\widehat{Y}_{n},$ in (\ref{predte}),   Fr\'echet functional  regression predictors (see, e.g., Theorem 2.1 in \cite{Afsari11}).

 \section{Weak-consistency of the  predictor}

\label{wc}

The main result of this paper, Theorem \ref{th1} below, provides the uniform weak--consistency of the proposed empirical Fr\'echet functional predictor. To derive  the proof of this result, the following additional conditions are assumed:

\begin{itemize}
\item[A.1] For each $x(\cdot)
\in \mathcal{X}_{\mathcal{C}_{\mathcal{M}}(\mathcal{T})},$
and,  for every  $\varepsilon >0,$
    \begin{eqnarray}
&&\hspace*{-1cm} \inf_{\sup_{t\in \mathcal{T}}d_{\mathcal{M}}\left(z(t),\widehat{Y}(x(\cdot))(t)\right)>\varepsilon}M\left(z(\cdot),x(\cdot )\right) > M\left(\widehat{Y}(x(\cdot)), x(\cdot )\right)\nonumber\\
&&\hspace*{-1.5cm}P\left( \inf_{\sup_{t\in \mathcal{T}}d_{\mathcal{M}}\left(z(t),\widehat{Y}_{n}(x(\cdot))(t)\right)
     >\varepsilon}\widehat{M}_{n}\left(z(\cdot), x(\cdot )\right)-\widehat{M}_{n}\left(\widehat{Y}_{n}(x(\cdot)),    x(\cdot )\right)\geq \zeta (\varepsilon )\right)\to 1,\nonumber\\
\label{p1}
\end{eqnarray}
\noindent  as $n\to \infty,$ for certain $\zeta (\varepsilon )>0.$

 \item[B.1] Assume that,  for every  $\varepsilon>0,$
     \begin{eqnarray}
    &&\hspace*{-1cm} \inf_{x(\cdot)\in \mathcal{X}_{\mathcal{C}_{\mathcal{M}}(\mathcal{T})}; \left\|\log_{\mu_{X_{0},\mathcal{M}}(\cdot)}(x(\cdot))\right\|_{\mathbb{H}}\leq B}\
    \inf_{\sup_{t\in \mathcal{T}}d_{\mathcal{M}}\left(z(t),\widehat{Y}(x(\cdot))(t)\right)
     >\varepsilon} M\left(z(\cdot), x(\cdot )\right)\nonumber\\ &&
       \hspace*{7cm} -M\left(\widehat{Y}(x(\cdot)), x(\cdot )\right)>0, \nonumber\\
           &&\hspace*{-1.5cm}  P\left(\inf_{x(\cdot)\in  \mathcal{X}_{\mathcal{C}_{\mathcal{M}}(\mathcal{T})};\ \left\|\log_{\widehat{\mu}_{X_{0},\mathcal{M}}(\cdot)}(x(\cdot))\right\|_{\mathbb{H}}\leq B}\
    \inf_{\sup_{t\in \mathcal{T}}d_{\mathcal{M}}\left(z(t),\widehat{Y}_{n}(x(\cdot))(t)\right)
     >\varepsilon}\widehat{M}_{n}\left(z(\cdot), x(\cdot )\right)\right.\nonumber\\ &&\hspace*{0.5cm}\left.
    \hspace*{5cm} -\widehat{M}_{n}\left(\widehat{Y}_{n}(x(\cdot)),    x(\cdot )\right)\geq \zeta (\varepsilon )\right)\to 1,\nonumber\\
     \label{prob}
     \end{eqnarray}
          \noindent as $n\to \infty,$ for certain $\zeta (\varepsilon )>0.$
          \begin{remark}
          \label{intervb1}
          Under condition (v),  the radius $B$ in equation (\ref{prob}) can be  the radius of the closed  ball in the ambient Hilbert space $\mathbb{H}$ containing  the log--mapped support of the probability measure $dP_{X_{0}(\cdot )}$ induced by $X_{0}(\cdot).$

          \end{remark}

\item[C.1]  The   functional moments of the log--mapped bivariate curve  process \linebreak $\left\{\left(\log_{\widehat{\mu}_{X_{0},\mathcal{M}}(\cdot)}(X_{s}(\cdot)),\log_{\widehat{\mu}_{X_{0},\mathcal{M}}(\cdot)}(Y_{s}(\cdot)) \right),\ s\in \mathbb{Z}\right\}$ satisfy the
following \linebreak summability conditions:
 \begin{eqnarray}&&\sum_{u\in \mathbb{Z}}E\left[\left\|\log_{\mu_{X_{0},\mathcal{M}}(\cdot)}(Y_{0}(\cdot))-
    \log_{\mu_{X_{0},\mathcal{M}}(\cdot)}(z(\cdot))\right\|_{\mathbb{H}}^{2}\right.\nonumber\\
    &&\hspace*{2cm}\left.\times
    \left\|\log_{\mu_{X_{0},\mathcal{M}}(\cdot)}(Y_{u}(\cdot))-
    \log_{\mu_{X_{0},\mathcal{M}}(\cdot)}(z(\cdot))\right\|_{\mathbb{H}}^{2}\right]<\infty.
    \nonumber\\
    &&\sum_{u\in \mathbb{Z}}E\left[\left\|\log_{\mu_{X_{0},\mathcal{M}}(\cdot)}(Y_{0}(\cdot))-
    \log_{\mu_{X_{0},\mathcal{M}}(\cdot)}(z(\cdot))\right\|_{\mathbb{H}}^{2}\right.\nonumber\\
    &&\hspace*{3cm}\left.\times  \left\|\log_{\mu_{X_{0},\mathcal{M}}(\cdot)}(Y_{u}(\cdot))-
    \log_{\mu_{X_{0},\mathcal{M}}(\cdot)}(z(\cdot))\right\|_{\mathbb{H}}^{2}\right.\nonumber\\
    &&\hspace*{3cm}\left.\times  \left\|\log_{\mu_{X_{0},\mathcal{M}}(\cdot)}(X_{u}(\cdot))-
     \mu\right\|_{\mathbb{H}}\right]<\infty,\nonumber\\
&&\sum_{u\in \mathbb{Z}}E\left[\left\|\log_{\mu_{X_{0},\mathcal{M}}(\cdot)}(Y_{0}(\cdot))-
    \log_{\mu_{X_{0},\mathcal{M}}(\cdot)}(z(\cdot))\right\|_{\mathbb{H}}^{2}\right.\nonumber\\
    &&\hspace*{3cm}\left.\times
    \left\|\log_{\mu_{X_{0},\mathcal{M}}(\cdot)}(Y_{u}(\cdot))-
    \log_{\mu_{X_{0},\mathcal{M}}(\cdot)}(z(\cdot))\right\|_{\mathbb{H}}^{2}\right.\nonumber\\
    &&\hspace*{3cm}\left.\times
    \left\|\log_{\mu_{X_{0},\mathcal{M}}(\cdot)}(X_{0}(\cdot))-
     \mu\right\|_{\mathbb{H}}\right.\nonumber\\
    &&\hspace*{3cm}\left.\times
    \left\|\log_{\mu_{X_{0},\mathcal{M}}(\cdot)}(X_{u}(\cdot))-
     \mu\right\|_{\mathbb{H}}\right]<\infty.\nonumber
    \end{eqnarray}

\end{itemize}
\begin{remark}
\label{r5}
In the proof of Theorem \ref{th1} below, condition C.1 allows proving  pointwise mean--square consistency  in the first argument of the empirical loss function, when the second--order  moments of the log--mapped curve  regressor process are totally specified. In the misspecified case, conditions (i)-(v)  also allow  proving consistency of the empirical loss function from this result, given the compactness of $\mathcal{T}$ and $\mathcal{M}.$ Condition A.1  is required to apply Corollary
3.2.3 in \cite{Vaart.1996}  to obtain  the weak--consistency of the empirical curve predictor in the supremum geodesic distance.
Condition B.1 leads to the  uniform weak--consistency of the empirical Fr\'echet functional  predictor,
applying  uniform equicontinuity   satisfied by the theoretical and empirical loss functions in the second argument under the assumed conditions in Sections \ref{s31}  and \ref{s32}.
\end{remark}

The following result provides the conditions ensuring the uniform weak--consistency, in the supremum geodesic distance, of the global empirical Fr\'echet curve   predictor $\widehat{Y}_{n}(x(\cdot)).$
 \begin{theorem}
 \label{th1} Under conditions  (i)--({\em v}) in Section \ref{s31}, and conditions  (a)--(b)  in Section  \ref{s32}, if  assumptions   A.1 and C.1 hold, then,
 \begin{equation}\sup_{t\in \mathcal{T}} d_{\mathcal{M}}\left(\widehat{Y}(x(\cdot))(t),\widehat{Y}_{n}(x(\cdot))(t)\right)=o_{P}(1),\label{c1w}
 \end{equation}
 \noindent
 \noindent  for each  $x(\cdot )\in \mathcal{X}_{\mathcal{C}_{\mathcal{M}}(\mathcal{T})}.$

  Additionally, under assumption    B.1 , \begin{equation}
  \sup_{\left\|\log_{\mu_{X_{0},\mathcal{M}}(\cdot)}(x(\cdot))\right\|_{\mathbb{H}}\leq B}\sup_{t\in \mathcal{T}} d_{\mathcal{M}}\left(\widehat{Y}(x(\cdot))(t),\widehat{Y}_{n}(x(\cdot))(t)\right)=o_{P}(1).\label{c1w2}
 \end{equation}
 \end{theorem}

 As commented in Remark \ref{r5}, the proof of Theorem \ref{th1} (see Appendix \ref{app1}) follows from conditions A.1, B.1 and C.1, applying
 Corollary
3.2.3 in \cite{Vaart.1996}, under the scenario defined by conditions   (i)--(v) in Section \ref{s31}, and conditions  (a)--(b)  in Section  \ref{s32}.
 Specifically, under condition C.1,
 \begin{equation}\widetilde{M}_{n}\left(z(\cdot), x(\cdot ) \right)-M\left(z(\cdot), x(\cdot )\right)=o_{P}(1),\ n\to \infty,\label{pr1}
 \end{equation}
  \noindent providing  the pointwise weak--consistency  in the first argument of the empirical loss function, when the second--order  moments of the log--mapped curve  regressor process are totally specified.
  Here,  \begin{eqnarray} &&\widetilde{M}_{n}\left(z(\cdot), x(\cdot ) \right)=\frac{1}{n}\sum_{i=1}^{n} \int_{\mathcal{T}}\left[d_{\mathcal{M}}\left( Y_{i }(t),z(t)\right)\right]^{2}dt
\nonumber\\
&&\hspace*{-1.5cm}\times \left[1+\left\langle \sqrt{\mathcal{K}}\left(\log_{\mu_{X_{0},\mathcal{M}}(\cdot)}\left(x(\cdot )\right)-\mu(\cdot )\right), \mathcal{R}_{X}^{-1}\left(\sqrt{\mathcal{K}}\left(\log_{\mu_{X_{0},\mathcal{M}}(\cdot)}\left(X_{i}(\cdot )\right)-\mu(\cdot )\right)\right)\right\rangle_{\mathbb{H}}\right].\nonumber\end{eqnarray}
 \noindent Conditions (i)--(v)  lead to \begin{equation}\widehat{M}_{n}\left(z(\cdot), x(\cdot )\right)-\widetilde{M}_{n}\left(z(\cdot), x(\cdot )\right)=o_{P}(1),\ n\to \infty.\label{cons}
 \end{equation}

Under conditions (i)--(v) and (a)--(b), keeping in mind  the compactness of  $\mathcal{T}$  and $\mathcal{M},$ and Remark  \ref{intervb1},
\begin{equation}\sup_{d_{\mathcal{C}_{\mathcal{M}}(\mathcal{T})}(z_{1}(\cdot),z_{2}(\cdot))\leq \delta } \left|\widehat{M}_{n}(z_{1}(\cdot),x(\cdot))-\widehat{M}_{n}(z_{2}(\cdot),x(\cdot))\right|=\mathcal{O}_{P}(\delta ).\label{consb}
 \end{equation}

From equations (\ref{pr1})--(\ref{consb}), under condition A.1,   Corollary
3.2.3 in \cite{Vaart.1996}  leads to the  weak--consistency of the  Fr\'echet  functional  predictor in the supremum geodesic distance.

 Condition B.1 allows proving uniform weak--consistency,  in the supremum geodesic distance,  of the empirical Fr\'echet  functional  predictor
 applying Theorem 1.5.4 in \cite{Vaart.1996}.  Specifically,  the proved equicontinuity of
  the theoretical loss function in the second argument, uniformly with respect to  the first argument,  leads to the uniform continuity of the theoretical predictor by  applying first part of B.1. In addition, the  equicontinuity in probability of the empirical loss function  in the second argument, uniformly with respect to  the first argument,  leads to the uniform continuity in probability of the empirical predictor by applying the second part of B.2. Applying triangle inequality, Theorem 1.5.4 in \cite{Vaart.1996}
  leads to the desired result on uniform weak--consistency of the  empirical Fr\'echet curve  predictor.

\begin{corollary}
\label{th1cor}   If there exist  positive constants $\mathcal{U}(Y)$ and $\mathcal{U}(X)$ such that
\begin{eqnarray}
&& P\left(\sup_{z(\cdot)\in \mathcal{Y}_{\mathcal{C}_{\mathcal{M}}(\mathcal{T})}}\left\|\log_{\mu_{X_{0},\mathcal{M}}(\cdot)}(Y_{0}(\cdot))-
    \log_{\mu_{X_{0},\mathcal{M}}(\cdot)}(z(\cdot))\right\|_{\mathbb{H}}\leq \mathcal{U}(Y)\right)=1,\nonumber\\
    &&P\left(\left\|\log_{\mu_{X_{0},\mathcal{M}}(\cdot)}(X_{0}(\cdot))-\mu \right\|_{\mathbb{H}}\leq \mathcal{U}(X)\right)=1,\nonumber\\
\label{efvcor}
\end{eqnarray}
\noindent which is the case under conditions (i)--(v), then,  condition C.1  in Theorem \ref{th1}  can be replaced by
\begin{eqnarray}
&& \sum_{u\in \mathbb{Z}}E\left[\left\|\log_{\mu_{X_{0},\mathcal{M}}(\cdot)}(Y_{u}(\cdot))-
    \log_{\mu_{X_{0},\mathcal{M}}(\cdot)}(z(\cdot))\right\|_{\mathbb{H}}^{2}\right]<\infty,\ \forall z(\cdot)\in \mathcal{Y}_{\mathcal{C}_{\mathcal{M}}(\mathcal{T})}.\nonumber\\
    \label{eqc1rep}
\end{eqnarray}
\end{corollary}
\begin{remark}
Note that, considering $\log_{\mu_{X_{0},\mathcal{M}}(\cdot)}(z(\cdot))=E\left[\log_{\mu_{X_{0},\mathcal{M}}(\cdot)}(Y_{0}(\cdot))\right]$
in equation (\ref{eqc1rep}),
  the log-mapped
curve response process \linebreak $\left\{\log_{\mu_{X_{0},\mathcal{M}}(\cdot)}(Y_{s}(\cdot)),\ s\in \mathbb{Z}\right\}$   displays Short Range Dependence (SRD)  (see, e.g., \cite{Panaretos13}). Thus, Theorem \ref{th1} provides uniform weak--consistency of the empirical Fr\'echet curve predictor $\widehat{Y}_{n}(\cdot)$ in the supremum geodesic distance under weak--dependent curve  data evaluated in $\mathcal{M}.$
\end{remark}
\section{Numerical examples}
\label{se}

 The finite functional sample size  performance of $\widehat{Y}_{n}(x(\cdot )),$ $x(\cdot )\in
 \mathcal{X}_{\mathcal{C}_{\mathcal{M}}(\mathcal{T})},$
   is numerically illustrated  in this section. We restrict our attention to  the sphere  $\mathbb{S}_{2}$ in $\mathbb{R}^{3},$
   and generate, at $1000$ temporal nodes, a bivariate curve sample of size $n= 100$  of time correlated random spherical curves.   We compute the  pointwise   quadratic  geodesic distances between the original values of the response, and their Fr\'echet curve predictions. These values are summarized in terms of the  empirical mean   (see the left--hand side of Figure \ref{fig:7}), and the corresponding  histogram  (see the left--hand side of Figure \ref{fig:8}).   The  empirical temporal means of the observed values of the  quadratic geodesic curve errors  at each sampled time are   also computed  (see the right--hand side of Figure \ref{fig:7}).  The corresponding histogram is displayed at the  right--hand side of Figure \ref{fig:8}.

    We have implemented, in MatLab language, a simulation algorithm based on  vector diffusion process subordination  by applying the inverse
    von Mises-Fisher distribution transform (see Algorithm 3 in the Supplementary Material of \cite{Terdik22}).
    We restrict our attention to the family of vector  diffusion processes  with  linear drift obeying the following stochastic differential equation:

\begin{equation}dX_{t}=\mu(t)X_{t}dt+\sigma (t,X_{t})dW_{t},\label{sde}
\end{equation}\noindent where $X_{t}$ defines the vector process modeling  the states of the system, $\mu(t)X_{t} $ represents the linear drift, and   $\sigma (t,X_{t})= D(t,X_{t})V(t)$ defines  the diffusion coefficient. Thus, the   coefficient  $\sigma (t,X_{t})$ is computed from  a diagonal matrix $D$, where each element along the main diagonal is the corresponding element of the state vector $X_{t}.$  The process $W$ is vectorial  Brownian motion with correlated components.  This model has been generated with \emph{sde}$(\cdot)$ MatLab function (see Simulation of Multidimensional Market Models in Financial Toolbox of MatLab).

\begin{figure}[!h]
\begin{center}
\includegraphics[trim= 100 280 100 230, clip, width=0.32\textwidth]{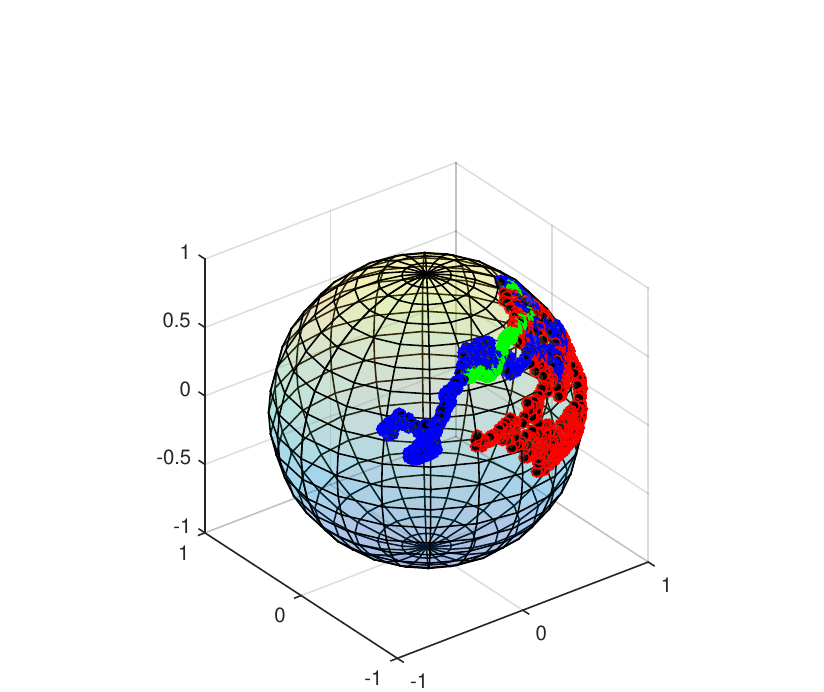}
\includegraphics[trim= 100 280 100 230, clip, width=0.32\textwidth]{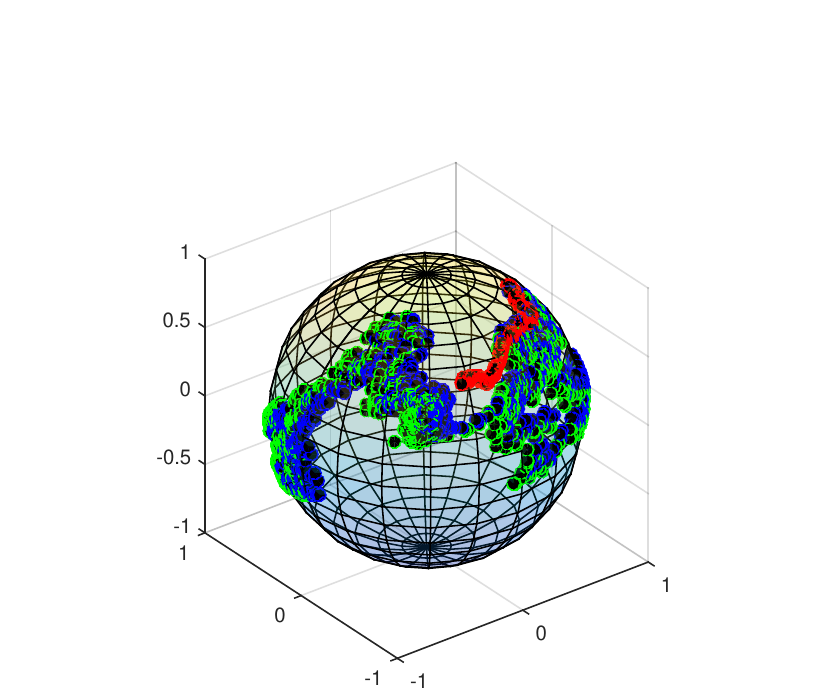}
\includegraphics[trim= 100 280 100 230, clip, width=0.32\textwidth]{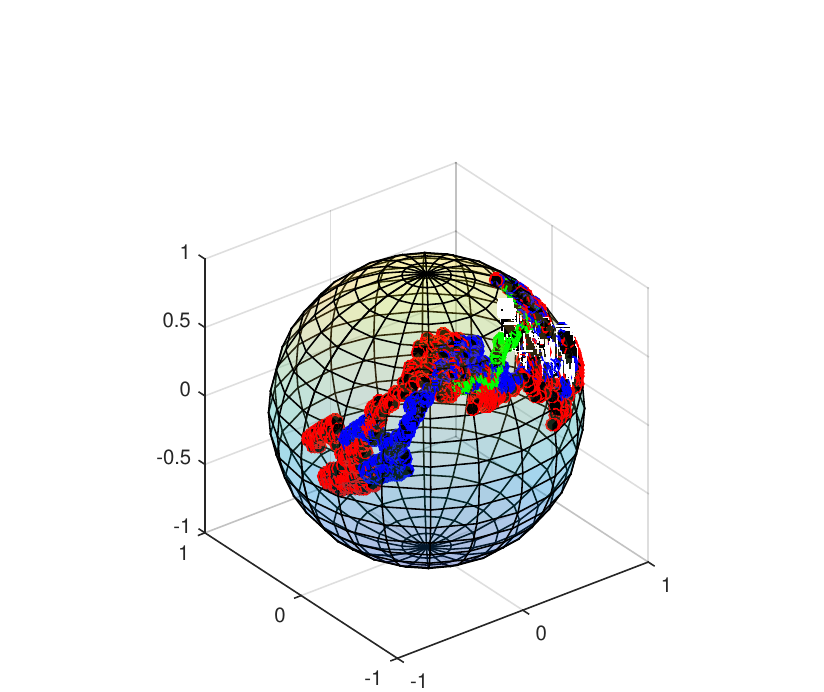}
\end{center}
\caption{ Spherical curve values of the generated  regressor process at times $s_{i}=10, 20, 30, 40, 50, 60,70 ,80, 90,100.$ Three sampling times  are represented on the left and center plots, and four on the right plot.}\label{fig:1}
\end{figure}

  A sample of $n$ strictly stationary and ergodic  vector diffusion processes correlated in time is generated. This sample is mapped into the unit ball of
the Banach space   $\mathcal{C}_{\mathcal{X}}(\mathcal{T})$ of continuous vector functions with compact support contained in the interval $\mathcal{T}$ taking values   in  the bounded set $\mathcal{X}\subset \mathbb{R}^{3}.$
  The inverse von Mises-Fisher transform is then pointwise applied to obtain  a sample of the regressor curve process $ X_{s_{i}}(t)=\left(X_{s_{i}}^{(1)}(t), X_{s_{i}}^{(2)}(t),  X_{s_{i}}^{(3)}(t)\right),$ $t\in \mathcal{T},$ $i=1,\dots,n$  (see  Figure \ref{fig:1}).
For each $s_{i},$ $i=1,\dots,n,$ and for every $t\in \mathcal{T},$ we compute  the logarithm map of the generated sample values of the spherical curve regressor process as follows: For each $t\in \mathcal{T},$
  \begin{eqnarray}&&\log_{\mu_{X_{0},\mathcal{M}}(t)}(X_{s_{i}}(t))=\frac{u(t,i)}{\|u(t,i)\|}d_{\mathbb{S}_{d}}
  \left(\mu_{X_{0},\mathcal{M}}(t),X_{s_{i}}(t)\right),\nonumber\\ &&u(t,i)=X_{s_{i}}(t)-([\mu_{X_{0},\mathcal{M}}(t)]^{T}X_{s_{i}}(t))\mu_{X_{0},\mathcal{M}}(t),\label{lms}\end{eqnarray}
\noindent where $\mu_{X_{0},\mathcal{M}}(\cdot)$ denotes, as before,  the Fr\'echet  functional   mean of $X_{0}.$
Finally, for  $i=1,\dots,n,$ we generate the response curve process as
\begin{equation} Y_{s_{i}}(t)= \exp_{\mu_{X_{0},\mathcal{M}}(t)}\left(\mathbf{\Gamma }\left(\log_{\mu_{X_{0},\mathcal{M}}(\cdot )}\left(X_{s_{i}}(\cdot )\right)\right)(t)+\varepsilon_{s_{i}}(t)\right),\quad  t\in \mathcal{T},\label{genresp}
\end{equation}
\noindent where $\mathbf{\Gamma }:\mathbb{H}\to \mathbb{H}$ is a bounded linear operator, with supremum norm less than one, and $\varepsilon= \left\{\varepsilon_{i}(\cdot),\ i\in \mathbb{Z}\right\}$  defines an $\mathbb{H}$--valued Gaussian  strong--white noise, uncorrelated with the log--mapped regressors.  Process $\varepsilon $  has been generated in terms of the Karhunen--Lo\'eve expansion (see \cite{Dai.18}).

Given the  geometrical characteristics  of the sphere, conditions (i)--(ii) hold. Condition (iii) is ensured by the sample path regularity properties
displayed by the generated vector diffusion processes, and by the regularity properties of the eigenfunctions, involved  in  the  Karhunen--Lo\'eve expansion of the strong Gaussian white noise process $\varepsilon$ in the time--varying tangent space in (\ref{genresp}). The  assumed mean--square ergodicity in condition (iv) follows  from the mixing conditions satisfied by the vector diffusion processes. The generation of the response process in the time--varying tangent space
in (\ref{genresp}) leads to its strictly stationarity from the strictly  stationarity of the vector diffusion processes, and of $\varepsilon $ in this example.  Our choice of the  concentration parameter involved in the implementation of  the inverse von Mises-Fisher transform, and of the correlation structure of the generated vector--valued diffusion process ensures that condition (v)   holds for the  corresponding regressor process $X.$
Condition (v)  is also satisfied by the generated response curve process, given  the dispersion of the   curve values of $\varepsilon$  in the time--varying tangent space is controlled by the truncated trace of its autocovariance operator, dominating the dispersion of the resulting exponential mapped curve values of the response process. Note that this trace  is computed from the eigenvalues defining the variance of the random coefficients in  the Karhunen--Lo\'eve expansion of $\varepsilon.$ Furthermore, the selected   eigenfunctions of the Dirichlet negative Laplace operator on a spatial interval are bounded and smooth.

\begin{figure}[!h]
\begin{center}
\includegraphics[trim= 100 280 100 220, clip, width=0.32\textwidth]{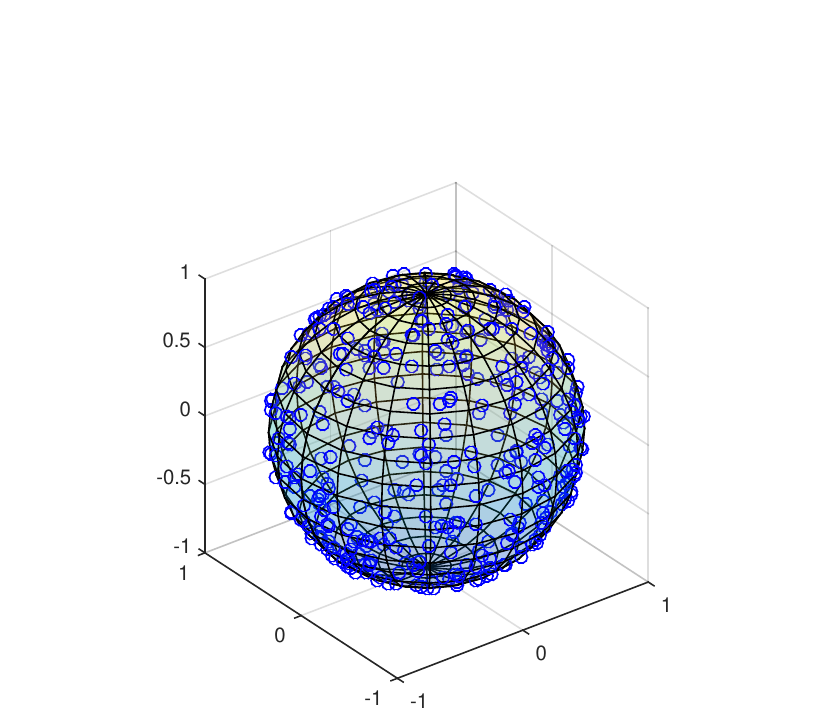}
\includegraphics[trim= 100 280 100 220, clip, width=0.32\textwidth]{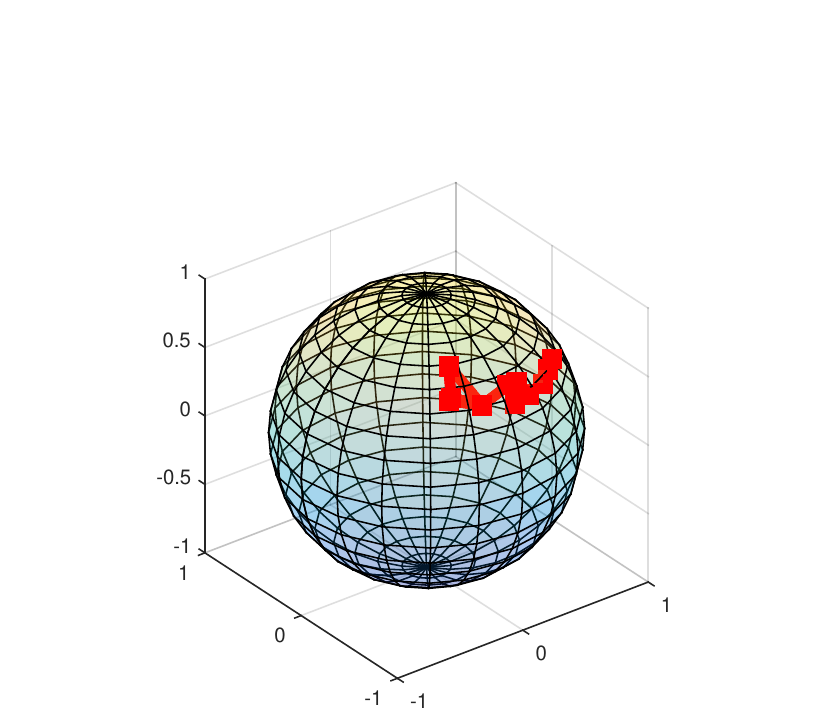}
\includegraphics[trim= 100 280 100 220, clip, width=0.32\textwidth]{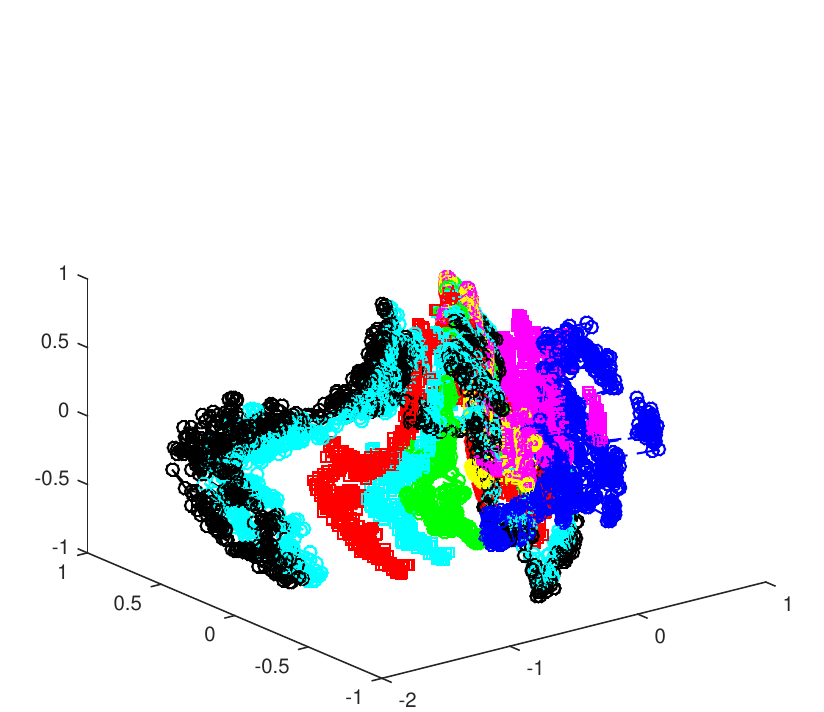}
 \end{center}
\caption{Uniform spherical grid with $400$  nodes (left--hand side),  the empirical curve Fr\'echet mean  $\mu_{\widehat{X}_{0},\mathcal{M}}(\cdot)$ (center), and the  response in the time--varying tangent space at times $s_{i}=10,20,30,40,50,60,70,80,90,100$ (right--hand side).}
\label{fig:4}
\end{figure}

\begin{figure}[!h]
\begin{center}
\includegraphics[trim= 100 280 100 220, clip, width=0.32\textwidth]{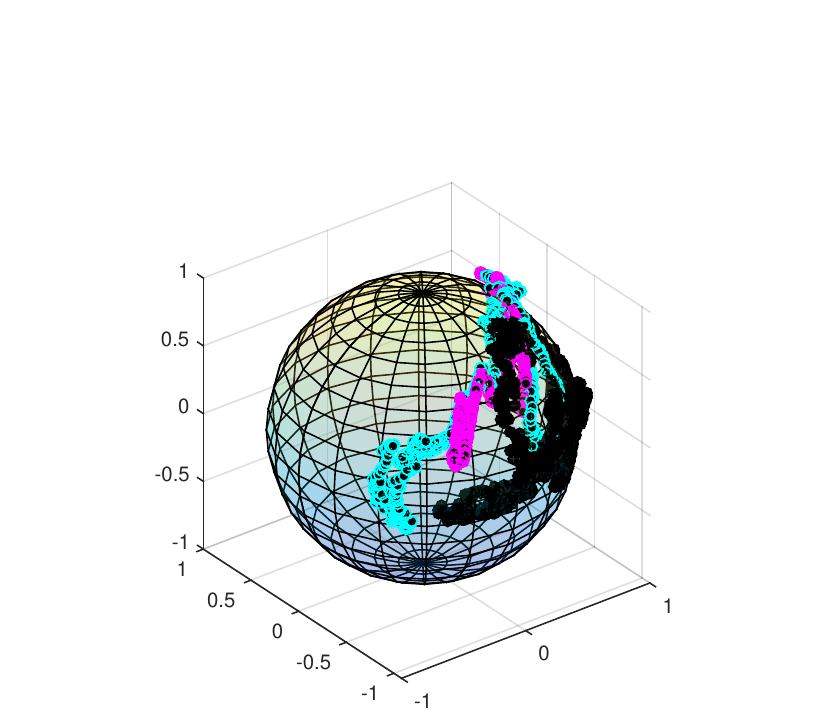}
\includegraphics[trim= 100 280 100 220, clip, width=0.32\textwidth]{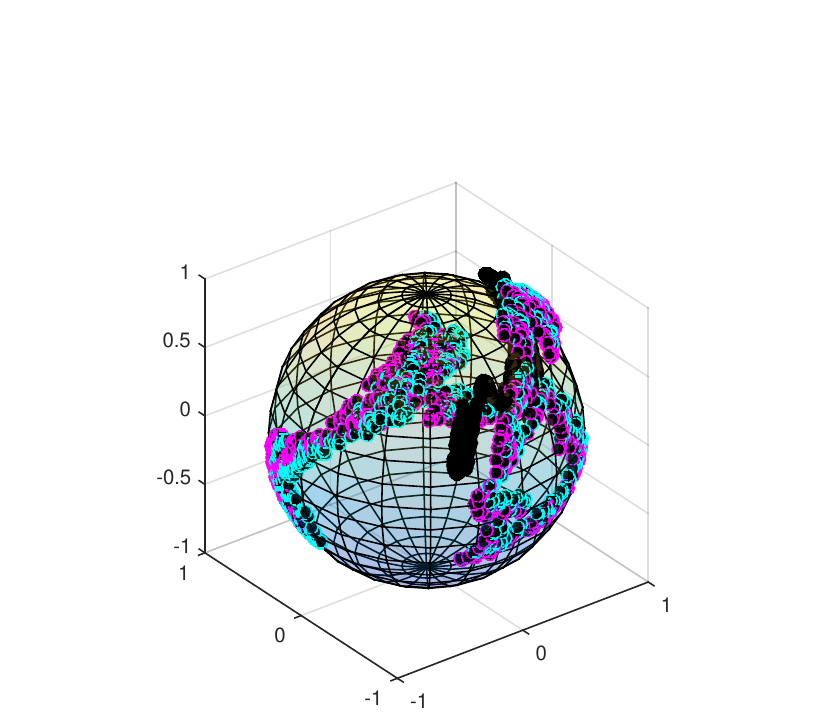}
\includegraphics[trim= 100 280 100 220, clip, width=0.32\textwidth]{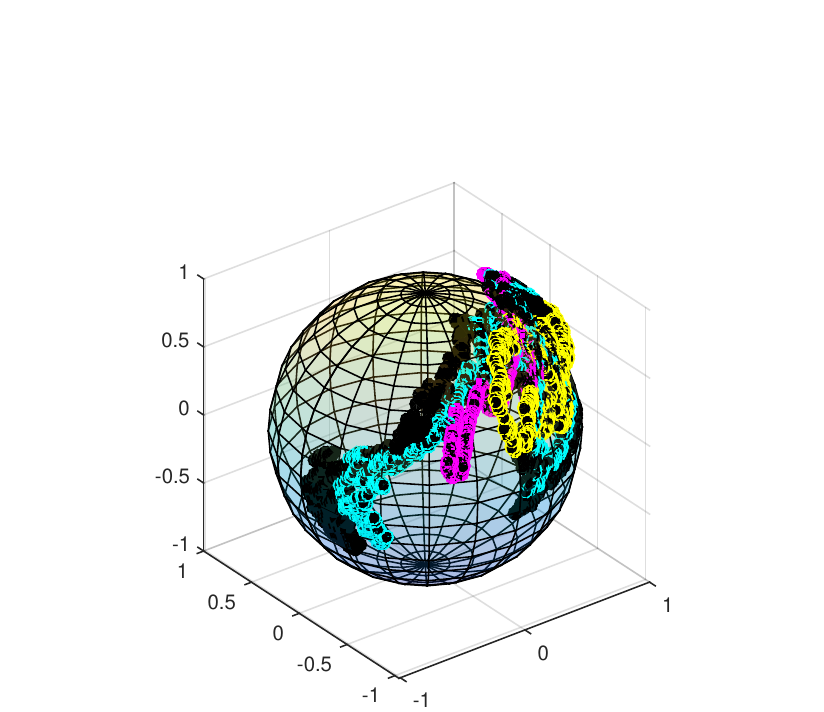}
\end{center}
\caption{Spherical curve values of the generated spherical curve response process at times $s_{i}=10, 20, 30, 40, 50, 60,70 ,80, 90,100.$}
\label{fig:2}
\end{figure}

\begin{figure}[!h]
\begin{center}
\includegraphics[trim= 100 280 100 220, clip, width=0.32\textwidth]{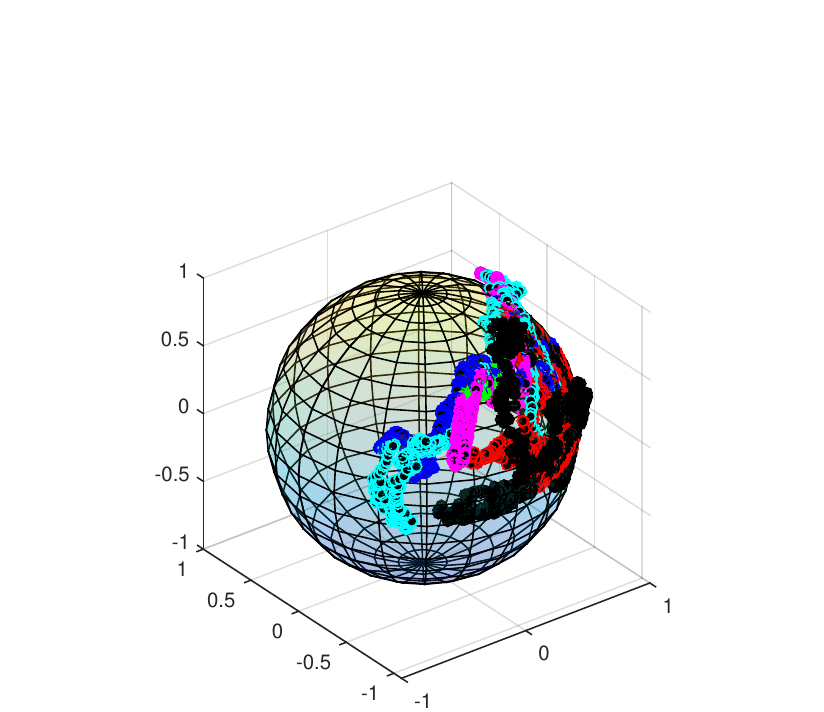}
\includegraphics[trim= 100 280 100 220, clip, width=0.32\textwidth]{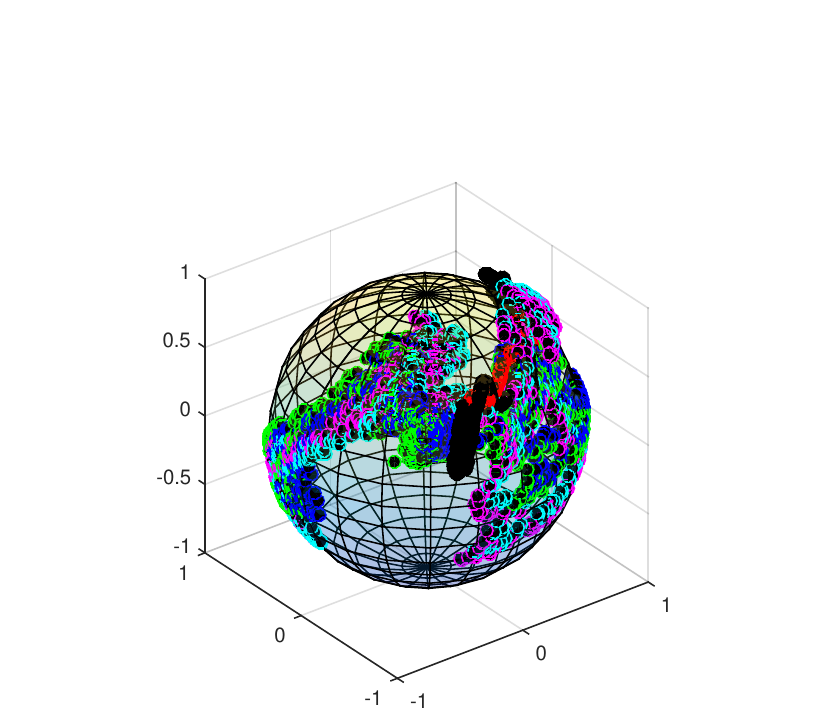}
\includegraphics[trim= 100 280 100 220, clip, width=0.32\textwidth]{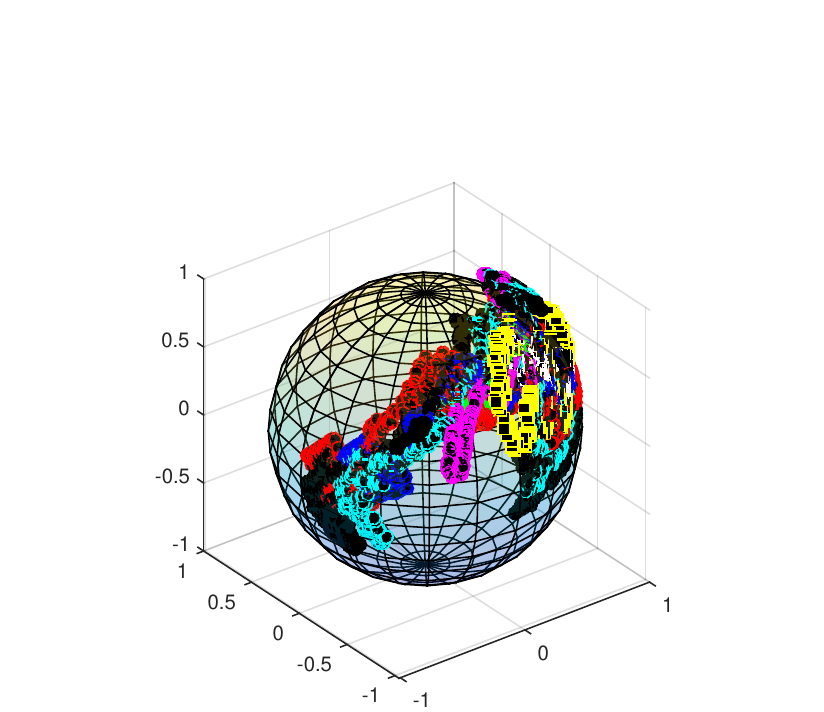}
\end{center}
\caption{Joint representation of spherical curve values of  the generated spherical curve response and regressor processes at times $s_{i}=10, 20, 30, 40, 50, 60,70 ,80, 90,100.$ Here, red, green, blue, and white colors are used for regressor spherical curve values, while black, magenta, cyan, and yellow colors are used  for response  spherical curve values.}
\label{fig:3}
\end{figure}

\begin{figure}[!h]
\begin{center}
 \includegraphics[height=10cm, width=12cm]{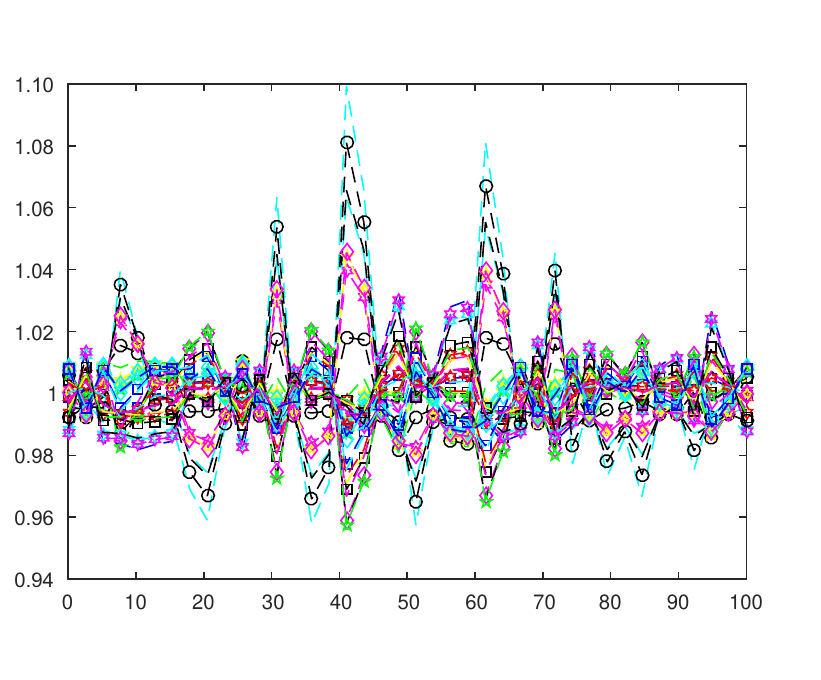}
\end{center}
\caption{The empirical Fr\'echet  weights are plotted. The dashed lines  displayed correspond to  their evaluation from generated log--mapped sample curve regressor  values at  $s_{i},$ $i=1,\dots,100,$  and  for  $40$ log--mapped $\mathcal{M}$--curve arguments of the predictor. }
\label{fig:5}
\end{figure}

\begin{figure}[!h]
\begin{center}
\includegraphics[trim= 100 280 100 220, clip, width=0.32\textwidth]{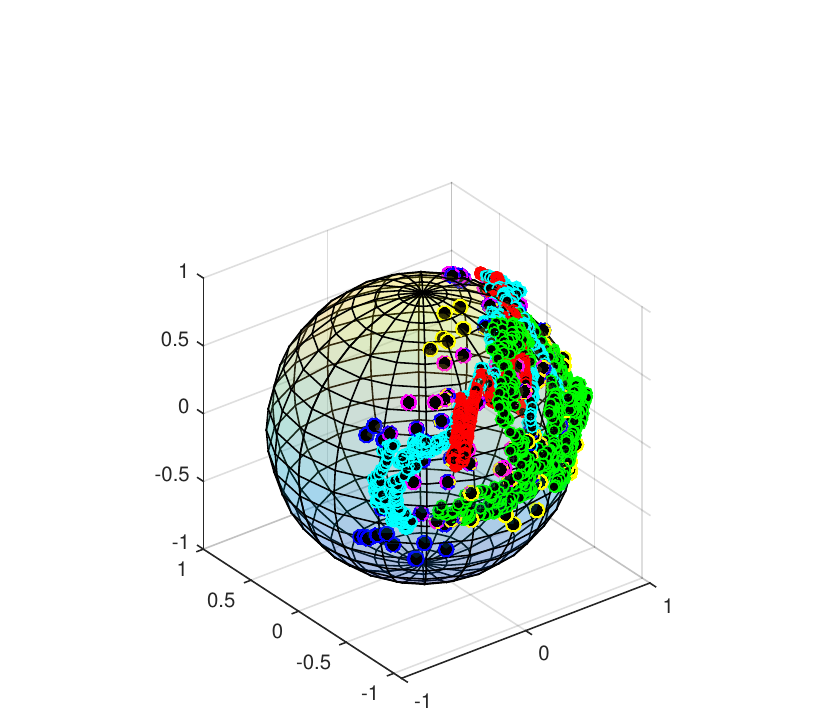}
\includegraphics[trim= 100 280 100 220, clip, width=0.32\textwidth]{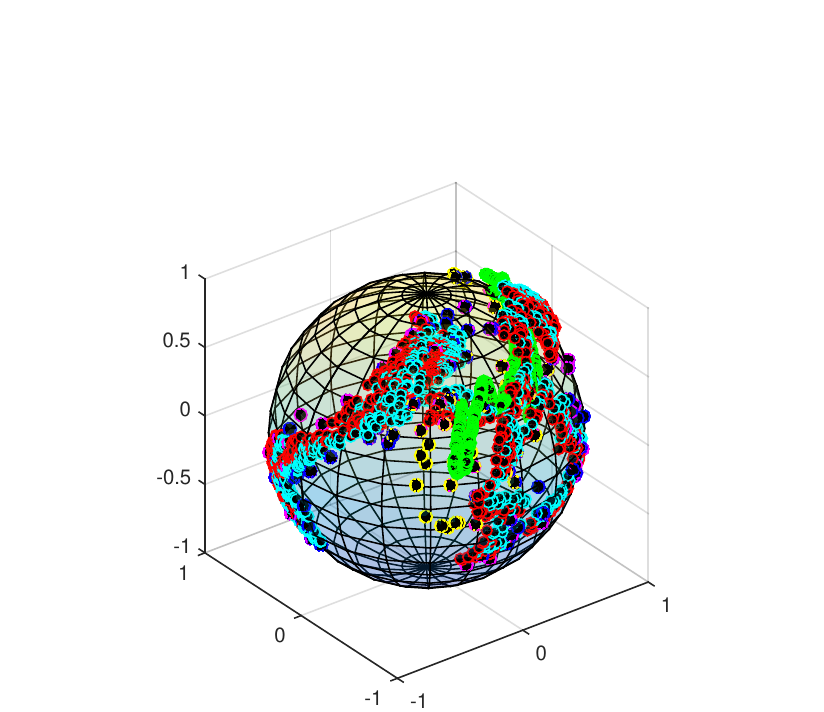}
\includegraphics[trim= 100 280 100 220, clip, width=0.32\textwidth]{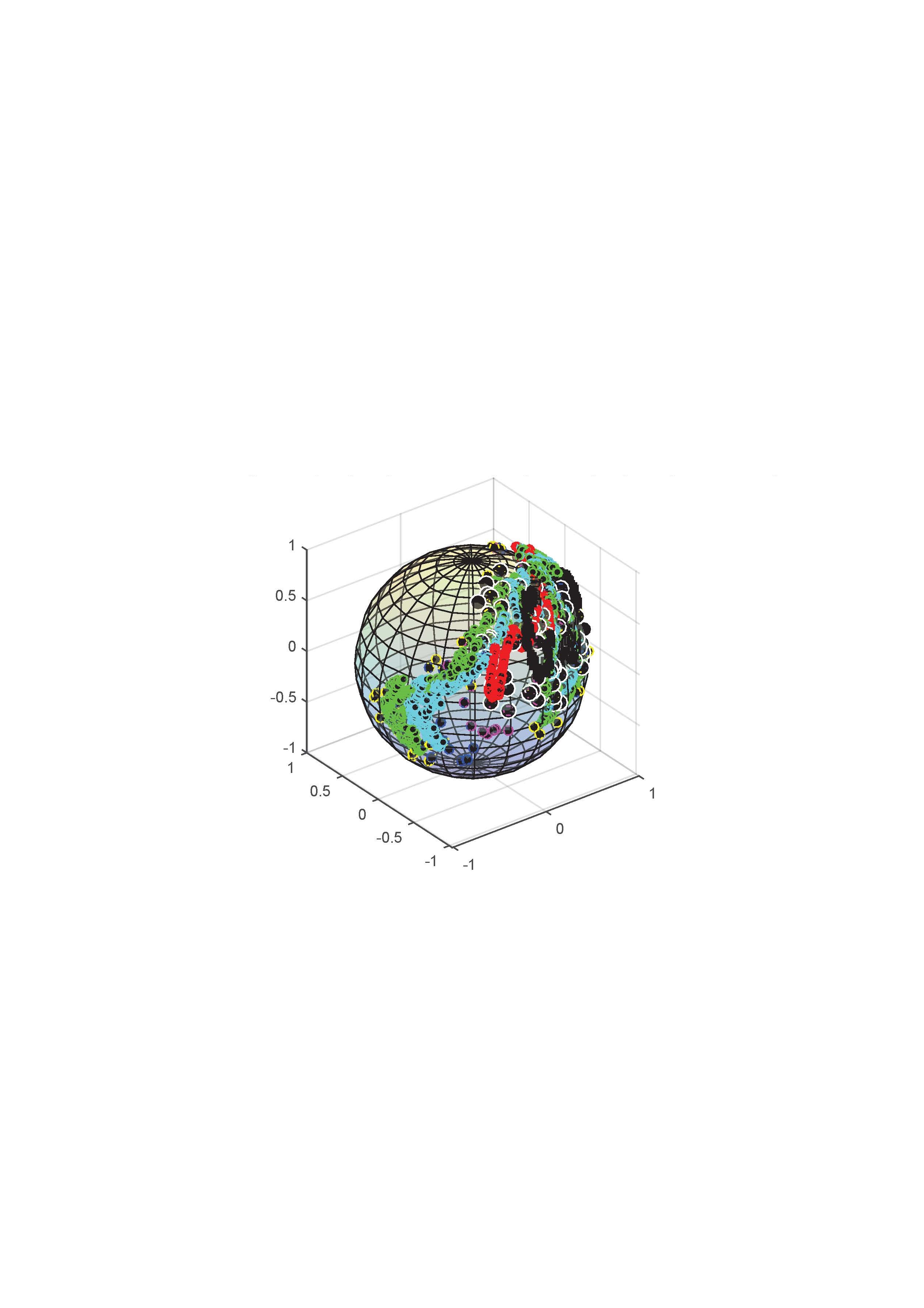}
\end{center}
\caption{The F\'echet predictor realizations, and the original  responses are displayed at times
$s_{i}=10, 20, 30$  (left--hand side), $s_{i}=40, 50, 60$ (center), and  $s_{i}=70 ,80, 90,100$ (right--hand side).}
\label{fig:6}
\end{figure}

\begin{figure}[!h]
\begin{center}
\includegraphics[height=5cm, width=5cm]{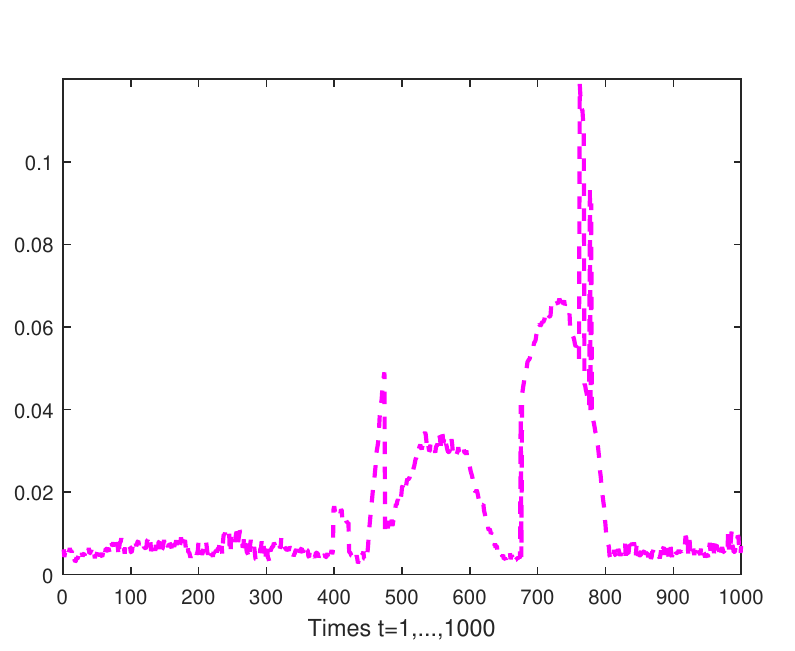}
\includegraphics[height=5cm, width=5cm]{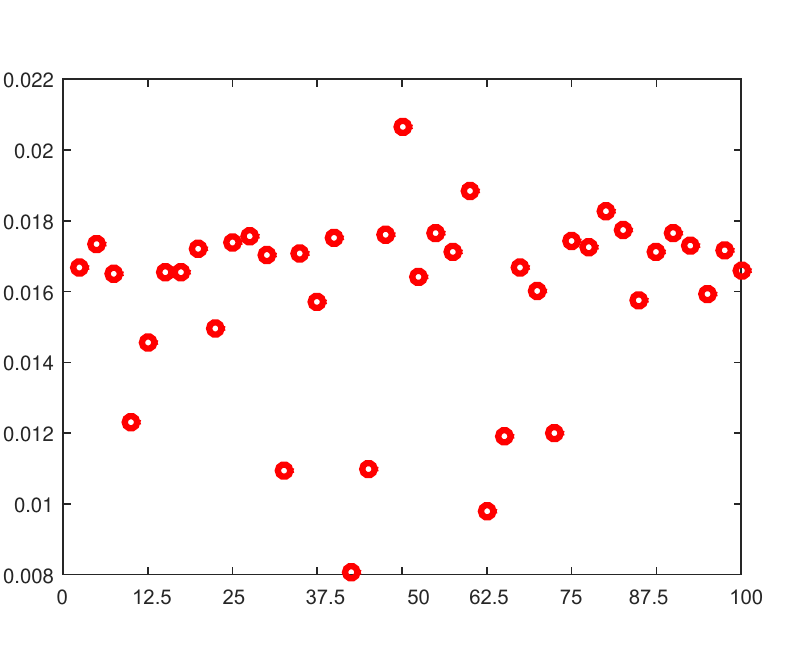}
\end{center}
\caption{The empirical mean of  the  quadratic  geodesic curve errors  (left--hand side), and the  empirical temporal means of their pointwise values at each sampled time
 $s_{i},$ $i=1,\dots,100$ (right--hand side).}
\label{fig:7}
\end{figure}

\begin{figure}[!h]
\begin{center}
\includegraphics[height=5cm, width=5cm]{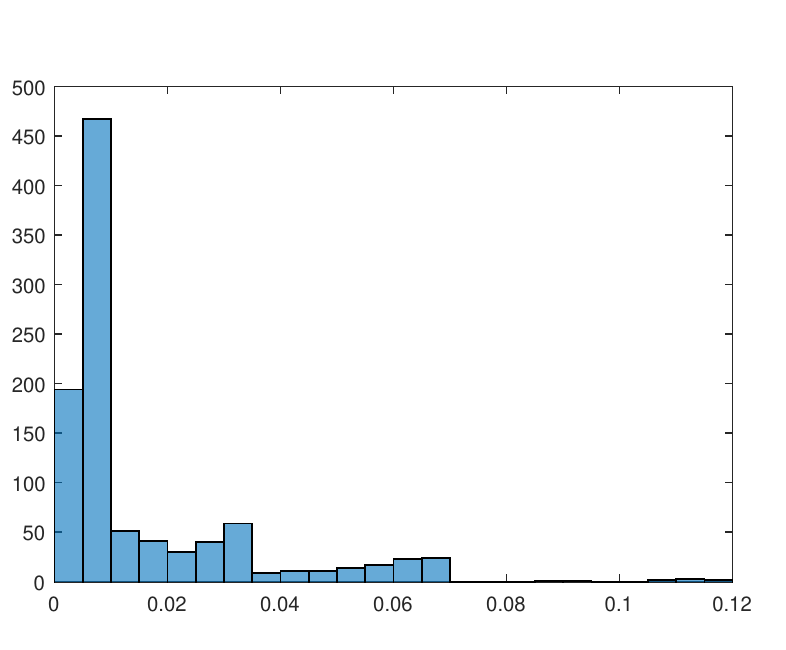}
\includegraphics[height=5cm, width=5cm]{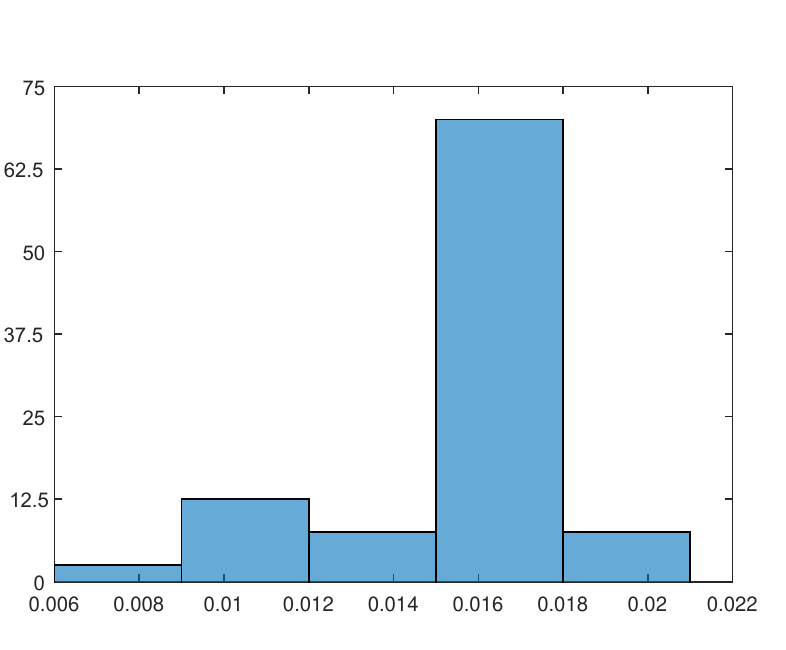}
\end{center}
\caption{Histogram of the  quadratic  geodesic curve error empirical mean values  at the left--hand side. The  histogram of the empirical  temporal   mean values  is also plotted at  the right--hand side.}
\label{fig:8}
\end{figure}

The empirical Fr\'echet mean $\widehat{\mu}_{X_{0},\mathcal{M}}(\cdot)$ in equation (\ref{efmvf}) is computed from a  uniform spherical grid of 400 nodes  (see Figure \ref{fig:4}). The  time--varying  Riemannian logarithm map with origin at $\widehat{\mu}_{X_{0},\mathcal{M}}(\cdot)$  is applied to the generated spherical curve regressor values. As commented,  these values are transformed via  the dynamic functional linear model
$$\mathbf{\Gamma }\left(\log_{\widehat{\mu}_{X_{0},\mathcal{M}}(t)}\left(X_{s_{i}}(\cdot )\right)\right)(t)+\varepsilon_{s_{i}}(t),\quad t\in \mathcal{T},\ i=1,\dots,n,$$
\noindent in the time--varying tangent space, leading  to the  log--mapped  curve values of the response in those tangent spaces
(see  plot at the right--hand side of Figure \ref{fig:4}). The time--varying exponential map at the same  origin  $\widehat{\mu}_{X_{0},\mathcal{M}}(\cdot)$ is then applied to them (see Figures \ref{fig:2}, and \ref{fig:3}).  The empirical Fr\'echet weights are also calculated as displayed in Figure \ref{fig:5}. Finally, the  empirical weighted Fr\'echet mean function is obtained. The corresponding Fr\'echet curve predictor is then computed as showed in Figure \ref{fig:6}.   A good performance is observed when a functional sample of size  $n=100$ is considered,  as displayed in Figures \ref{fig:7} and \ref{fig:8},  where Fr\'echet curve prediction errors are summarized.

\section{Real--data example}
\label{rde}
 It is well known that Earth's magnetic field protects Earth from solar wind that emanates from the sun. The three dimensional structure of this field is inferred from launched  satellite measurements using three--axis magnetometers. World magnetic models are generated from these remote sensors combined over the last few decades. In particular, navigation and heading referencing systems  can be improved from   accurate information on geomagnetic field. Data  from NASA's National Space Science Data Center are available in the period 02/11/1979--06/05/1980, recorded every half second, and correspond to the first satellite NASA's MAGSAT spacecraft, which orbited the earth  every 88 minutes during seven months at around 400 km altitude.   The available measurements during  the  days 3, 4, 5  of each month in the period 02/11/1979--06/05/1980  allow us to construct seven functional samples  of size 82, 84, 83, 84, 85, 72 and  52, respectively. The elements of these samples are discretely observed at  6000 temporal nodes.

 The
5--fold   cross validation technique  is implemented  from these data sets to assess the performance of the proposed   Fr\'echet   functional regression prediction methodology.    This section displays the results, based on  a functional sample  of size 82,    reflecting  the geocentric  latitude and longitude of the spacecraft  at 82 consecutive temporal intervals, containing 6000 equally spaced  temporal nodes, during the  days 3, 4, 5  of November, 1979, and reflecting the time--varying spherical coordinates of the magnetic field vector   at  the same temporal intervals and nodes. Both, the regressor and response  spherical curve observations, share the azimuthal angle, and display different time--varying polar angles, as given at the  top--right--hand--side plot in Figure 4 of Section 9.1 in \cite{Marzio.14}. Note that, in \cite{Marzio.14},  a  purely spatial analysis is carried out ignoring time information.  The 5--fold  cross validation results for  the remaining months are displayed in  Appendix \ref{app2}.

\begin{figure}[!h]
\begin{center}
\includegraphics[width=0.32\textwidth]{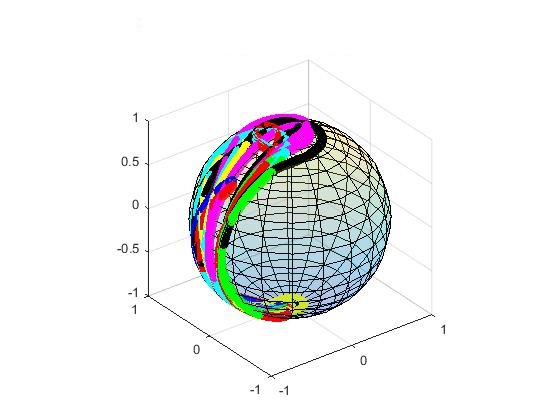}
\includegraphics[width=0.32\textwidth]{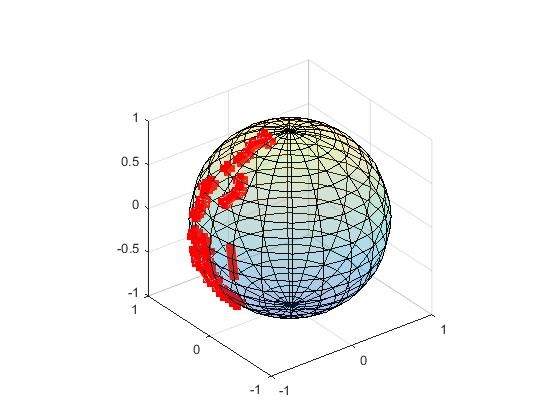}
\includegraphics[width=0.32\textwidth]{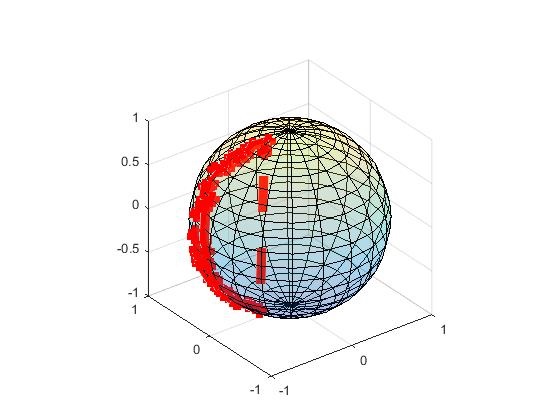}
\end{center}
\caption{Spherical bivariate curve data at times $t=1,11, 21, 31, 41, 51, 61, 71, 81$  (left--hand side). Empirical Fr\'echet curve mean of regressors (center) and of response (right--hand side).}\label{fig:0rda}
\end{figure}

\begin{figure}[!h]
\begin{center}
\includegraphics[width=0.32\textwidth]{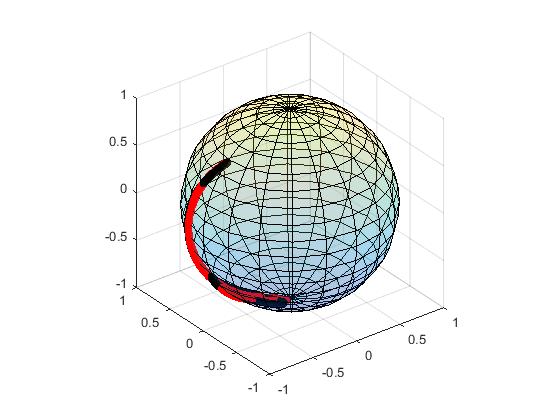}
\includegraphics[width=0.32\textwidth]{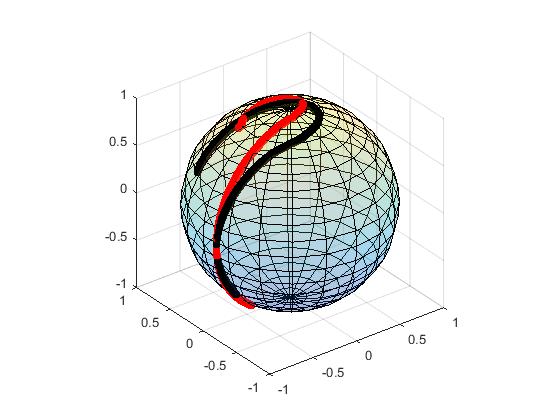}
\includegraphics[width=0.32\textwidth]{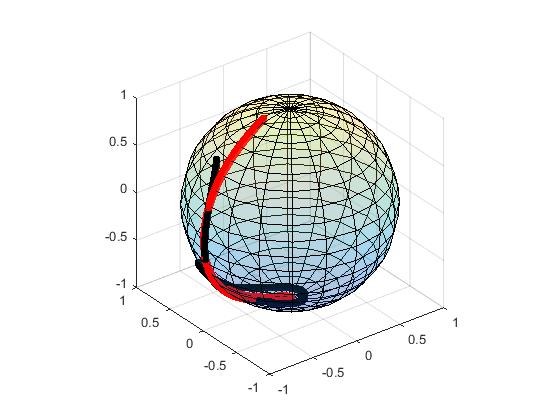}
\includegraphics[width=0.32\textwidth]{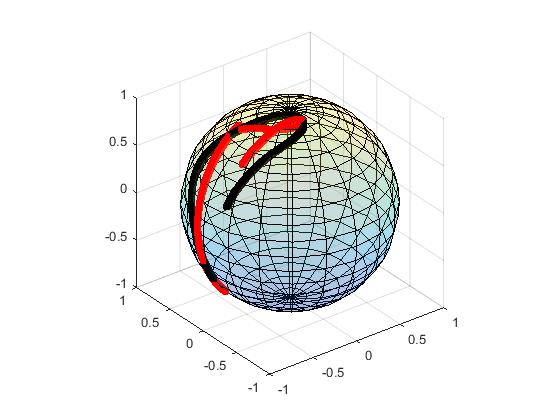}
\includegraphics[width=0.32\textwidth]{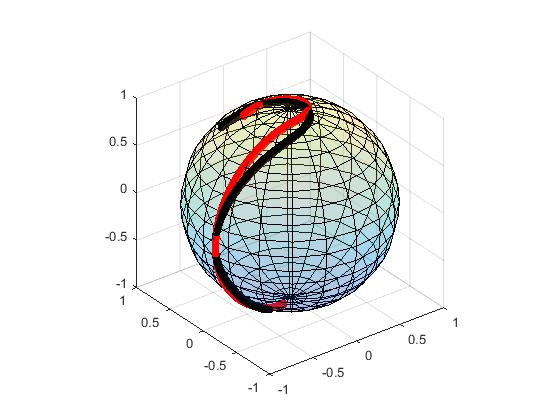}
\includegraphics[width=0.32\textwidth]{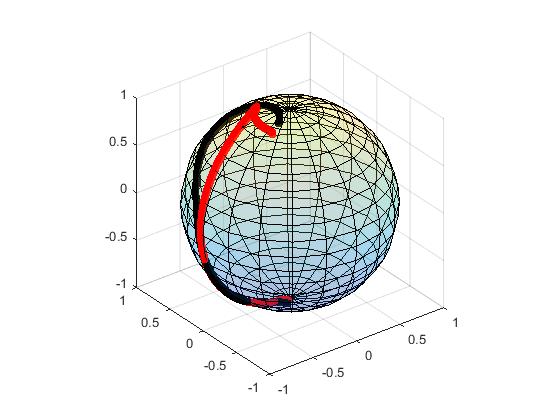}
\includegraphics[width=0.32\textwidth]{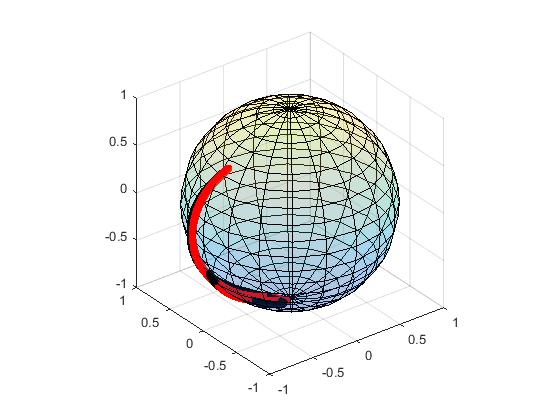}
\includegraphics[width=0.32\textwidth]{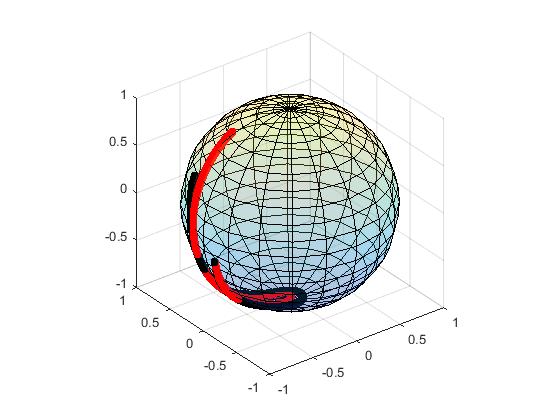}
\includegraphics[width=0.32\textwidth]{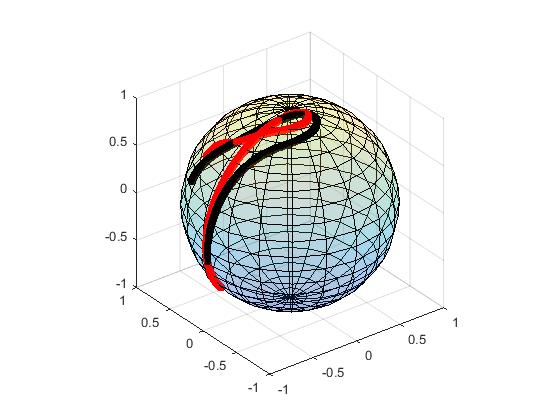}
\end{center}
\caption{Spherical bivariate curve data at times $t=1,12, 14, 25, 37, 49, 61, 73, 82$ (November, 1979). NASA's MAGSAT spacecraft (black curve), and   magnetic field vector (red curve) spherical coordinates are displayed at  $6000$ temporal nodes for every sampled time.}\label{fig:1rda}
\end{figure}

\begin{figure}[!h]
\begin{center}
\includegraphics[width=0.32\textwidth]{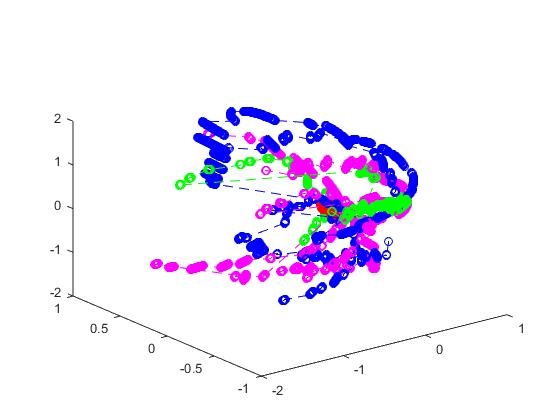}
\includegraphics[width=0.32\textwidth]{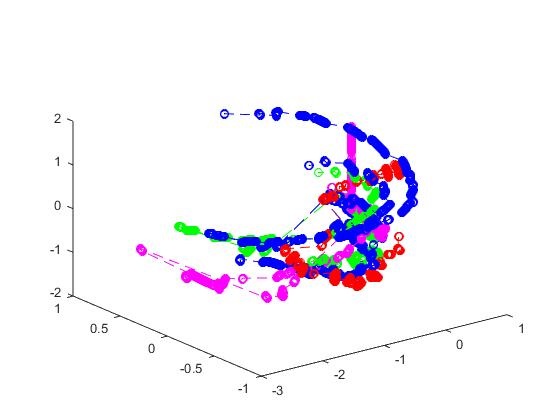}
\includegraphics[width=0.32\textwidth]{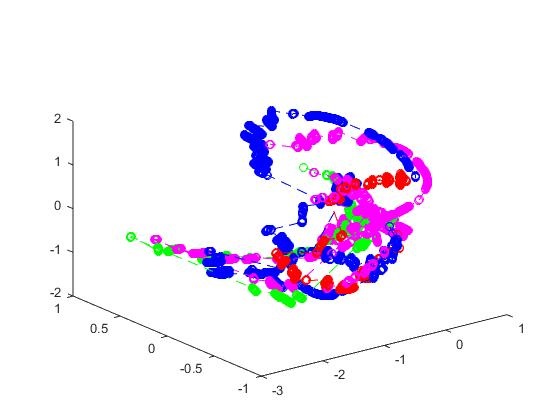}
\includegraphics[width=0.32\textwidth]{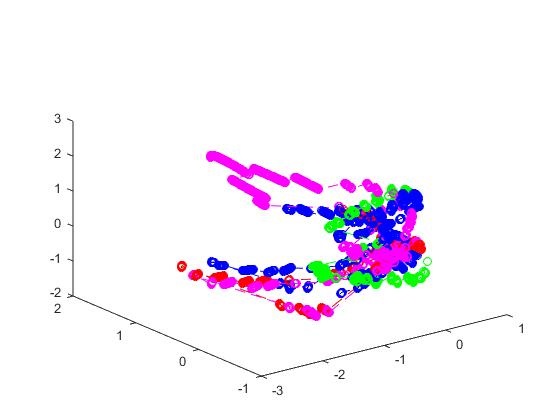}
\includegraphics[width=0.32\textwidth]{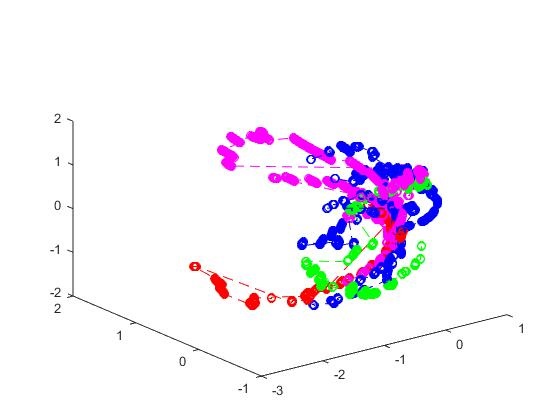}
\end{center}
\caption{Target   log-mapped regressor curve observations at times $t=1, 2, 4, 6, 8, 10$ for iteration one, at times  $t=1, 4, 7, 10, 18, 21$ for iteration 2, at times  $t=1,4, 7, 10, 13, 16$ for iteration 3, at times $t=1, 2, 3, 7, 10, 13$ for iteration 4, and at times $t=3, 5, 7, 10, 15, 21$ for iteration 5.}\label{fig:2rda}
\end{figure}

\begin{figure}[!h]
\begin{center}
\includegraphics[width=0.32\textwidth]{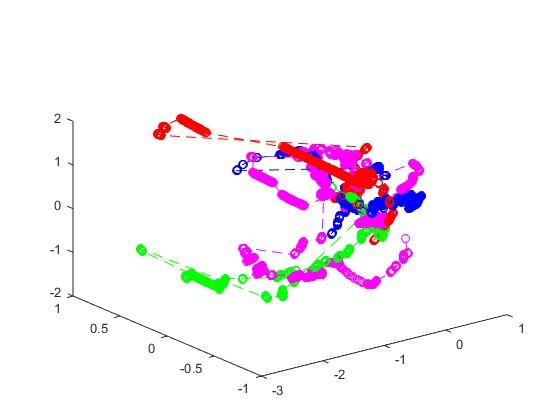}
\includegraphics[width=0.32\textwidth]{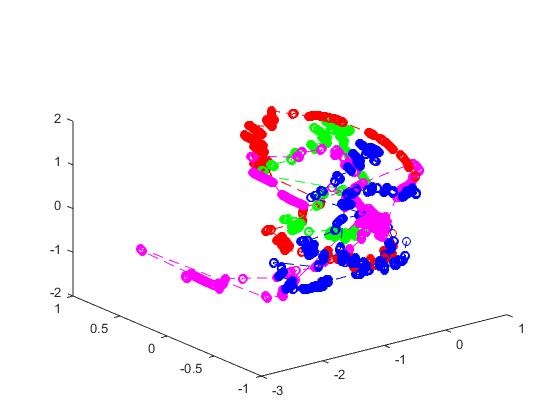}
\includegraphics[width=0.32\textwidth]{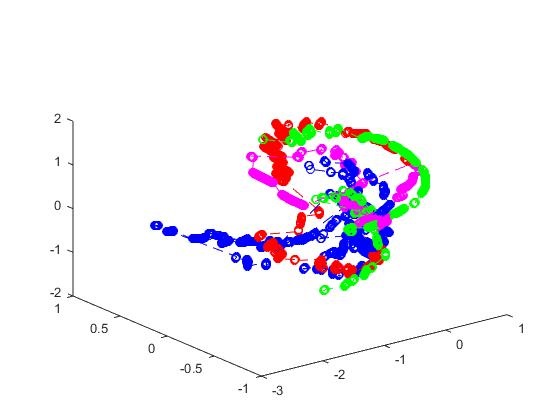}
\includegraphics[width=0.32\textwidth]{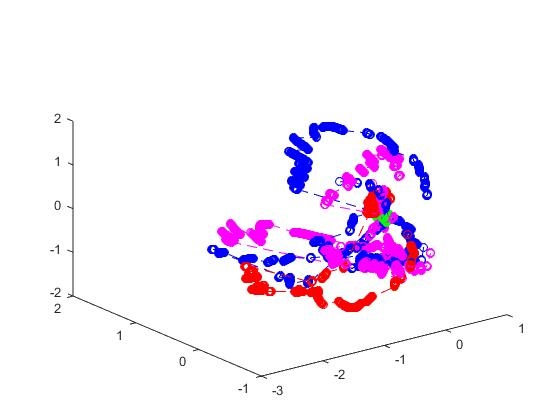}
\includegraphics[width=0.32\textwidth]{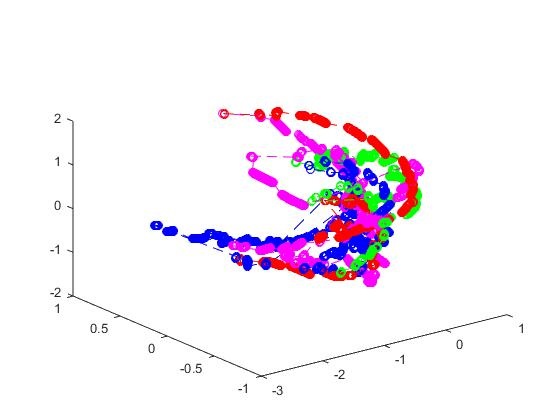}
\end{center}
\caption{Training  log-mapped regressor curve observations at times $t=1, 3, 6, 9, 12, 72$ for iteration one, at times  $t=1, 3, 6, 9, 12, 61$ for iteration 2, at times  $t=1, 3, 6, 9, 12, 66$ for iteration 3, at times $t=10, 20, 30, 40, 50, 69$ for iteration 4, and at times $t=10, 20, 30, 40, 50, 60$ for iteration 5.}\label{fig:3rda}
\end{figure}

 As in the previous section, conditions (i)-(ii) are satisfied by $\mathcal{M}=\mathbb{S}_{2}\subset\mathbb{R}^{3}.$  Figure  \ref{fig:0rda} displays the spherical bivariate   curve  observations at
times $t=1,11, 21, 31, 41, 51, 61, 71, 81$ at the left--hand side (see also Figure \ref{fig:1rda}, where  some bivariate curve observations are displayed for  different times). From these plots, condition (iii)--(v) seems to be satisfied  (see Remark \ref{remcondint}). The
  empirical intrinsic Fr\'echet curve means, based on  the  regressors and the response curve observations,  are respectively shown at the center, and at the right--hand side of Figure  \ref{fig:0rda}. These plots also agree with the assumption of  a common  theoretical Fr\'echet curve  mean for the marginal probability measure of the response and regressor process, as given in condition (v).  The 5--fold   cross validation technique is implemented from  the log--mapped  training and target spherical curve regressor subsamples in the time--varying tangent space  (see Figures \ref{fig:2rda} and  \ref{fig:3rda}), and from the corresponding training and target spherical curve response subsamples in $\mathcal{M}.$
   The displayed running of the  5--fold  cross validation   algorithm involves   the  training samples of sizes    $72, 61,    66,    69,    60,$
  and target samples of sizes $10, 21, 16, 13,  22,$ respectively.

\begin{figure}[!h]
\begin{center}
\includegraphics[width=0.32\textwidth]{GFCMREG}
\includegraphics[width=0.32\textwidth]{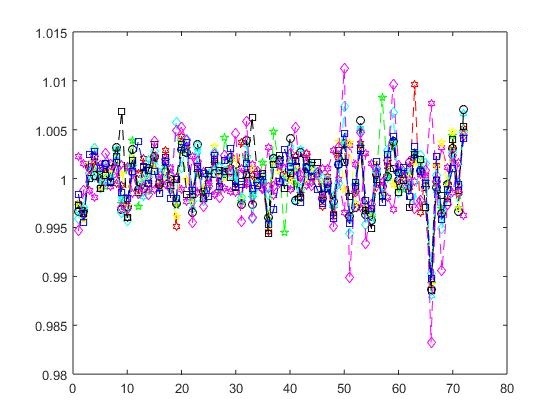}
\includegraphics[width=0.32\textwidth]{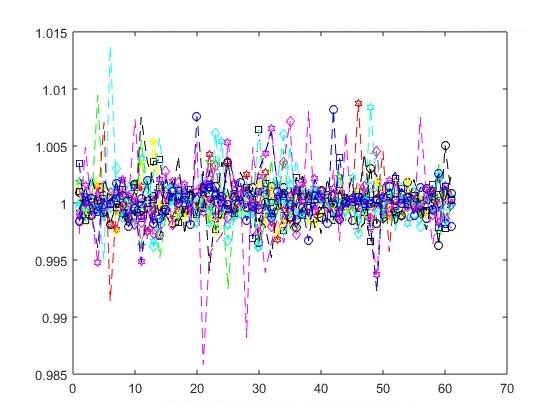}
\includegraphics[width=0.32\textwidth]{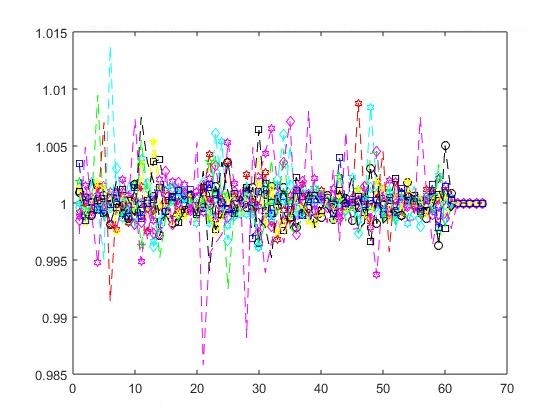}
\includegraphics[width=0.32\textwidth]{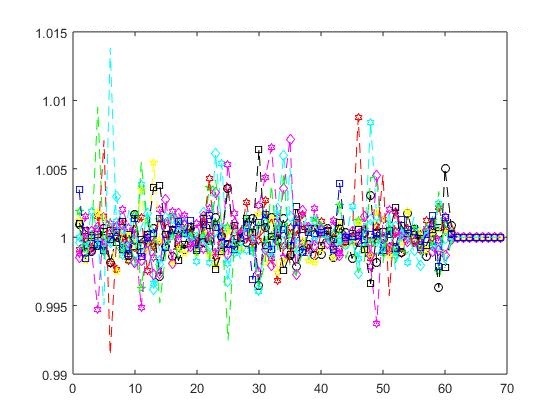}
\includegraphics[width=0.32\textwidth]{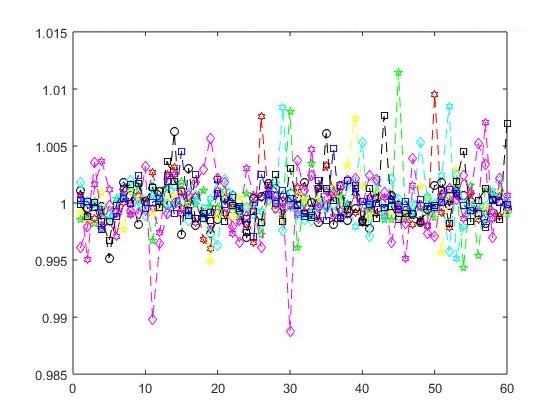}
\end{center}
\caption{Empirical Fr\'echet curve mean based on  $82$  spherical curve regressor observations,  computed from a uniform spherical grid with 1000 nodes  (top--left--hand--side).
The remaining plots provide Fr\'echet weights at each one of the five iterations of the 5--fold  cross validation algorithm implemented.}\label{fig:4rda}
\end{figure}

\begin{figure}[!h]
\begin{center}
\includegraphics[width=0.32\textwidth]{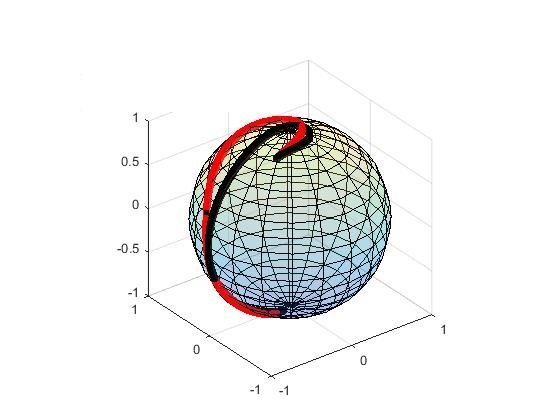}
\includegraphics[width=0.32\textwidth]{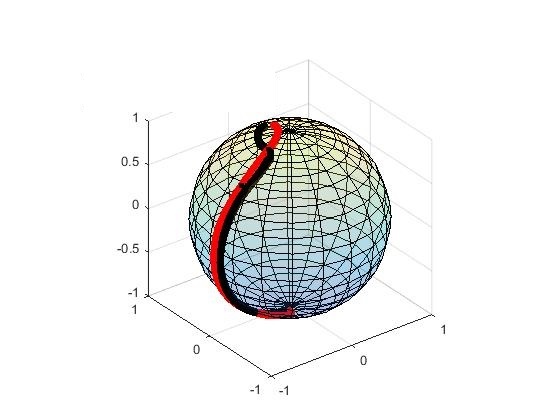}
\includegraphics[width=0.32\textwidth]{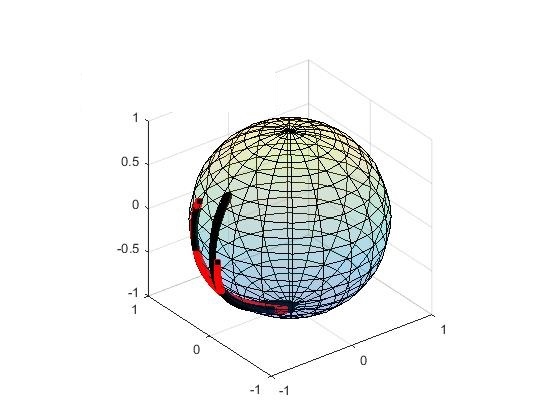}
\end{center}
\caption{Spherical Fr\'echet  functional   predictor (red curve) and response (black curve) at times  $t_{i},$ $i=1,2,3$ (corresponding to original observed times $t=  4, 11, 75$) of the  target subsample  at the first iteration of the 5--fold cross validation algorithm. }\label{fig:5rda}
\end{figure}

\begin{figure}[!h]
\begin{center}
\includegraphics[width=0.32\textwidth]{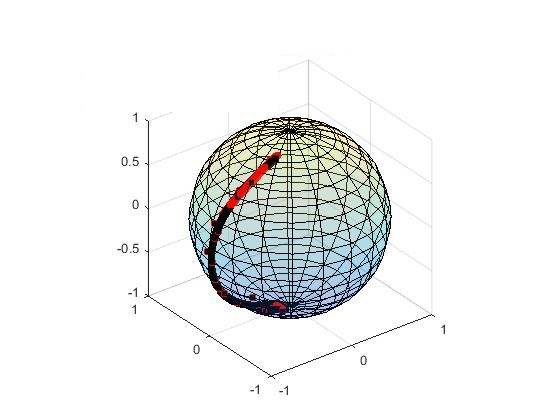}
\includegraphics[width=0.32\textwidth]{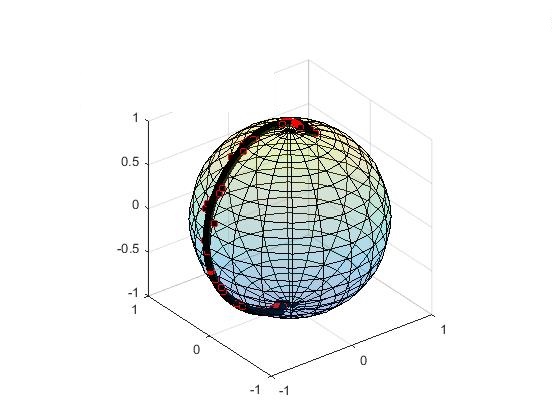}
\includegraphics[width=0.32\textwidth]{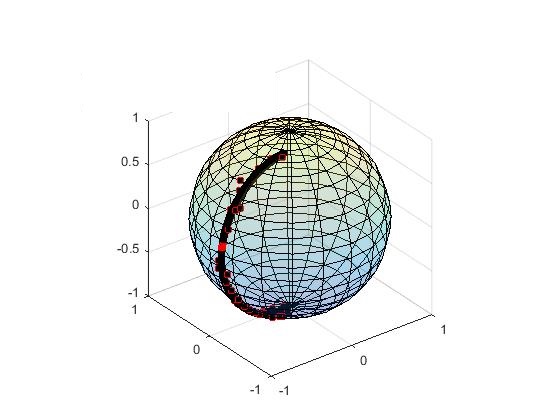}
\end{center}
\caption{Spherical Fr\'echet  functional   predictor (red curve) and response (black curve) at times  $t_{i},$ $i=1,2,3$ (corresponding to original observed times $t=  3, 6, 2$) of the  target subsample  at the second iteration of the 5--fold cross validation algorithm. }
\label{fig:6rda}
\end{figure}

\begin{figure}[!h]
\begin{center}
\includegraphics[width=0.32\textwidth]{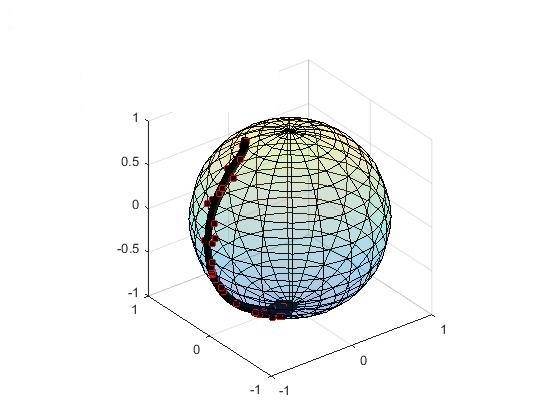}
\includegraphics[width=0.32\textwidth]{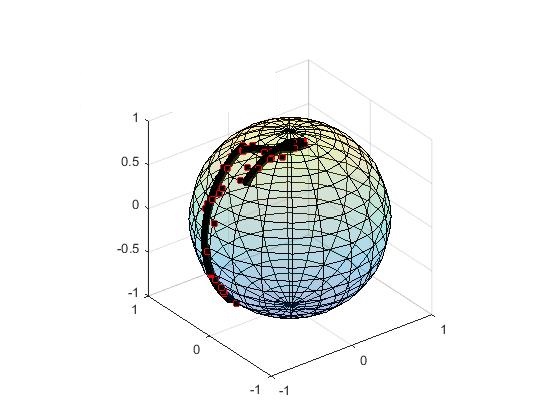}
\includegraphics[width=0.32\textwidth]{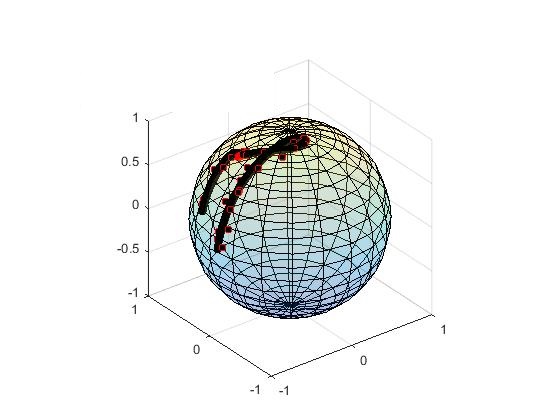}
\end{center}
\caption{Spherical Fr\'echet functional   predictor (red curve) and response (black curve) at times  $t_{i},$ $i=1,2,3$ (corresponding to original observed times $t=  10 ,25, 58$) of the  target subsample  at the third iteration of the 5--fold cross validation algorithm. }
\label{fig:7rda}
\end{figure}

\begin{figure}[!h]
\begin{center}
\includegraphics[width=0.32\textwidth]{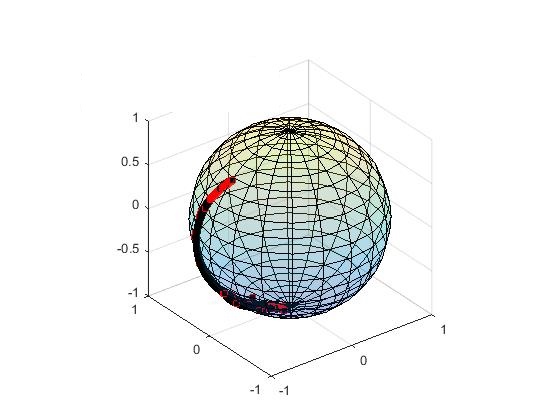}
\includegraphics[width=0.32\textwidth]{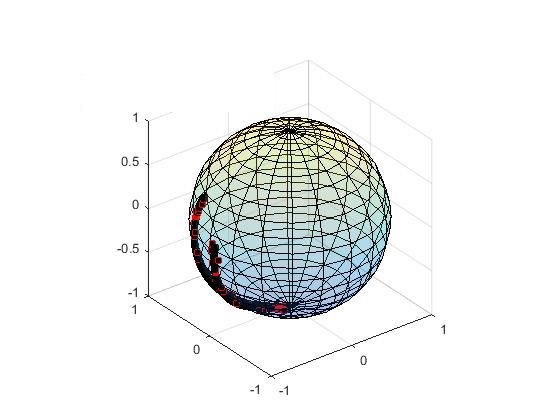}
\includegraphics[width=0.32\textwidth]{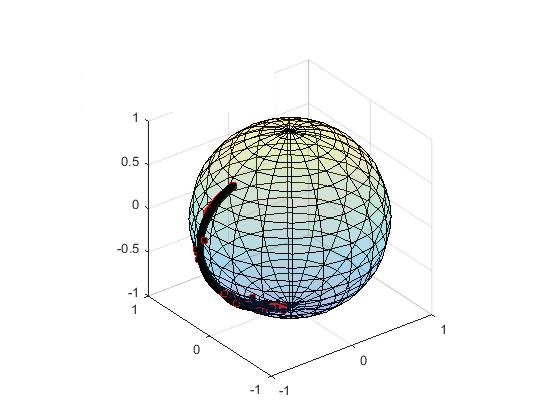}
\end{center}
\caption{Spherical Fr\'echet functional   predictor (red curve) and response (black curve) at times  $t_{i},$ $i=1,2,3$ (corresponding to original observed times $t=  1, 28, 61$) of the  target subsample  at the fourth iteration of the 5--fold cross validation algorithm. }
\label{fig:8rda}
\end{figure}

\begin{figure}[!h]
\begin{center}
\includegraphics[width=0.32\textwidth]{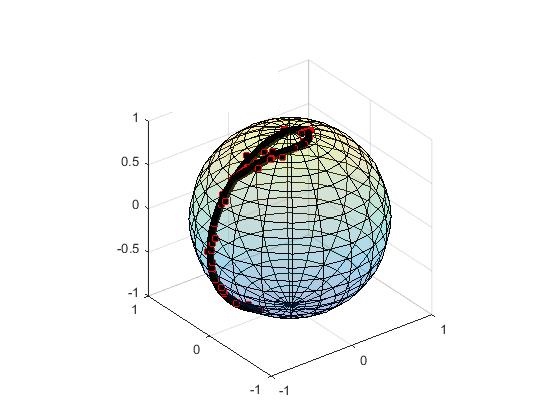}
\includegraphics[width=0.32\textwidth]{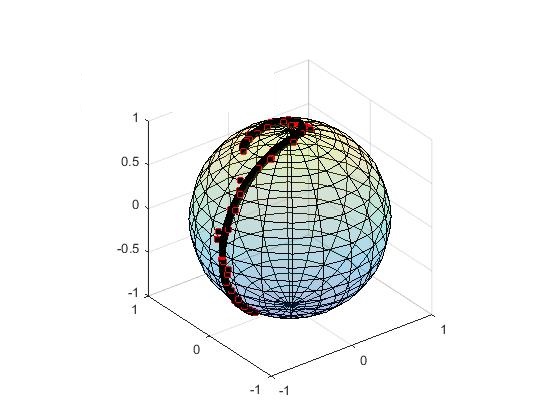}
\includegraphics[width=0.32\textwidth]{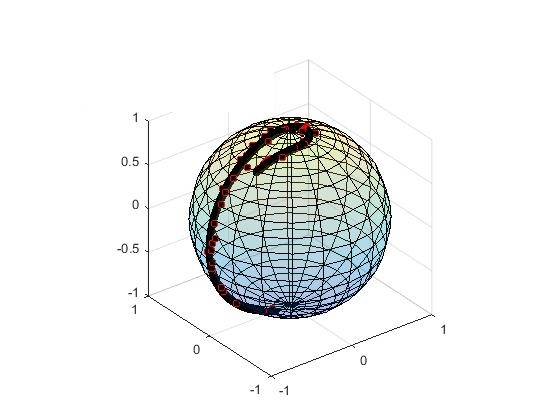}
\end{center}
\caption{Spherical  Fr\'echet  functional predictor (red curve) and response (black curve) at times  $t_{i},$ $i=1,2,3,$ (corresponding to original observed times $t=  2, 12, 31$) of the  target subsample  at the fifth iteration of the 5--fold cross validation algorithm. }
\label{fig:9rda}
\end{figure}

 Fr\'echet weights are computed  from the training and target log--mapped regressor subsamples   (see  Figure \ref{fig:4rda}). Specifically, a discretized version of   the empirical matrix covariance operator of the log--mapped curve regressors is computed from the   training regressor subsample  at each iteration of the 5--fold  cross validation. Note that the row and column input vectors  of the   block matrix, approximating  this empirical operator,  are respectively evaluated in the training and target regressor subsamples.
At each one of the 5 iterations of the 5--fold  cross validation algorithm, the original spherical curve response values, and the corresponding spherical  Fr\'echet functional  predictions    are  plotted in Figures \ref{fig:5rda}--\ref{fig:9rda} at the first three target times. Finally, the average over the 6000 temporal nodes of  the 5--fold cross validation functional  error   is    displayed in Table  \ref{T1},  for the $22$ Fr\'echet spherical curve  predictions  at   target times  $t=2,  8,  12, 15,  22, 24,   26,    29,31,  33,    40,    43,    45,    47,    53,    55,    56,    65,    67,    71,    76,  77.$  This \linebreak   average is denoted as TPAF5fCVE  in Table  \ref{T1}.

\clearpage

    \begin{table}[!h]
\caption{\textbf{\emph{Temporal  pointwise average of 5--fold   cross validation functional  error} }}\label{T1}
\begin{center}
\begin{tabular}{|c|c|c|c|}
\hline
 {\bf TIMES}& {\bf TPAF5fCVE} & {\bf TIMES} & {\bf TPAF5fCVE} \\
  \hline
T1 & 0.0037 & T12 &0.0027\\
T2 &  0.0036  &T13 & 0.0030\\
T3 & 0.0034 & T14 & 0.0023\\
T4 &0.0037 &T15 & 0.0022  \\
T5 & 0.0031&T16&  0.0021 \\
T6 & 0.0039 &T17&  0.0013\\
T7 & 0.0039&T18&  0.0016\\
T8 &  0.0036 &T19&  0.0015 \\
T9 & 0.0040 &T20& 0.0016\\
T10 &  0.0036&T21&  0.0014\\
T11 & 0.0030&T22& 0.0006 \\
\hline
 \end{tabular}
\end{center}
\end{table}

\section{Final Comments}
\label{corl}
This paper deals with  global  intrinsic  Fr\'echet regression in the  framework of $\mathcal{C}_{\mathcal{M}}(\mathcal{T})$--valued   bivariate   curve processes, adopting the weighted Fr\'echet mean formulation introduced
in  \cite{Petersen.19}.  Thus, we extend this formulation to the case of time correlated  bivariate curve data evaluated in a compact Riemannian manifold.  In particular, we extend the Euclidean regressor setting in  \cite{Petersen.19} to the $\mathcal{M}$--valued curve process framework, considering the metric  space $\left(\mathcal{C}_{\mathcal{M}}(\mathcal{T}),d_{\mathcal{C}_{\mathcal{M}}(\mathcal{T})}\right).$  Local linear Fr\'echet curve regression will be addressed in the manifold--valued bivariate curve process framework in a subsequent paper  (see  \cite{Petersen.19}, and  \cite{Marzio.14} in the  context of the local polynomial regression framework).

\newpage

\section*{Appendices}
The proof of the main   result, Theorem  \ref{th1},  is provided in Appendix \ref{app1}.  Appendix \ref{app2} shows  the 5--fold cross validation results obtained since December, 1979, to May, 1980.
\begin{appendix}
\section{Proof of Theorem \ref{th1}}
\label{app1}
\begin{proof}
For a fixed $x(\cdot)\in \mathcal{X}_{\mathcal{C}_{\mathcal{M}}(\mathcal{T})},$ under condition A.1,  from Corollary
3.2.3 in \cite{Vaart.1996}, it is sufficient to prove the convergence to zero in probability of
\begin{equation}
\sup_{z(\cdot)\in \mathcal{Y}_{\mathcal{C}_{\mathcal{M}}(\mathcal{T})}}\left|\widehat{M}_{n}\left(z(\cdot), x(\cdot )\right)-
M\left(z(\cdot),x(\cdot )\right)\right|.\label{ie}
\end{equation}
\noindent Indeed,   from Theorem 1.5.4  in \cite{Vaart.1996}, it is sufficient to prove
\begin{itemize}
\item[(1)]  For each $x\in \mathcal{X}_{\mathcal{C}_{\mathcal{M}}(\mathcal{T})},$ and for every $z(\cdot)\in \mathcal{Y}_{\mathcal{C}_{\mathcal{M}}(\mathcal{T})},$   $$\widehat{M}_{n}(z(\cdot ), x(\cdot ))-M(z(\cdot), x(\cdot ))=o_{P}(1),\ n\to \infty.$$
\item[(2)]  For each $x\in \mathcal{X}_{\mathcal{C}_{\mathcal{M}}(\mathcal{T})},$ and $z(\cdot)\in \mathcal{Y}_{\mathcal{C}_{\mathcal{M}}(\mathcal{T})},$  $\widehat{M}_{n}(z(\cdot ), x(\cdot ))$ is asymptotically equicontinuos in probability, i.e., for all $\varepsilon >0$  and $\eta,$ there exists a $\delta >0$ which does not depend on $z(\cdot)$ such that
    \begin{eqnarray}
    \limsup_{n} P\left(\sup_{d_{\mathcal{C}_{\mathcal{M}}(\mathcal{T})}(z(\cdot),y(\cdot))\leq \delta }\left|\widehat{M}_{n}(z(\cdot ), x(\cdot ))-\widehat{M}_{n}(y(\cdot), x(\cdot ))\right|>\varepsilon\right)<\eta .\nonumber\\
    \label{eqc}
    \end{eqnarray}\end{itemize}

    To prove (1), we consider \begin{eqnarray}&&
\widetilde{M}_{n}\left(z(\cdot), x(\cdot )\right)=\frac{1}{n}\sum_{i=1}^{n} \int_{\mathcal{T}}\left[d_{\mathcal{M}}\left( Y_{i }(t),z(t)\right)\right]^{2}dt
\nonumber\\
&&\hspace*{-0.5cm}\times \left[1+\left\langle \sqrt{\mathcal{K}}\left(\log_{\mu_{X_{0},\mathcal{M}}(\cdot)}\left(x(\cdot )\right)-\mu(\cdot )\right), \mathcal{R}_{X}^{-1}\left(\sqrt{\mathcal{K}}\left(\log_{\mu_{X_{0},\mathcal{M}}(\cdot)}\left(X_{i}(\cdot )\right)-\mu(\cdot )\right)\right)\right\rangle_{\mathbb{H}}\right].\nonumber
\end{eqnarray}

Under (iv), applying strictly stationarity,  for each $x(\cdot)\in \mathcal{X}_{\mathcal{C}_{\mathcal{M}}(\mathcal{T})},$ and every $z(\cdot)\in \mathcal{Y}_{\mathcal{C}_{\mathcal{M}}(\mathcal{T})},$ $E\left[\widetilde{M}_{n}\left(z(\cdot), x(\cdot )\right)\right]=M(z(\cdot), x(\cdot )).$ We also obtain
\begin{eqnarray}&&
\mbox{Var}\left(\widetilde{M}_{n}\left(z(\cdot), x(\cdot )\right)\right)=\frac{1}{n}\sum_{u\in \{-(n-1),\dots, n-1\}}\left(1-\frac{|u|}{n}\right)
\nonumber\\
&&\hspace*{0.5cm}\times \left[\mathcal{R}_{u}^{(g(Y),g(Y))}+\mathcal{R}_{u}^{(g(Y),g(Y)h(X))}+\mathcal{R}_{u}^{(g(Y)h(X),g(Y))}
+\mathcal{R}_{u}^{(g(Y)h(X),g(Y)h(X))}\right]
\nonumber\\
&&\leq \frac{1}{n}\sum_{u\in \mathbb{Z}}\left[\mathcal{R}_{u}^{(g(Y),g(Y))}+\mathcal{R}_{u}^{(g(Y),g(Y)h(X))}+\mathcal{R}_{u}^{(g(Y)h(X),g(Y))}
+\mathcal{R}_{u}^{(g(Y)h(X),g(Y)h(X))}\right],\nonumber\\
\label{fth1proof}
\end{eqnarray}
\noindent where
\begin{eqnarray}&&
\mathcal{R}_{u}^{(g(Y),g(Y))}= E\left[\left(\int_{\mathcal{T}}\left[d_{\mathcal{M}}\left(Y_{0}(t),z(t)\right)\right]^{2}dt\right)
\left(\int_{\mathcal{T}}\left[d_{\mathcal{M}}\left(Y_{u}(t),z(t)\right)\right]^{2}dt\right)\right]\nonumber\\
&&\mathcal{R}_{u}^{(g(Y),g(Y)h(X))}=E\left[\left(\int_{\mathcal{T}}\left[d_{\mathcal{M}}\left(Y_{0}(t),z(t)\right)\right]^{2}dt\right)
\left(\int_{\mathcal{T}}\left[d_{\mathcal{M}}\left(Y_{u}(t),z(t)\right)\right]^{2}dt\right)\right.\nonumber\\
&& \left. \times  \left\langle \sqrt{\mathcal{K}}\left(\log_{\mu_{X_{0},\mathcal{M}}(\cdot)}\left(x(\cdot )\right)-\mu (\cdot)\right), \sqrt{\mathcal{K}}\left(\log_{\mu_{X_{0},\mathcal{M}}(\cdot)}\left(X_{u}(\cdot )\right)-\mu (\cdot)\right)\right\rangle_{\widetilde{H}}\right] \nonumber\\
&&\mathcal{R}_{u}^{(g(Y)h(X),g(Y))}=E\left[\left(\int_{\mathcal{T}}\left[d_{\mathcal{M}}\left(Y_{0}(t),z(t)\right)\right]^{2}dt\right)
\left(\int_{\mathcal{T}}\left[d_{\mathcal{M}}\left(Y_{u}(t),z(t)\right)\right]^{2}dt\right)\right.\nonumber\\
&& \left. \times  \left\langle \sqrt{\mathcal{K}}\left(\log_{\mu_{X_{0},\mathcal{M}}(\cdot)}\left(x(\cdot )\right)-\mu (\cdot)\right), \sqrt{\mathcal{K}}\left(\log_{\mu_{X_{0},\mathcal{M}}(\cdot)}\left(X_{0}(\cdot )\right)-\mu (\cdot)\right)\right\rangle_{\widetilde{H}}\right]
\nonumber\\
&&\mathcal{R}_{u}^{(g(Y)h(X),g(Y)h(X))}=E\left[\left(\int_{\mathcal{T}}\left[d_{\mathcal{M}}\left(Y_{0}(t),z(t)\right)\right]^{2}dt\right)
\left(\int_{\mathcal{T}}\left[d_{\mathcal{M}}\left(Y_{u}(t),z(t)\right)\right]^{2}dt\right)\right.\nonumber\\
&& \left. \times  \left\langle \sqrt{\mathcal{K}}\left(\log_{\mu_{X_{0},\mathcal{M}}(\cdot)}\left(x(\cdot )\right)-\mu (\cdot)\right), \sqrt{\mathcal{K}}\left(\log_{\mu_{X_{0},\mathcal{M}}(\cdot)}\left(X_{0}(\cdot )\right)-\mu (\cdot)\right)\right\rangle_{\widetilde{H}}\right.
\nonumber\\
&&\left. \times  \left\langle \sqrt{\mathcal{K}}\left(\log_{\mu_{X_{0},\mathcal{M}}(\cdot)}\left(x(\cdot )\right)-\mu (\cdot)\right), \sqrt{\mathcal{K}}\left(\log_{\mu_{X_{0},\mathcal{M}}(\cdot)}\left(X_{u}(\cdot )\right)-\mu (\cdot)\right)\right\rangle_{\widetilde{H}}\right].
\label{eeq2pth1}
\end{eqnarray}

\noindent  Under (i)--(v),  keeping in mind Remark \ref{rem1},   the following inequalities hold:
\begin{eqnarray}&&
\mathcal{R}_{u}^{(g(Y),g(Y))}\leq E\left[\left\|\log_{\mu_{X_{0},\mathcal{M}}(\cdot)}(Y_{0}(\cdot))-
    \log_{\mu_{X_{0},\mathcal{M}}(\cdot)}(z(\cdot))\right\|_{\mathbb{H}}^{2}\right.
    \nonumber\\
&&
    \hspace*{2cm} \left. \times \left\|\log_{\mu_{X_{0},\mathcal{M}}(\cdot)}(Y_{u}(\cdot))-
    \log_{\mu_{X_{0},\mathcal{M}}(\cdot)}(z(\cdot))\right\|_{\mathbb{H}}^{2}\right]\label{ubth10}\\
    &&\mathcal{R}_{u}^{(g(Y),g(Y)h(X))}\leq \|\sqrt{\mathcal{K}}\mathcal{R}_{X}^{-1}\sqrt{\mathcal{K}}\|_{\mathcal{L}(\mathbb{H})}
    N(x(\cdot))\nonumber\\
   && \hspace*{2cm} \times  E\left[\left\|\log_{\mu_{X_{0},\mathcal{M}}(\cdot)}(Y_{0}(\cdot))-
    \log_{\mu_{X_{0},\mathcal{M}}(\cdot)}(z(\cdot))\right\|_{\mathbb{H}}^{2}\right.
    \nonumber\\
&&
    \hspace*{2cm} \left. \times \left\|\log_{\mu_{X_{0},\mathcal{M}}(\cdot)}(Y_{u}(\cdot))-
    \log_{\mu_{X_{0},\mathcal{M}}(\cdot)}(z(\cdot))\right\|_{\mathbb{H}}^{2}\right.\nonumber\\
    && \hspace*{2cm} \left. \times
   \left\| \log_{\mu_{X_{0},\mathcal{M}}(\cdot)}\left(X_{u}(\cdot )\right)-\mu (\cdot)\right\|_{\mathbb{H}}
     \right]\label{ubth11}\\
     &&\mathcal{R}_{u}^{(g(Y)h(X),g(Y))}\leq \|\sqrt{\mathcal{K}}\mathcal{R}_{X}^{-1}\sqrt{\mathcal{K}}\|_{\mathcal{L}(\mathbb{H})}
    N(x(\cdot))\nonumber\\
   && \hspace*{2cm} \times  E\left[\left\|\log_{\mu_{X_{0},\mathcal{M}}(\cdot)}(Y_{0}(\cdot))-
    \log_{\mu_{X_{0},\mathcal{M}}(\cdot)}(z(\cdot))\right\|_{\mathbb{H}}^{2}\right.
    \nonumber\\
&&
    \hspace*{2cm} \left. \times \left\|\log_{\mu_{X_{0},\mathcal{M}}(\cdot)}(Y_{u}(\cdot))-
    \log_{\mu_{X_{0},\mathcal{M}}(\cdot)}(z(\cdot))\right\|_{\mathbb{H}}^{2}\right.\nonumber\\
    && \hspace*{2cm} \left. \times
   \left\| \log_{\mu_{X_{0},\mathcal{M}}(\cdot)}\left(X_{0}(\cdot )\right)-\mu (\cdot)\right\|_{\mathbb{H}}\right]\label{ubth12}\\
   &&\mathcal{R}_{u}^{(g(Y)h(X),g(Y)h(X))}\leq
   \left[\|\sqrt{\mathcal{K}}\mathcal{R}_{X}^{-1}\sqrt{\mathcal{K}}\|_{\mathcal{L}(\mathbb{H})}
    N(x(\cdot))\right]^{2}\nonumber\\
   && \hspace*{2cm} \times  E\left[\left\|\log_{\mu_{X_{0},\mathcal{M}}(\cdot)}(Y_{0}(\cdot))-
    \log_{\mu_{X_{0},\mathcal{M}}(\cdot)}(z(\cdot))\right\|_{\mathbb{H}}^{2}\right.
    \nonumber\\
&&
    \hspace*{2cm} \left. \times \left\|\log_{\mu_{X_{0},\mathcal{M}}(\cdot)}(Y_{u}(\cdot))-
    \log_{\mu_{X_{0},\mathcal{M}}(\cdot)}(z(\cdot))\right\|_{\mathbb{H}}^{2}\right.\nonumber\\
    && \hspace*{2cm} \left. \times
   \left\| \log_{\mu_{X_{0},\mathcal{M}}(\cdot)}\left(X_{0}(\cdot )\right)-\mu (\cdot)\right\|_{\mathbb{H}}\right.\nonumber\\
   && \hspace*{2cm} \left. \times
   \left\| \log_{\mu_{X_{0},\mathcal{M}}(\cdot)}\left(X_{u}(\cdot )\right)-\mu (\cdot)\right\|_{\mathbb{H}}\right],
\label{eq3pth1}
\end{eqnarray}
\noindent where  Cauchy--Schwarz inequality has been applied in the derivation of the upper bounds (\ref{ubth11})--(\ref{eq3pth1}), keeping in mind condition (a)--(b) in Section \ref{s32}, ensuring the bounded operator norm  $\|\sqrt{\mathcal{K}}\mathcal{R}_{X}^{-1}\sqrt{\mathcal{K}}\|_{\mathcal{L}(\mathbb{H})}$ on $\mathbb{H}$ is finite. Here,
 $N(x(\cdot))=  \left\|\log_{\mu_{X_{0},\mathcal{M}}(\cdot)}\left(x(\cdot )\right)-\mu (\cdot)\right\|_{\mathbb{H}}.$
   Under C.1, from (\ref{fth1proof}), equations (\ref{ubth10})--(\ref{eq3pth1}) mean that $$\widetilde{M}_{n}\left(z(\cdot), x(\cdot ) \right)-M\left(z(\cdot), x(\cdot )\right)=o_{P}(1),\ n\to \infty.$$

Under conditions (i)-(v),  conditions (A1) and (B1)--(B4) assumed in Proposition 2 in  \cite{Dai.18}   hold.
Weak--consistency in the supremum geodesic distance of the empirical Fr\'echet   functional   mean, under weak--dependent $\mathcal{M}$--valued curve  data,  can then be derived  in a similar way to Proposition 2 in \cite{Dai.18}, applying the mean--square ergodicity of the log--mapped regressor process in condition (iv). Furthermore, under conditions (i)-(v),
\begin{eqnarray}&&
\widehat{M}_{n}\left(z(\cdot), x(\cdot )\right)-\widetilde{M}_{n}\left(z(\cdot), x(\cdot )\right)\nonumber\\
&&=\frac{1}{n}\sum_{i=1}^{n}\left[\int_{\mathcal{T}}\left[d^{2}_{\mathcal{M}}\left( Y_{i }(t),z(t)\right)\right]dt\right]
\left[\left\langle \sqrt{\mathcal{K}}\left(\boldsymbol{\gamma}_{x(\cdot ),\overline{X}_{n}(\cdot )}\right), \widehat{\mathcal{R}}_{X}^{-1}\left(\sqrt{\mathcal{K}}\left(\boldsymbol{\gamma}_{X_{i}(\cdot ),\overline{X}_{n}(\cdot )}\right)\right)\right\rangle_{\mathbb{H}}\right.\nonumber\\
&&-\left.\left\langle \sqrt{\mathcal{K}}\left(\boldsymbol{\gamma}_{x(\cdot ), \mu (\cdot )}\right), \mathcal{R}_{X}^{-1}\left(\sqrt{\mathcal{K}}\left(\boldsymbol{\gamma}_{X_{i}(\cdot ),\mu (\cdot )}\right)\right)\right\rangle_{\mathbb{H}}\right]\nonumber\\
&&\leq [\mbox{diam}(\mathcal{M})]^{2}|\mathcal{T}|
\left[\left\langle \sqrt{\mathcal{K}}\left(\boldsymbol{\gamma}_{x(\cdot ),\overline{X}_{n}(\cdot )}\right), \widehat{\mathcal{R}}_{X}^{-1}\left(\sqrt{\mathcal{K}}\left(\frac{1}{n}\sum_{i=1}^{n}\boldsymbol{\gamma}_{X_{i}(\cdot ),\overline{X}_{n}(\cdot )}\right)\right)\right\rangle_{\mathbb{H}}\right.\nonumber\\
&&-\left.\left\langle \sqrt{\mathcal{K}}\left(\boldsymbol{\gamma}_{x(\cdot ), \mu (\cdot )}\right), \mathcal{R}_{X}^{-1}\left(\sqrt{\mathcal{K}}\left(\frac{1}{n}\sum_{i=1}^{n}\boldsymbol{\gamma}_{X_{i}(\cdot ),\mu (\cdot )}\right)\right)\right\rangle_{\mathbb{H}}\right]=o_{P}(1),\nonumber\\
\label{innerprod}
\end{eqnarray}
 \noindent where $|\mathcal{T}|=\int_{\mathcal{T}}dt,$ and
 \begin{eqnarray}&&
 \boldsymbol{\gamma}_{x(\cdot ),\mu (\cdot )}=
\log_{\mu_{X_{0},\mathcal{M}}(\cdot)}\left(x(\cdot )\right)-\mu (\cdot ),\quad x(\cdot )\in \mathcal{X}_{\mathcal{C}_{\mathcal{M}}(\mathcal{T})},\nonumber\\
&&\boldsymbol{\gamma}_{X_{i}(\cdot ),\mu (\cdot )}=\log_{\mu_{X_{0},\mathcal{M}}(\cdot)}\left(X_{i}(\cdot )\right)-\mu(\cdot ),\quad i=1,\dots,n,
 \label{innerprod}
\end{eqnarray}
 \noindent with $\boldsymbol{\gamma}_{x(\cdot ),\overline{X}_{n}(\cdot )}$ and $\boldsymbol{\gamma}_{X_{i}(\cdot ),\overline{X}_{n}(\cdot )}$
 being  introduced in equation  (\ref{elfbb3}).

 Now to prove  (2), under conditions (iv)--(v) in Section  \ref{s31} and (a)--(b) in Section \ref{s32}, applying, in particular, men--square  ergodicity of the  log--mapped regressor process $X$ in condition (iv),  for $n$ sufficiently large, and for $z_{1}(\cdot), z_{2}(\cdot)\in \mathcal{Y}_{\mathcal{C}_{\mathcal{M}}(\mathcal{T})},$
 \begin{eqnarray}&&
 \left|\widehat{M}_{n}(z_{1}(\cdot),x(\cdot))-\widehat{M}_{n}(z_{2}(\cdot),x(\cdot))\right|
 \nonumber\\
 &&\leq \frac{1}{n}\sum_{i=1}^{n}\left|\left\langle \sqrt{\mathcal{K}}\left(\boldsymbol{\gamma}_{x(\cdot ),\overline{X}_{n}(\cdot )}\right), \widehat{\mathcal{R}}_{X}^{-1}\left(\sqrt{\mathcal{K}}\left(\boldsymbol{\gamma}_{X_{i}(\cdot ),\overline{X}_{n}(\cdot )}\right)\right)\right\rangle_{\mathbb{H}}\right|\nonumber\\
 &&\times \left|\int_{\mathcal{T}}\left[d^{2}_{\mathcal{M}}\left( Y_{i }(t),z_{1}(t)\right)-d^{2}_{\mathcal{M}}\left( Y_{i }(t),z_{2}(t)\right)\right]dt
 \right|\nonumber\\
 &&=\frac{1}{n}\sum_{i=1}^{n}\left|\left\langle \sqrt{\mathcal{K}}\left(\boldsymbol{\gamma}_{x(\cdot ),\overline{X}_{n}(\cdot )}\right), \widehat{\mathcal{R}}_{X}^{-1}\left(\sqrt{\mathcal{K}}\left(\boldsymbol{\gamma}_{X_{i}(\cdot ),\overline{X}_{n}(\cdot )}\right)\right)\right\rangle_{\mathbb{H}}\right|\nonumber\\
 &&\times \left|\int_{\mathcal{T}}\left[d_{\mathcal{M}}\left( Y_{i }(t),z_{1}(t)\right)-d_{\mathcal{M}}\left( Y_{i }(t),z_{2}(t)\right)\right]\right.
 \nonumber\\
 &&\hspace*{3cm} \left.\times \left[d_{\mathcal{M}}\left( Y_{i }(t),z_{1}(t)\right)+d_{\mathcal{M}}\left( Y_{i }(t),z_{2}(t)\right)\right]
 dt\right|\nonumber\\
 &&\leq 2\mbox{diam}(\mathcal{M})\sup_{t\in \mathcal{T}} d_{\mathcal{M}}(z_{1}(t),z_{2}(t))|\mathcal{T}|\nonumber\\
 &&\hspace*{1cm}\times \frac{1}{n}\sum_{i=1}^{n} \left|\left\langle \sqrt{\mathcal{K}}\left(\boldsymbol{\gamma}_{x(\cdot ),\overline{X}_{n}(\cdot )}\right), \widehat{\mathcal{R}}_{X}^{-1}\left(\sqrt{\mathcal{K}}\left(\boldsymbol{\gamma}_{X_{i}(\cdot ),\overline{X}_{n}(\cdot )}\right)\right)\right\rangle_{\mathbb{H}}\right|
 \nonumber\\
&& \leq 2\mbox{diam}(\mathcal{M})\sup_{t\in \mathcal{T}} d_{\mathcal{M}}(z_{1}(t),z_{2}(t))|\mathcal{T}|\nonumber\\
 &&\times
 \sup_{i\in \mathbb{Z}}\left\|\boldsymbol{\gamma}_{X_{i}(\cdot ),\overline{X}_{n}(\cdot )}\right\|_{\mathbb{H}}
 \left\|\sqrt{\mathcal{K}}\widehat{\mathcal{R}}_{X}^{-1}\sqrt{\mathcal{K}}\right\|_{\mathcal{L}(\mathbb{H})}
    \left\|\boldsymbol{\gamma}_{x(\cdot ),\overline{X}_{n}(\cdot )}\right\|_{\mathbb{H}}\nonumber\\
 &&=\mathcal{O}_{P}\left(\sup_{t\in \mathcal{T}} d_{\mathcal{M}}(z_{1}(t),z_{2}(t))\right).
 \end{eqnarray}
\noindent   Note that,  under condition (v),  $\sup_{i\in \mathbb{Z}}\left\|\boldsymbol{\gamma}_{X_{i}(\cdot ),\overline{X}_{n}(\cdot )}\right\|_{\mathbb{H}}<\infty$ (see Remark \ref{intervb1}).
  Hence,
 \begin{eqnarray}&&
\sup_{d_{\mathcal{C}_{\mathcal{M}}(\mathcal{T})}(z_{1}(\cdot),z_{2}(\cdot))\leq \delta } \left|\widehat{M}_{n}(z_{1}(\cdot),x(\cdot))-\widehat{M}_{n}(z_{2}(\cdot),x(\cdot))\right|=\mathcal{O}_{P}(\delta ).
  \label{assii}
 \end{eqnarray}

 From Corollary
3.2.3 in \cite{Vaart.1996}, under assumption A.1 equation (\ref{c1w}) holds.  In particular,
the generalized   process $$\left\{Z_{n}(x(\cdot))=\sup_{t\in \mathcal{T}} d_{\mathcal{M}}\left(\widehat{Y}(x(\cdot))(t),\widehat{Y}_{n}(x(\cdot))(t)\right),\ x(\cdot)\in \mathcal{X}_{\mathcal{C}_{\mathcal{M}}(\mathcal{T})}\right\}$$ \noindent    satisfies $Z_{n}(x(\cdot))=o_{P}(1),$ for each $x (\cdot)\in \mathcal{X}_{\mathcal{C}_{\mathcal{M}}(\mathcal{T})}.$

To prove (\ref{c1w2}) under  assumption B.1, since we are evaluating our predictors on the exponential map $\exp_{\mu_{X_{0},\mathcal{M}}(\cdot)}\left(\mathcal{B}_{\mathbb{H}}(0,B)\right)$  of the ball  $\mathcal{B}_{\mathbb{H}}(0,B)$ with center $0$ and radius $B>0$ in $\mathbb{H}$  (see also Remark \ref{intervb1}),
    from Theorem 1.5.4 in \cite{Vaart.1996},  it is sufficient to prove that, for  any $S>0$ and as $\delta \to 0,$

\begin{eqnarray}
\limsup_{n\to \infty} P\left(\sup_{\left\|\log_{\mu_{X_{0},\mathcal{M}}(\cdot)}(x(\cdot)) -\log_{\mu_{X_{0},\mathcal{M}}(\cdot)}(y(\cdot))
 \right\|_{\mathbb{H}}\leq \delta }
\left|Z_{n}(x(\cdot))-Z_{n}(y(\cdot))\right|>2S\right)\to 0.\nonumber\\
\label{equcws}\end{eqnarray}

Hence, since from triangle inequality,   
 \begin{eqnarray}&&\left|Z_{n}(x(\cdot))-Z_{n}(y(\cdot))\right|\leq \sup_{t\in \mathcal{T}} d_{\mathcal{M}}\left(\widehat{Y}(x(\cdot))(t),\widehat{Y}(y(\cdot))(t)\right)\nonumber \\ &&\hspace*{3.5cm}+\sup_{t\in \mathcal{T}} d_{\mathcal{M}}\left(\widehat{Y}_{n}(x(\cdot))(t),\widehat{Y}_{n}(y(\cdot))(t)\right),\  \mbox{a.s},\label{dgep}
\end{eqnarray}
\noindent it is sufficient to prove that  $\widehat{Y}(x(\cdot))$
is uniformly continuous over  $\exp_{\mu_{X_{0},\mathcal{M}}(\cdot)}\left(\mathcal{B}_{\mathbb{H}}(0,B)\right),$ and that
\begin{eqnarray}
&&\hspace*{-1cm}\limsup_{n\to \infty} P\left(\sup_{y(\cdot)\in \widetilde{\mathcal{B}}_{x(\cdot)}(\delta^{\prime}),\ x(\cdot)\in \exp_{\mu_{X_{0},\mathcal{M}}(\cdot)}\left(\mathcal{B}_{\mathbb{H}}(0,B)\right)}
\sup_{t\in \mathcal{T}}  d_{\mathcal{M}}\left(\widehat{Y}_{n}(x(\cdot))(t),\widehat{Y}_{n}(y(\cdot))(t)\right)>S\right)\to 0,\nonumber\\
\label{dgep2bb}
\end{eqnarray}
\noindent keeping in mind that under conditions  (ii) and (v)  in Section \ref{s31},  given a  $\delta >0,$ there exists $\delta^{\prime}\leq \delta,$ such that \begin{eqnarray}&&
\exp_{\mu_{X_{0},\mathcal{M}}(\cdot)}\left(\mathcal{B}_{\mathbb{H}}(\log_{\mu_{X_{0},\mathcal{M}}(\cdot)}(x(\cdot)),\delta )\right)=\widetilde{\mathcal{B}}_{\mathcal{C}_{\mathcal{M}}(\mathcal{T}),
d_{\mathcal{C}_{\mathcal{M}}(\mathcal{T})}}\left(x(\cdot),\delta^{\prime}\right)=\widetilde{\mathcal{B}}_{x(\cdot)}(\delta^{\prime})
\nonumber\\
&&\hspace*{1cm}:=\left\{ y(\cdot)\in \mathcal{X}_{\mathcal{C}_{\mathcal{M}}(\mathcal{T})};\ \left\|\log_{\mu_{X_{0},\mathcal{M}}(\cdot)}(x(\cdot))-
\log_{\mu_{X_{0},\mathcal{M}}(\cdot)}(y(\cdot))
 \right\|_{\mathbb{H}}\leq \delta \right\}. \label{dgep2b}
\end{eqnarray}

Let  $\delta >0,$ and consider
$x(\cdot)\in \exp_{\mu_{X_{0},\mathcal{M}}(\cdot)}\left(\mathcal{B}_{\mathbb{H}}(0,B)\right),$ and $y(\cdot)\in \widetilde{\mathcal{B}}_{x(\cdot)}(\delta^{\prime})$ (see equation (\ref{dgep2b})).
 Under conditions  (i)-(v), and (a) and (b) (see also equation (\ref{emb}) in Section \ref{s32}),   by the form of the loss function $M,$ applying continuity of operator $\sqrt{\mathcal{K}}\mathcal{R}_{X}^{-1} \sqrt{\mathcal{K}},$  as  $\delta \to 0,$
  \begin{eqnarray}&&\sup_{z(\cdot)\in \mathcal{Y}_{\mathcal{C}_{\mathcal{M}}(\mathcal{T})}}\left|M\left(z(\cdot), x(\cdot )\right)-M\left(z(\cdot), y(\cdot )\right)\right|\nonumber\\
  &&\leq \sup_{z(\cdot)\in \mathcal{Y}_{\mathcal{C}_{\mathcal{M}}(\mathcal{T})}}
  E\left[\left(\int_{\mathcal{T}}d^{2}_{\mathcal{M}}\left(Y_{0}(t),z(t)\right)dt\right)\|\boldsymbol{\gamma}_{X_{0}(\cdot ),\mu (\cdot )}\|_{\mathbb{H}}\right]\left\|\sqrt{\mathcal{K}}\mathcal{R}_{X}^{-1}\sqrt{\mathcal{K}}\right\|_{\mathcal{L}(\mathbb{H})}\nonumber\\
  &&\hspace*{0.5cm}\times
  \left\|\boldsymbol{\gamma}_{x(\cdot ),\mu (\cdot )}-\boldsymbol{\gamma}_{y(\cdot ),\mu (\cdot )}\right\|_{\mathbb{H}}
    \leq [\mbox{diam}(\mathcal{M})]^{2}|\mathcal{T}|E\left[\|\boldsymbol{\gamma}_{X_{0}(\cdot ),\mu (\cdot )}\|_{\mathbb{H}}\right]
    \nonumber\\
  &&    \hspace*{0.5cm}\times
   \left\|\sqrt{\mathcal{K}}\mathcal{R}_{X}^{-1}\sqrt{\mathcal{K}}\right\|_{\mathcal{L}(\mathbb{H})}
    \left\|\boldsymbol{\gamma}_{x(\cdot ),\mu (\cdot )}-\boldsymbol{\gamma}_{y(\cdot ),\mu (\cdot )}\right\|_{\mathbb{H}}
        \to 0,
        \label{equicontthlf}\end{eqnarray}
\noindent where $\boldsymbol{\gamma}_{X_{0}(\cdot ),\mu (\cdot )}$ and $\boldsymbol{\gamma}_{x(\cdot ),\mu (\cdot )}$
have been introduced in equation (\ref{innerprod}). The first part of Assumption B.1 then implies that $\widehat{Y}$ is
  continuous over  the ball  \linebreak $\exp_{\mu_{X_{0},\mathcal{M}}(\cdot)}\left(\mathcal{B}_{\mathbb{H}}(0,B)\right)$
 in  the supremum geodesic distance, hence, uniformly continuous.

 To prove  (\ref{dgep2bb}), applying conditions (iv)--(v) (see also Remark \ref{intervb1}),
 we first have to note that   if $d_{\mathcal{C}_{\mathcal{M}}(\mathcal{T})}\left(\widehat{Y}_{n}(x(\cdot)),\widehat{Y}_{n}(y(\cdot))\right)>\varepsilon,$
 under B.1, the following inequality holds in probability
 \begin{eqnarray}&&
 \zeta \leq \sup_{\left\|\log_{\mu_{X_{0},\mathcal{M}}(\cdot)}\left(x(\cdot )\right)-\log_{\mu_{X_{0},\mathcal{M}}(\cdot)}\left( y(\cdot )\right)\right\|<\delta  }\sup_{z(\cdot)\in \mathcal{Y}_{\mathcal{C}_{\mathcal{M}}(\mathcal{T})}}\left|\widehat{M}_{n}(z(\cdot), x(\cdot ))-\widehat{M}_{n}(z(\cdot), y(\cdot ))\right|,
 \nonumber\\
 \label{equb}\end{eqnarray}
 \noindent for certain $\zeta>0.$ Furthermore, from the expression of $\widehat{M}_{n}$ in equation (\ref{elfbb}),
 for $x(\cdot ), y(\cdot),$ such that  $\left\|\log_{\mu_{X_{0},\mathcal{M}}(\cdot)}\left(x(\cdot )\right)-\log_{\mu_{X_{0},\mathcal{M}}(\cdot)}\left( y(\cdot )\right)\right\|<\delta ,$ for $n$ sufficiently large,
  \begin{eqnarray}&&
  \sup_{z(\cdot)\in \mathcal{Y}_{\mathcal{C}_{\mathcal{M}}(\mathcal{T})}}\left|\widehat{M}_{n}(z(\cdot), x(\cdot ))-\widehat{M}_{n}(z(\cdot), y(\cdot ))\right|\leq [\mbox{diam}(\mathcal{M})]^{2}|\mathcal{T}|\sup_{i\in \mathbb{Z}}\left\|\boldsymbol{\gamma}_{X_{i}(\cdot ),\overline{X}_{n}(\cdot )}\right\|_{\mathbb{H}}
  \nonumber\\
    && \hspace*{2cm}\times  \left\|\sqrt{\mathcal{K}}\widehat{\mathcal{R}}_{X}^{-1}\sqrt{\mathcal{K}}\right\|_{\mathcal{L}(\mathbb{H})}
    \left\|\boldsymbol{\gamma}_{x(\cdot ),\overline{X}_{n}(\cdot )}-\boldsymbol{\gamma}_{y(\cdot ),\overline{X}_{n}(\cdot )}\right\|_{\mathbb{H}}=\mathcal{O}_{P}(\delta ).
     \label{eqbound}
 \end{eqnarray}
\noindent Equation   (\ref{dgep2bb}) then follows from equations (\ref{equb}) and (\ref{eqbound}),
when $\delta\to 0,$ keeping in mind equation (\ref{dgep2b}), and the second part of condition B.1
(see again Corollary
3.2.3 in \cite{Vaart.1996}).

\end{proof}

\section{Analysis of the remaining months December, 1979--May, 1980}
\label{app2}

\begin{figure}[!h]
\begin{center}
\includegraphics[width=0.32\textwidth]{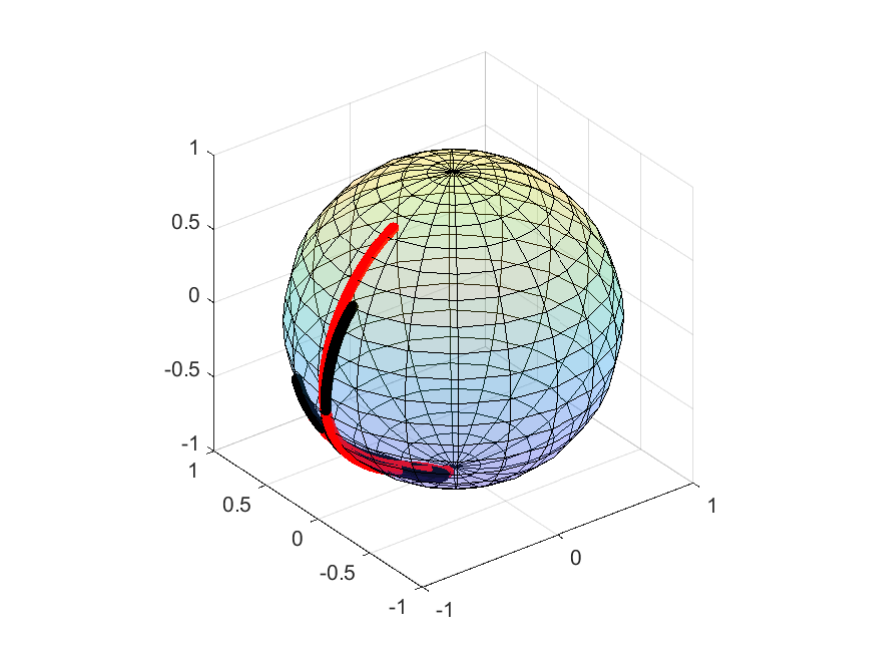}
\includegraphics[width=0.32\textwidth]{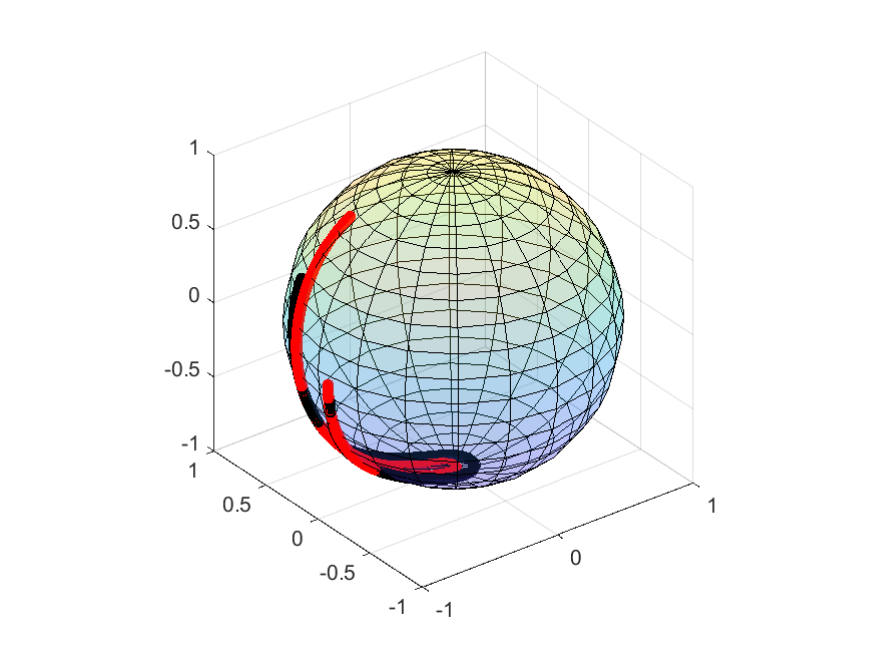}
\includegraphics[width=0.32\textwidth]{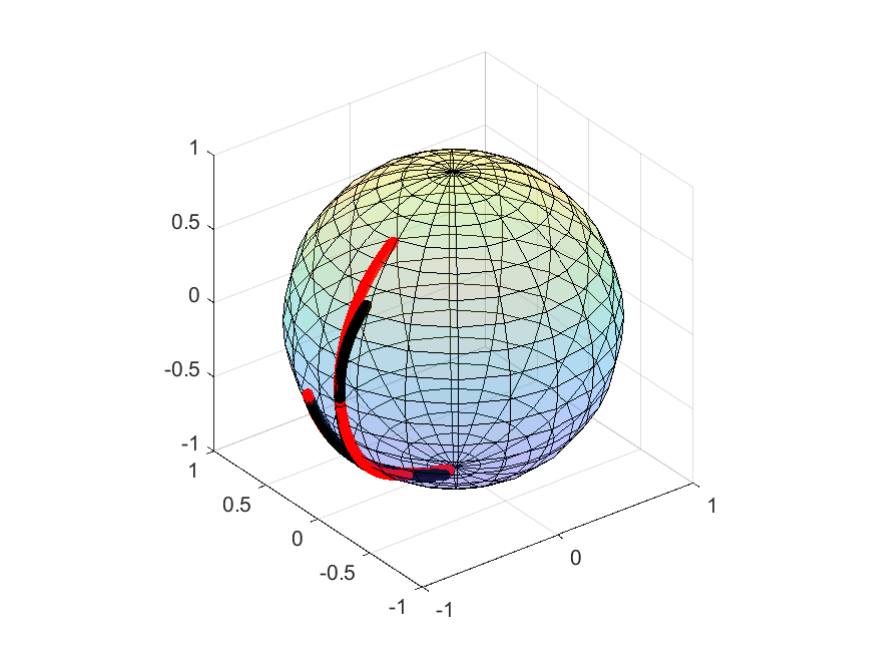}
\includegraphics[width=0.32\textwidth]{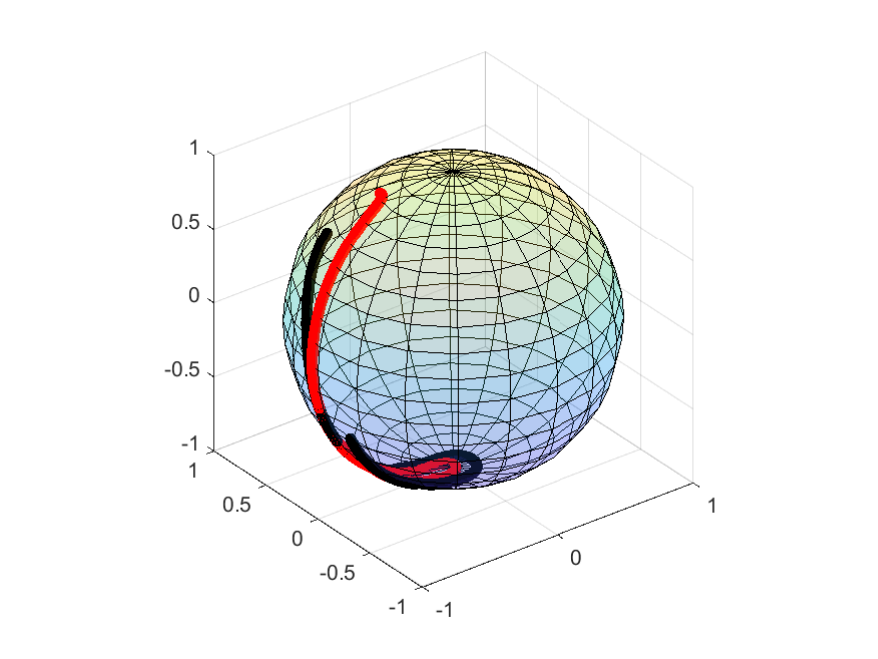}
\includegraphics[width=0.32\textwidth]{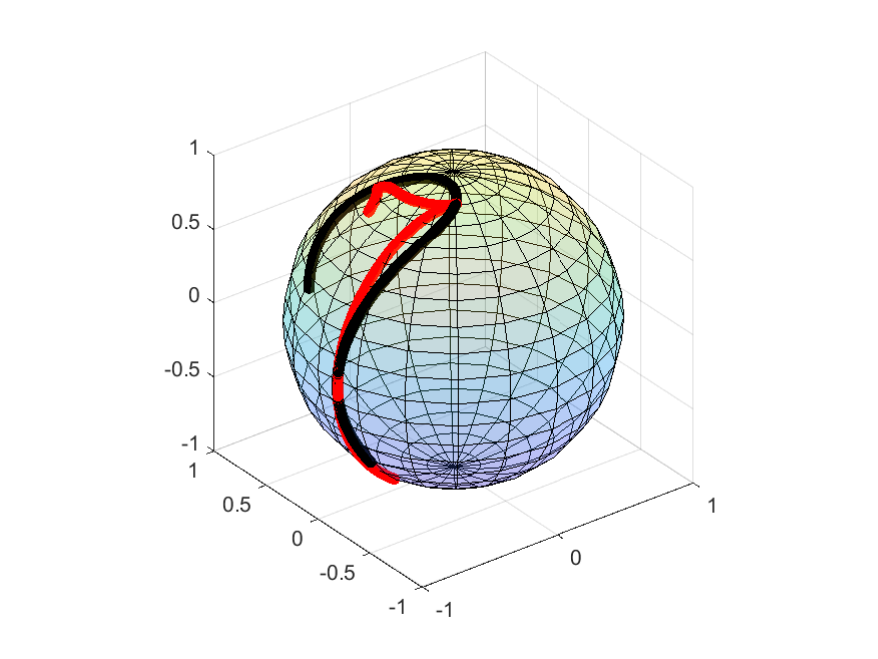}
\includegraphics[width=0.32\textwidth]{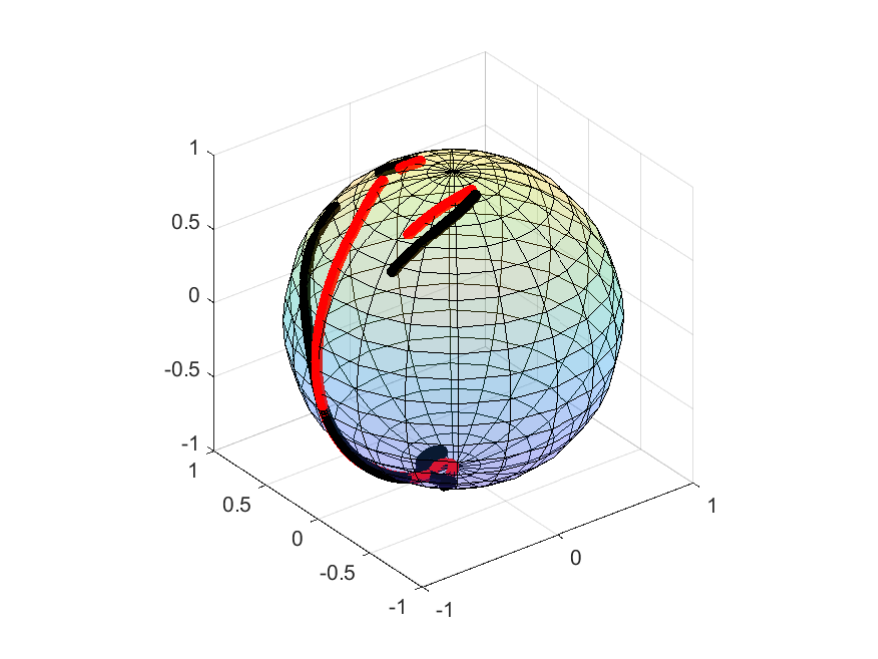}
\includegraphics[width=0.32\textwidth]{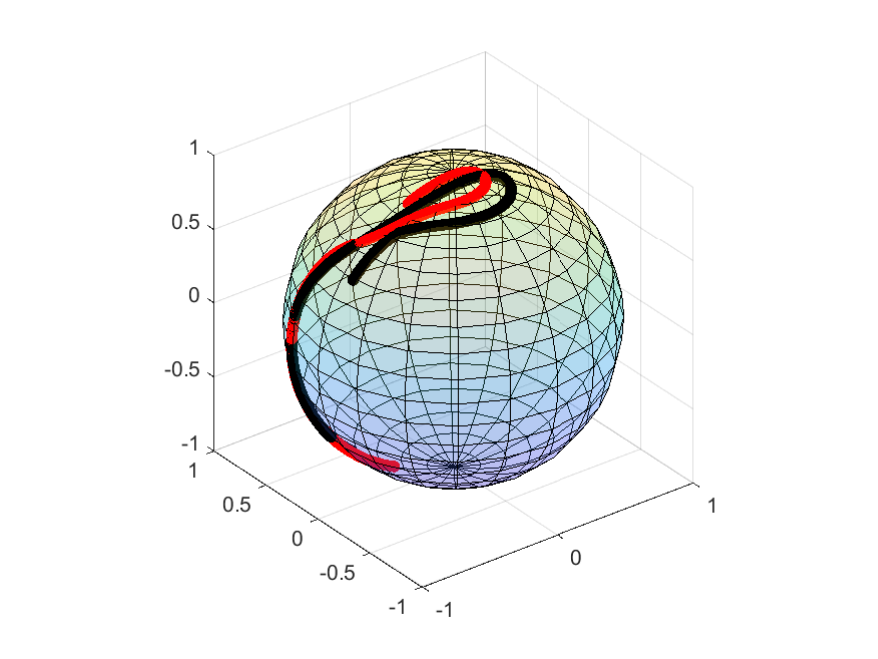}
\includegraphics[width=0.32\textwidth]{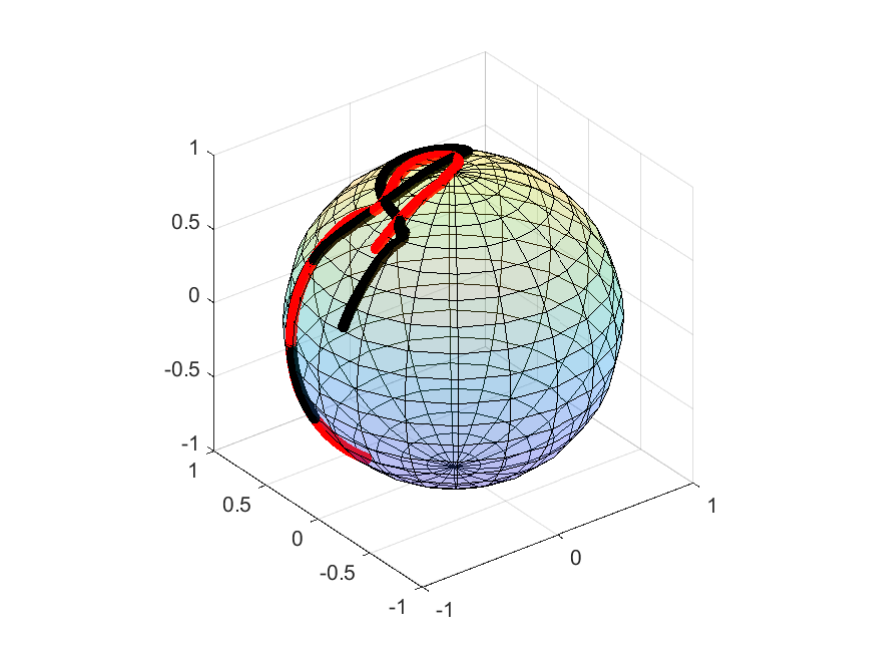}
\includegraphics[width=0.32\textwidth]{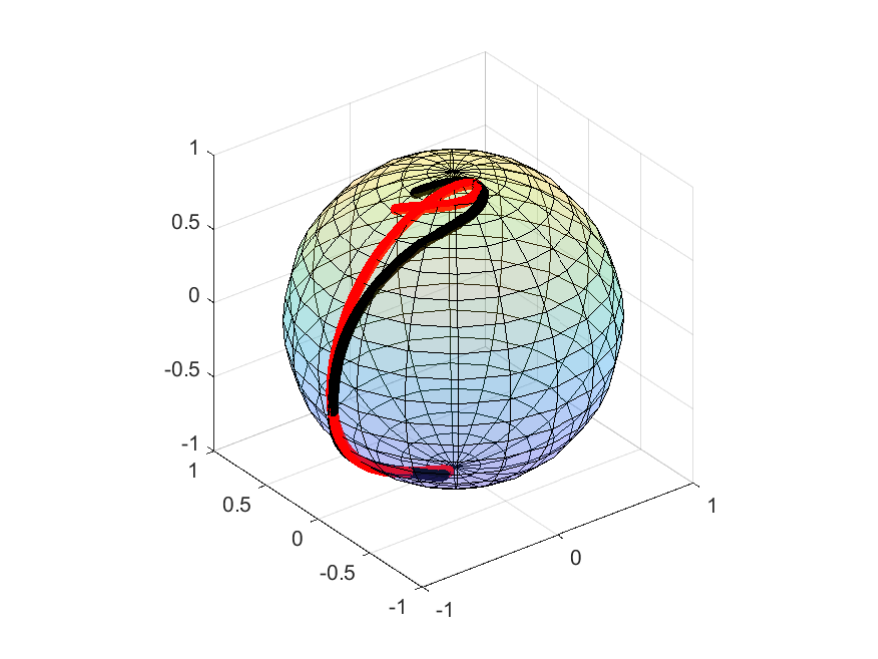}
\includegraphics[width=0.32\textwidth]{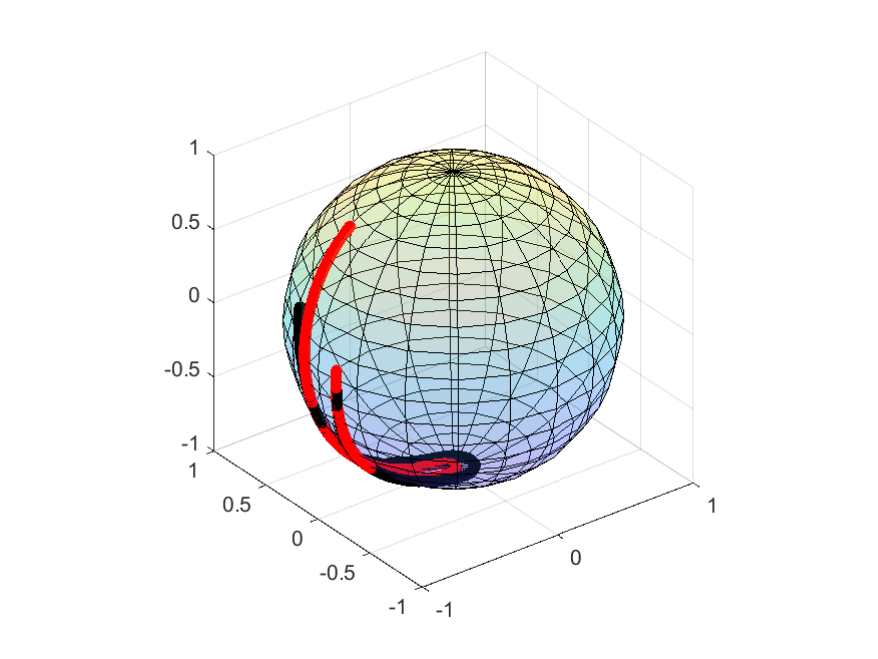}
\includegraphics[width=0.32\textwidth]{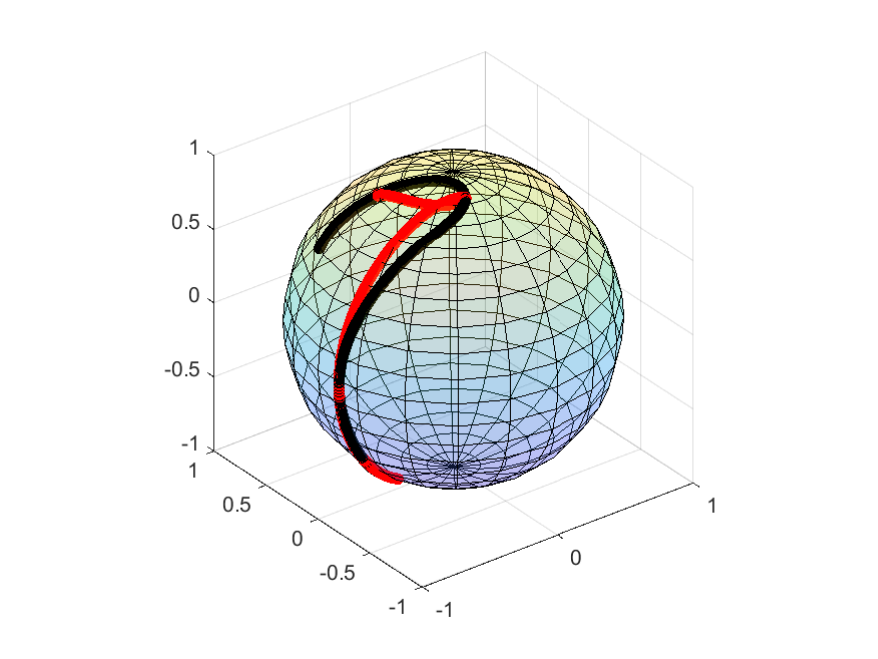}
\includegraphics[width=0.32\textwidth]{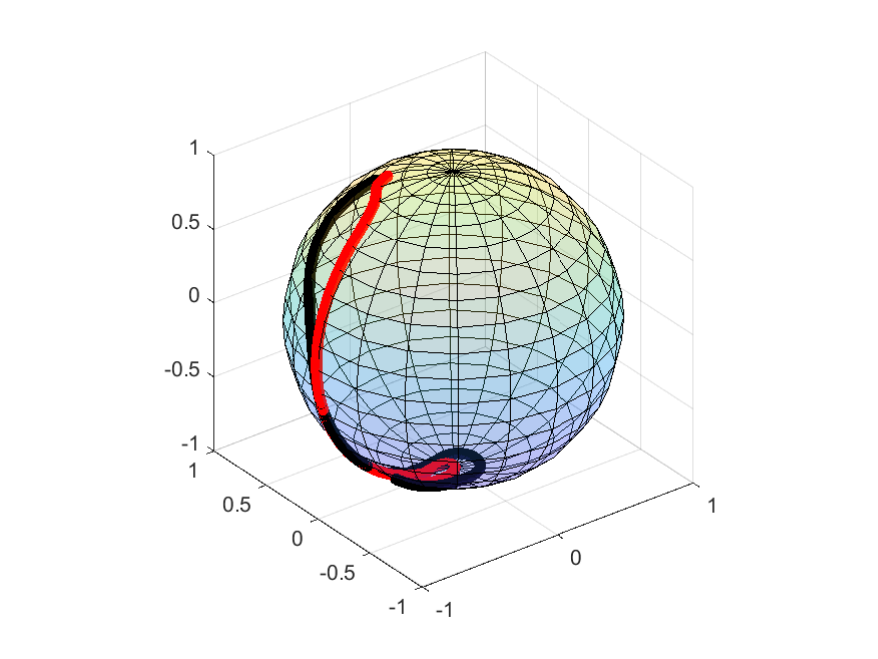}
\includegraphics[width=0.32\textwidth]{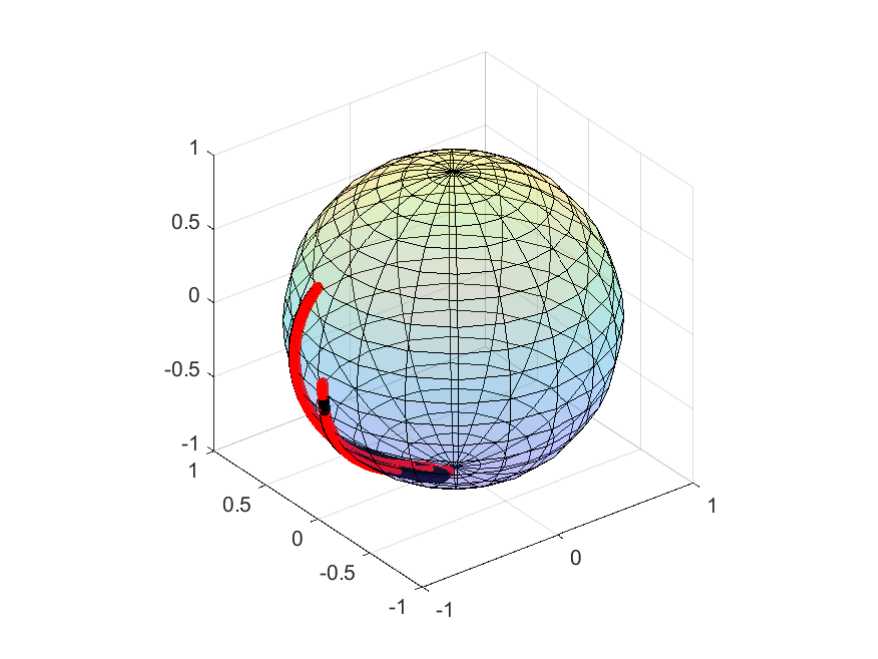}
\includegraphics[width=0.32\textwidth]{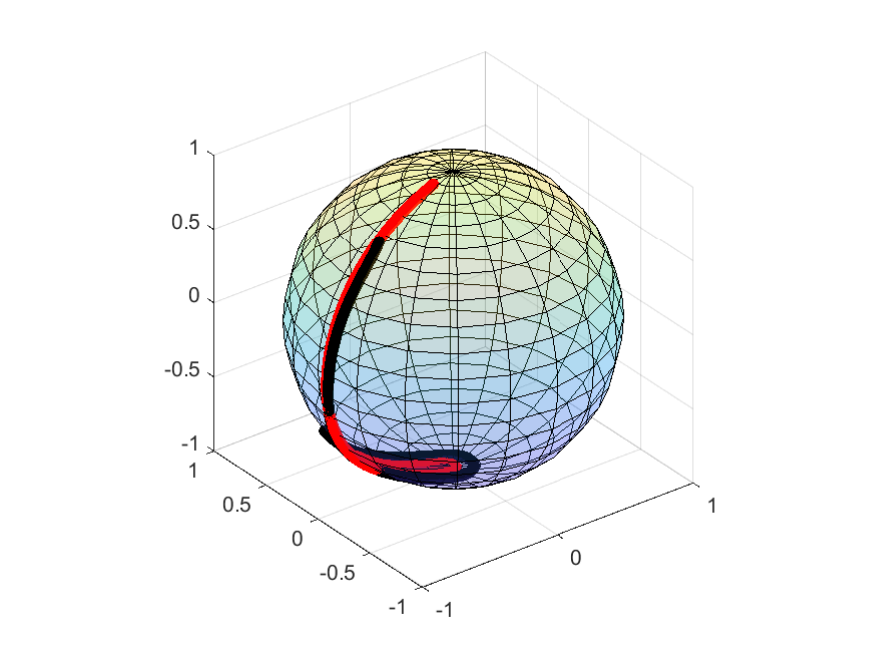}
\includegraphics[width=0.32\textwidth]{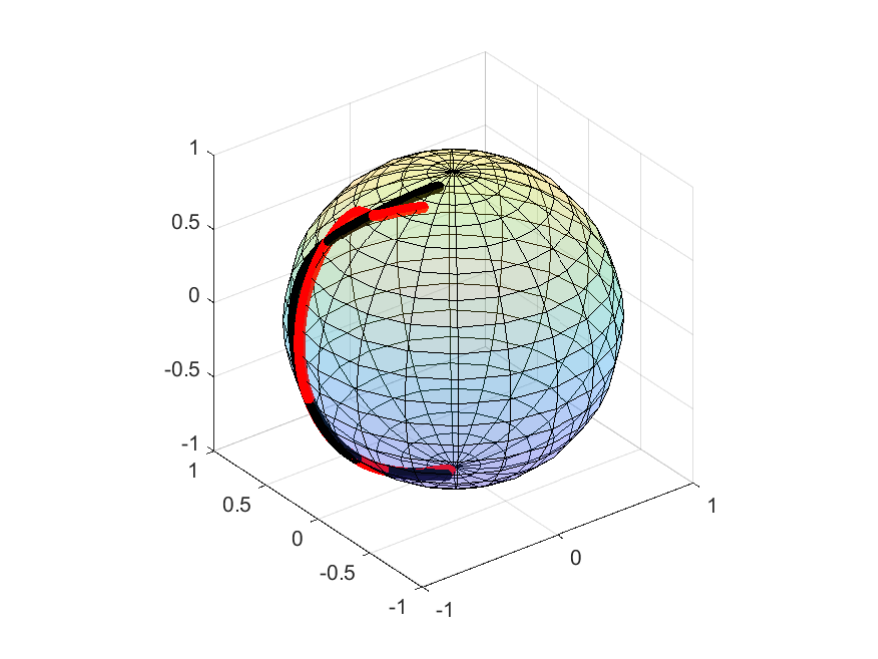}
\includegraphics[width=0.32\textwidth]{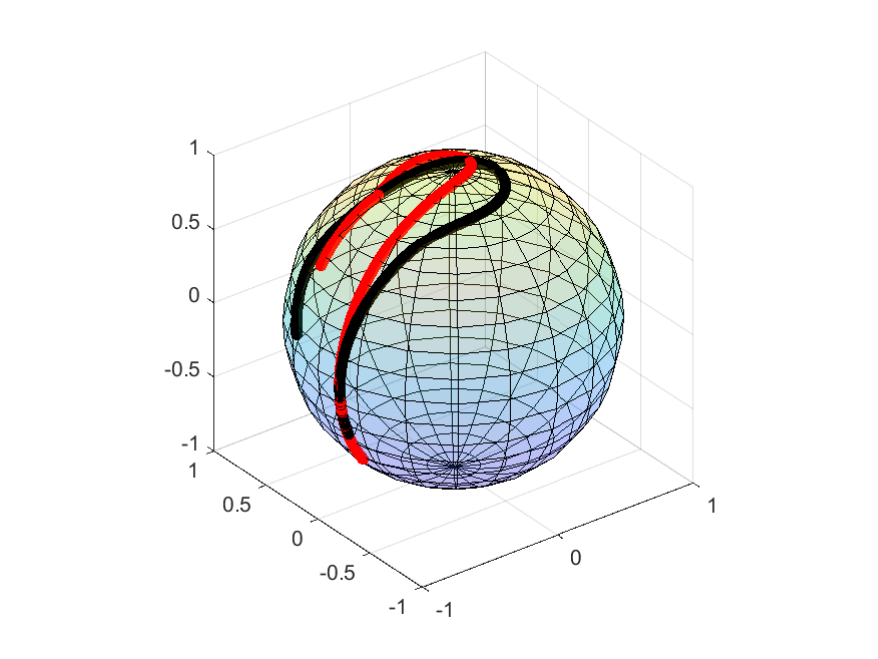}
\includegraphics[width=0.32\textwidth]{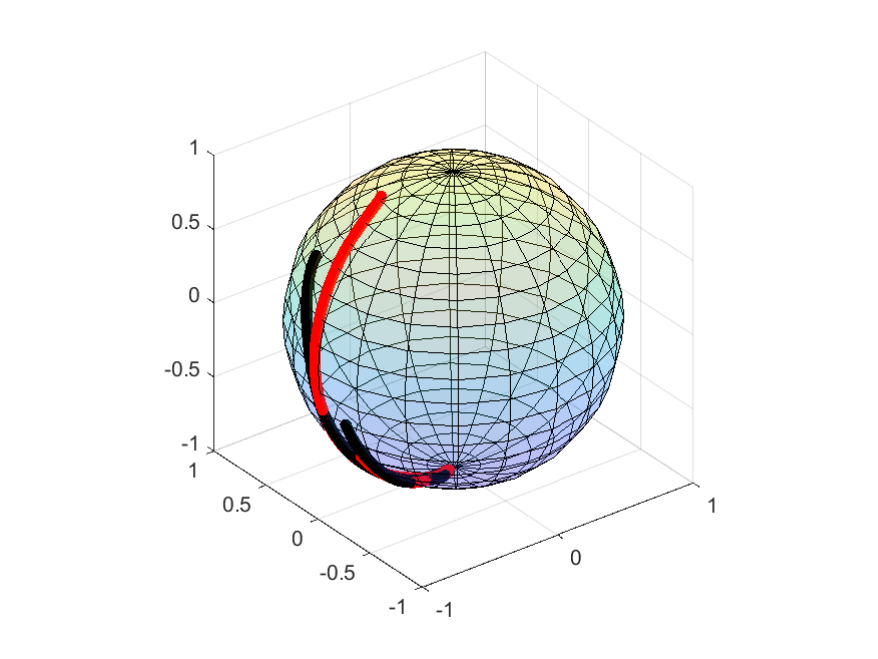}
\includegraphics[width=0.32\textwidth]{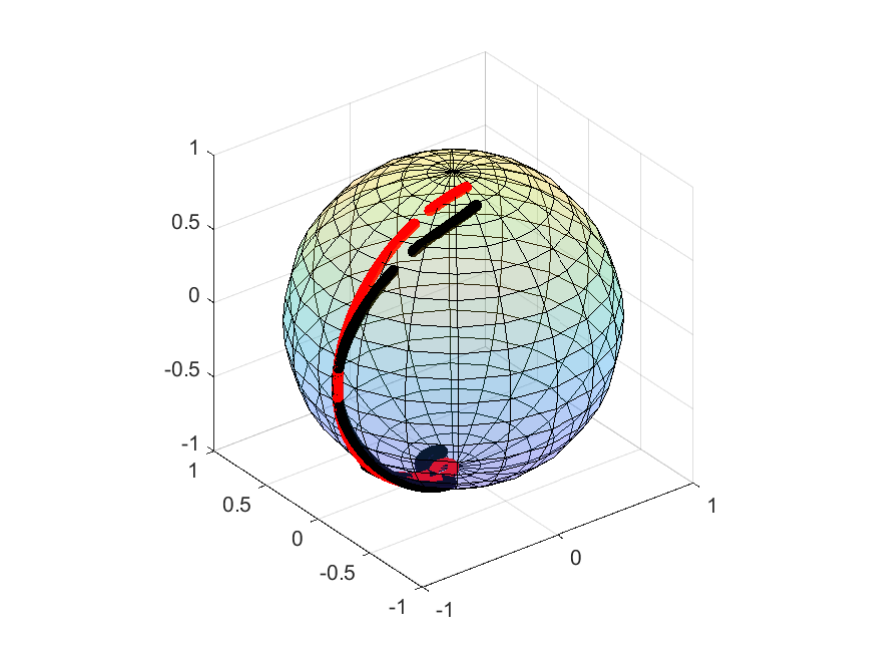}
\end{center}
\caption{Spherical bivariate curve data at times $t=1,15, 29, 43, 57,  71$ for December (first two lines), January (second two lines), and  February  (last two lines), representing  NASA's MAGSAT spacecraft (black curve), and   magnetic field vector (red curve). The spherical coordinates are displayed at  $6000$ temporal nodes for every sampled time.}
\label{fig:1rdam1}
\end{figure}

\begin{figure}[!h]
\begin{center}
\includegraphics[width=0.32\textwidth]{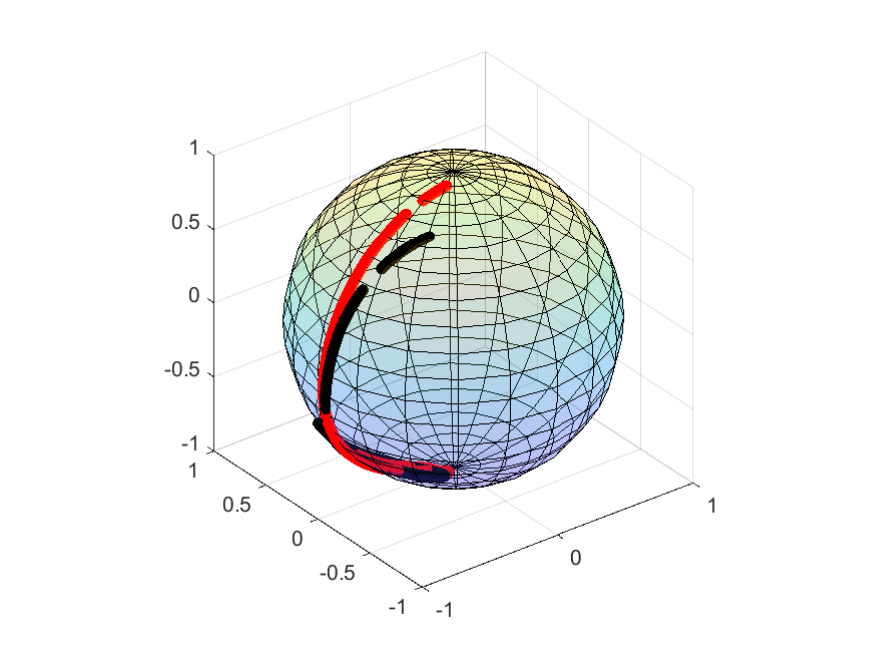}
\includegraphics[width=0.32\textwidth]{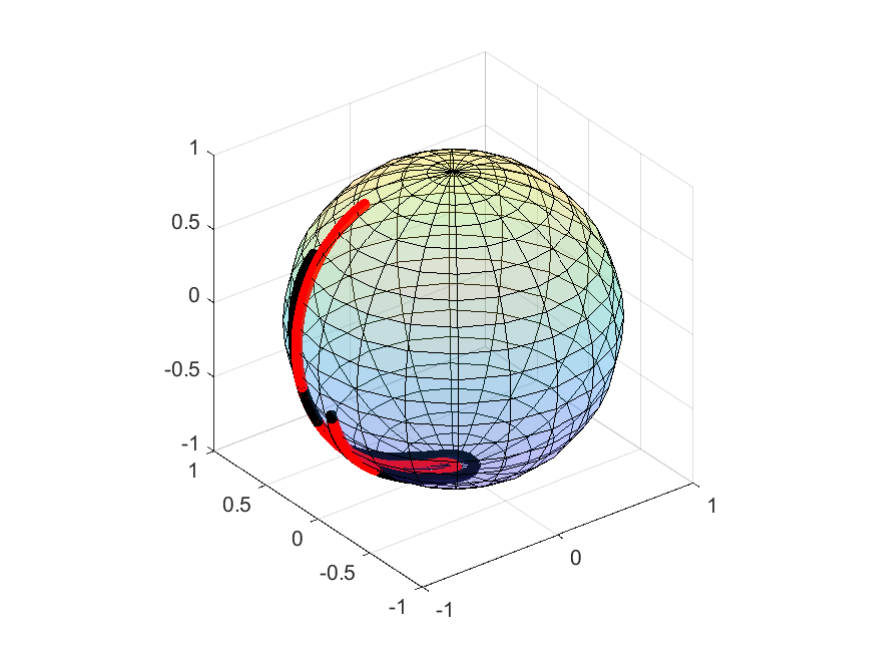}
\includegraphics[width=0.32\textwidth]{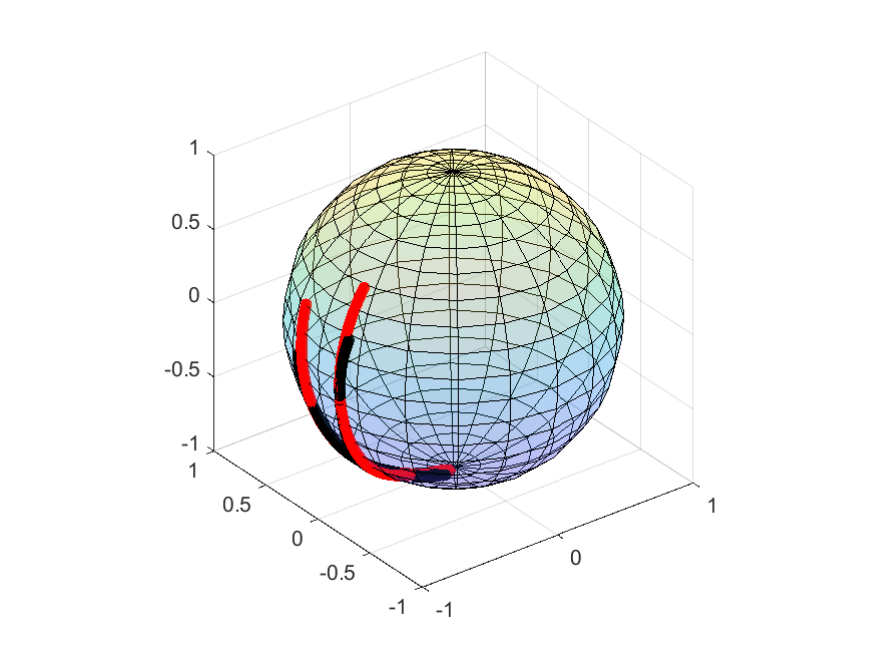}
\includegraphics[width=0.32\textwidth]{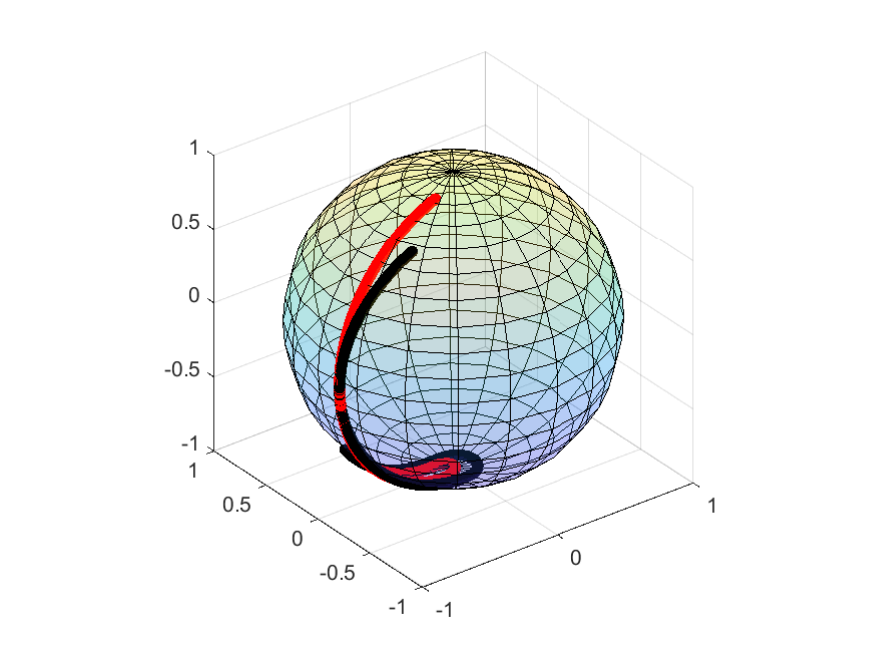}
\includegraphics[width=0.32\textwidth]{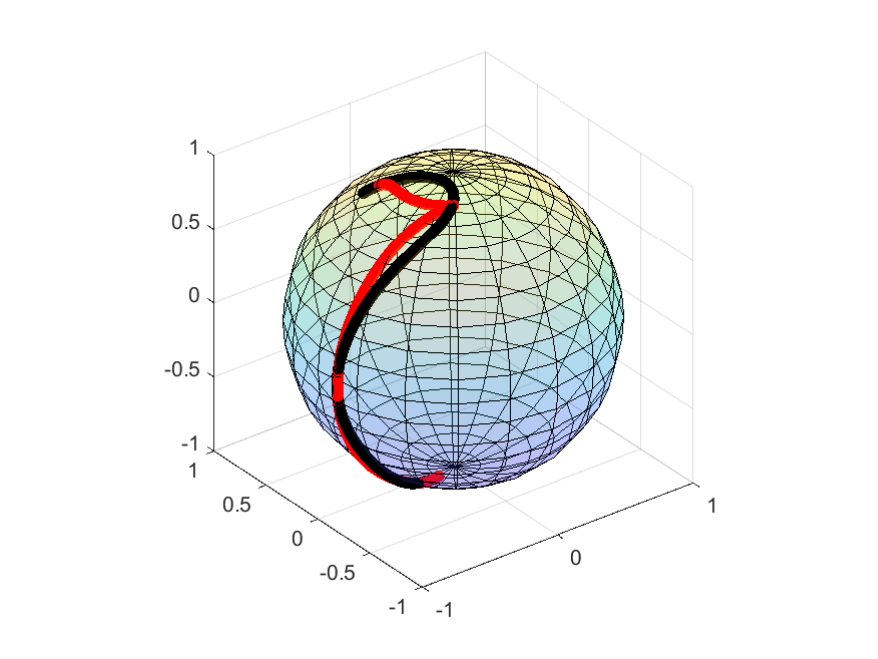}
\includegraphics[width=0.32\textwidth]{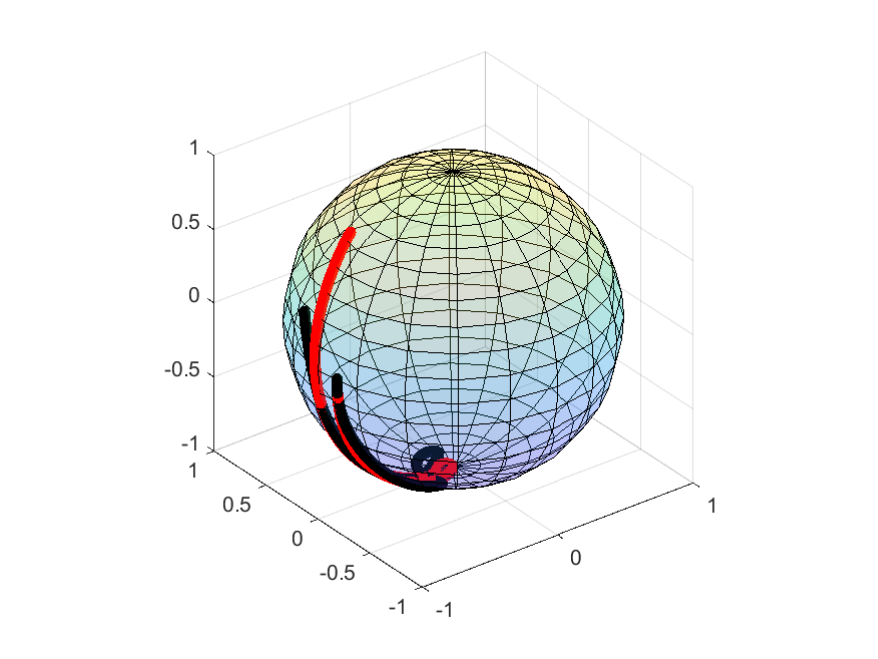}
\includegraphics[width=0.32\textwidth]{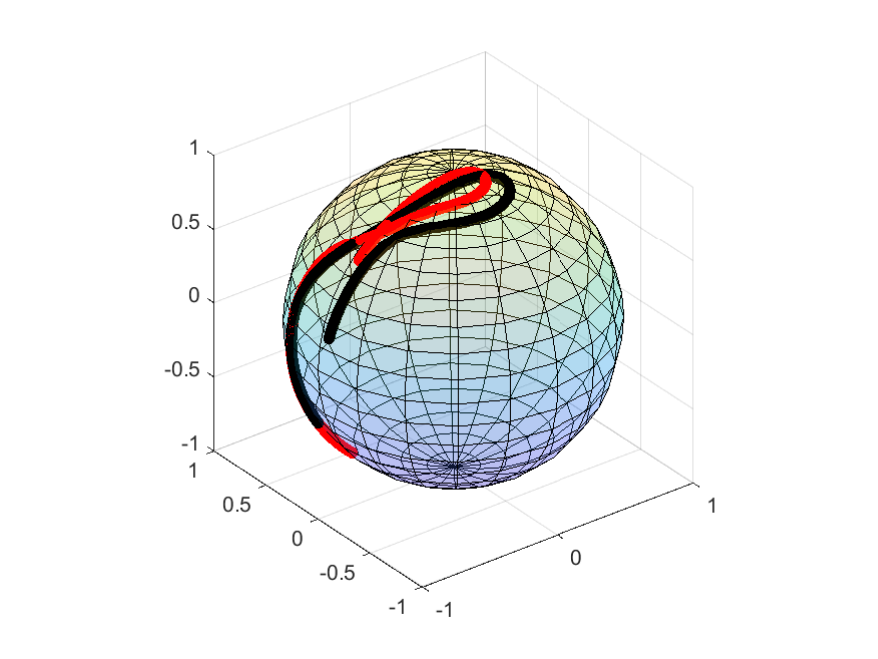}
\includegraphics[width=0.32\textwidth]{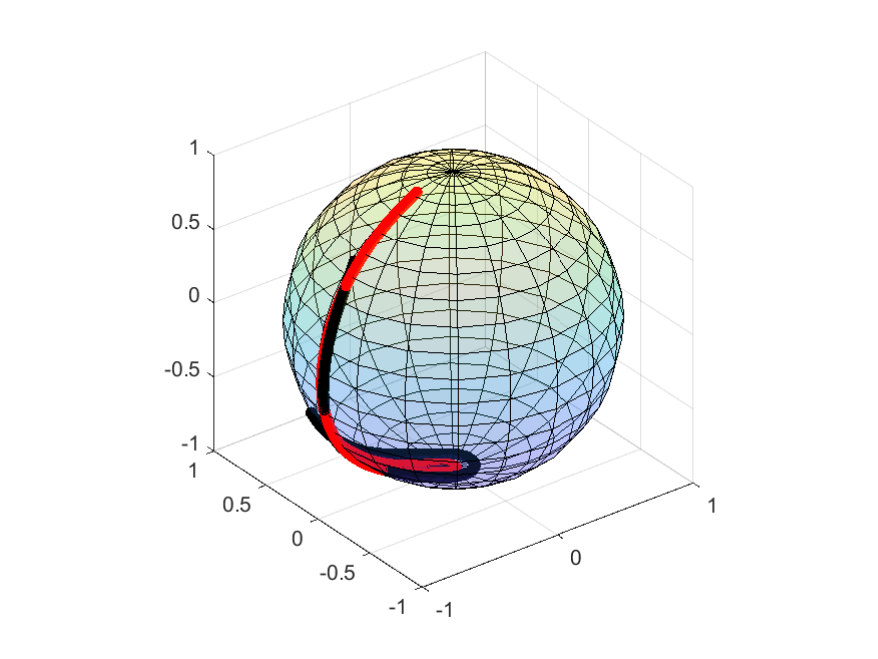}
\includegraphics[width=0.32\textwidth]{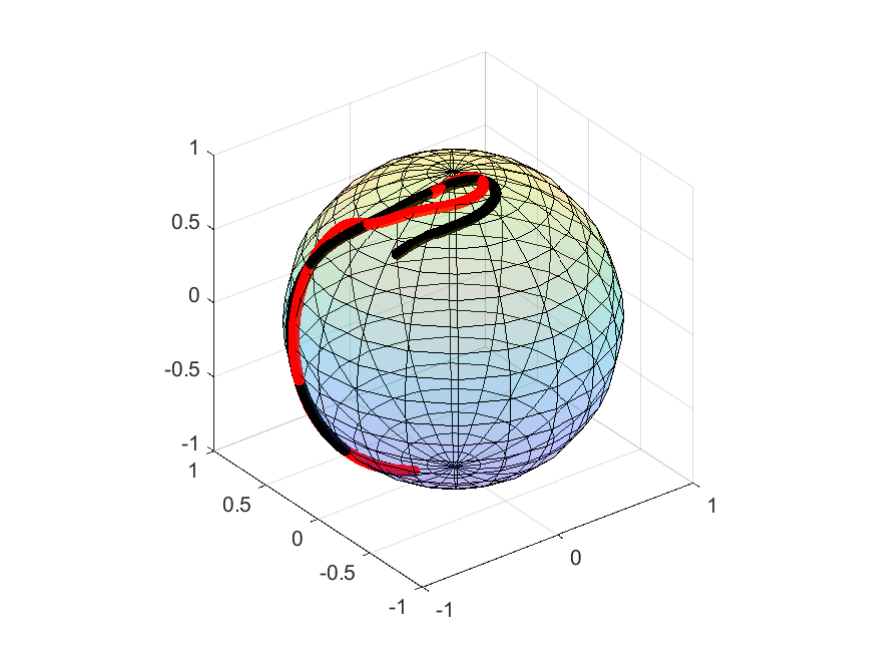}
\includegraphics[width=0.32\textwidth]{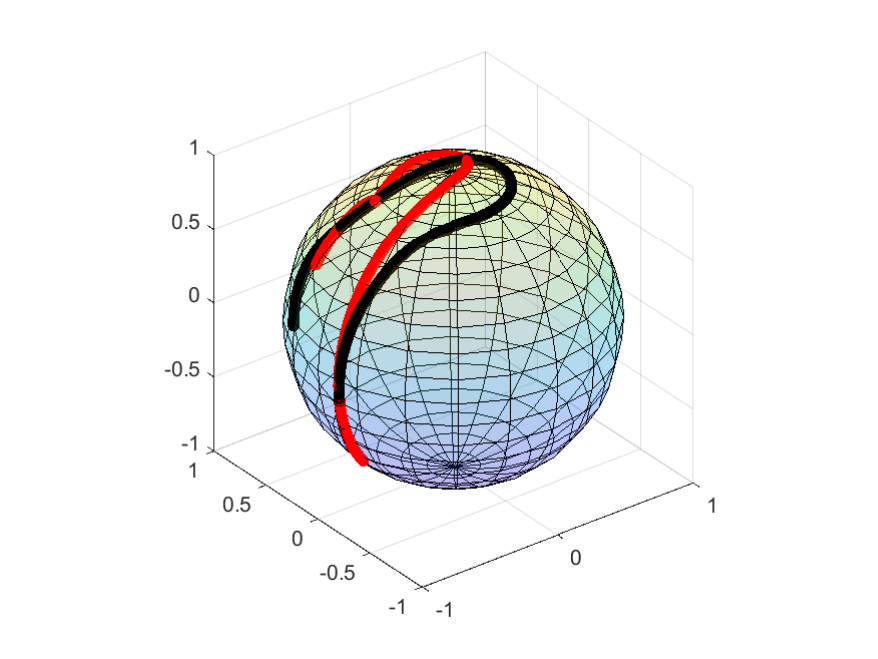}
\includegraphics[width=0.32\textwidth]{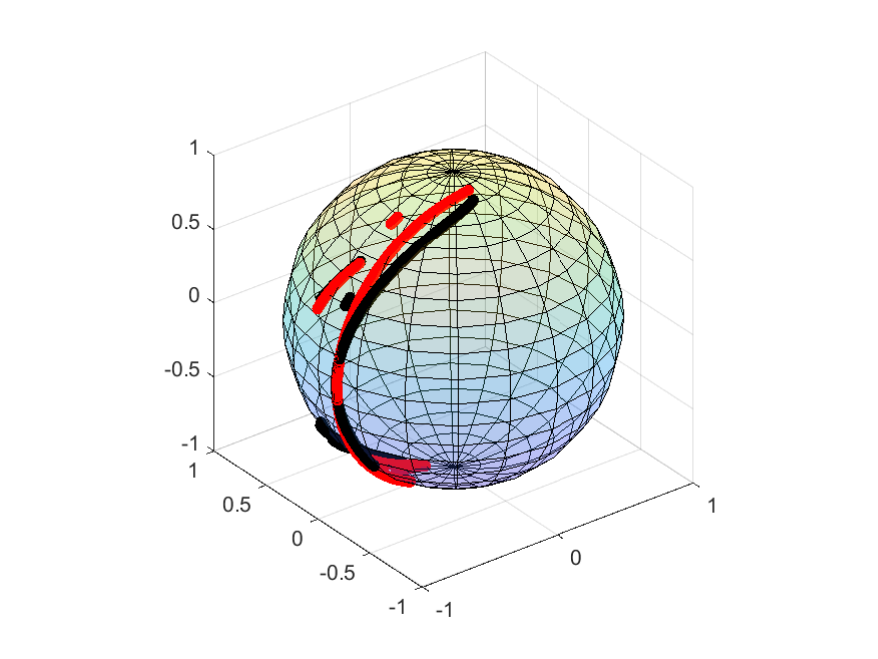}
\includegraphics[width=0.32\textwidth]{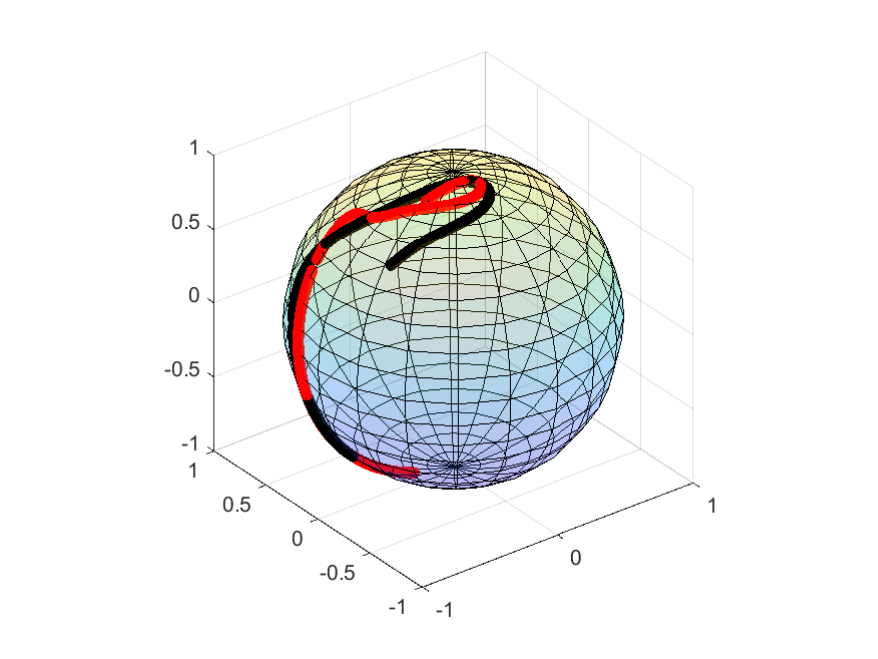}
\includegraphics[width=0.32\textwidth]{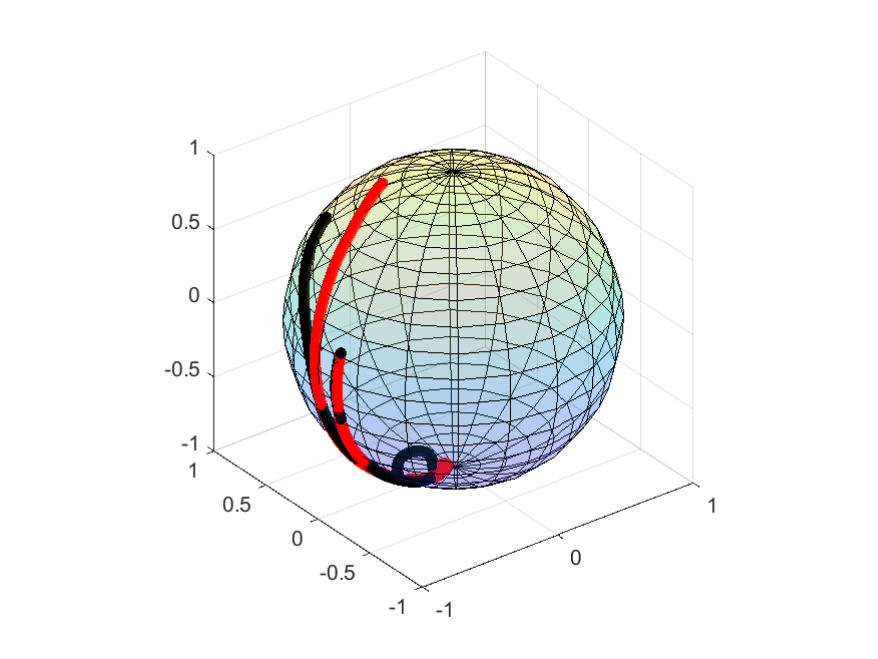}
\includegraphics[width=0.32\textwidth]{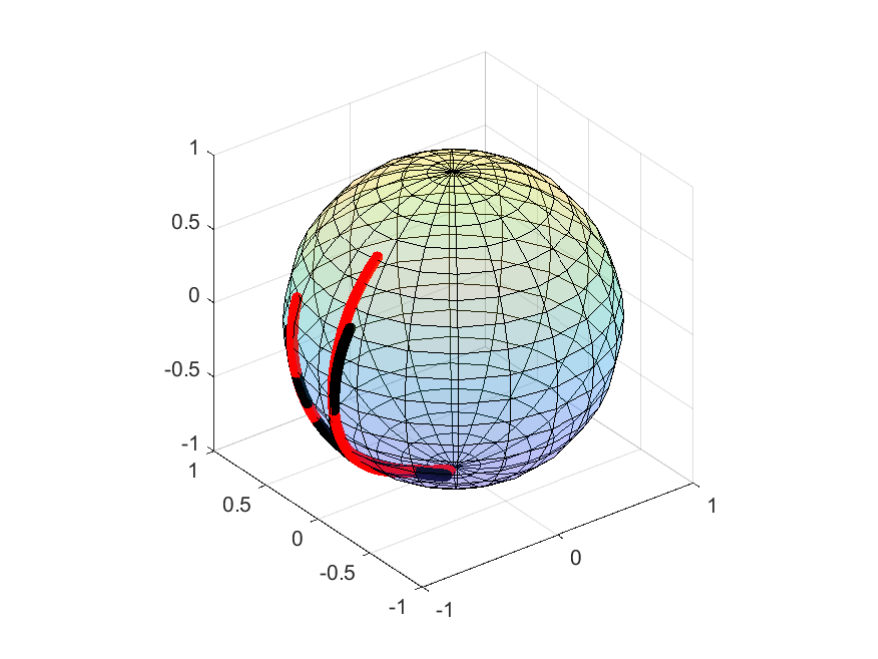}
\includegraphics[width=0.32\textwidth]{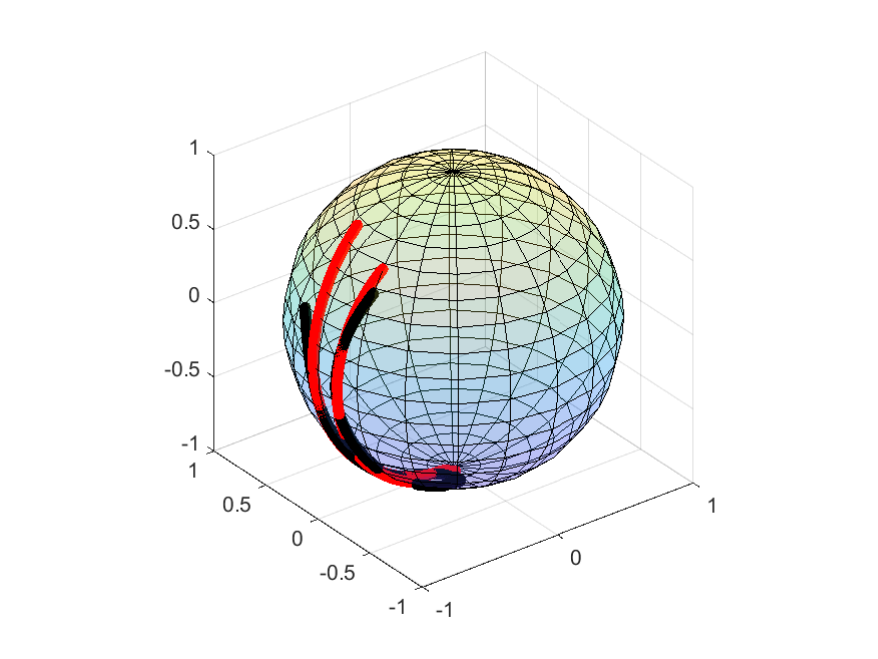}
\includegraphics[width=0.32\textwidth]{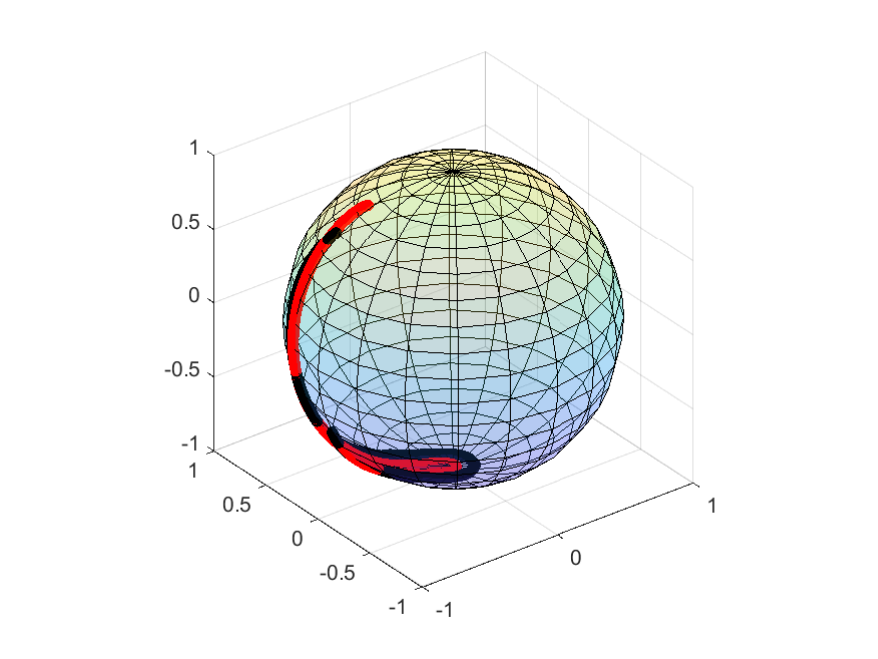}
\includegraphics[width=0.32\textwidth]{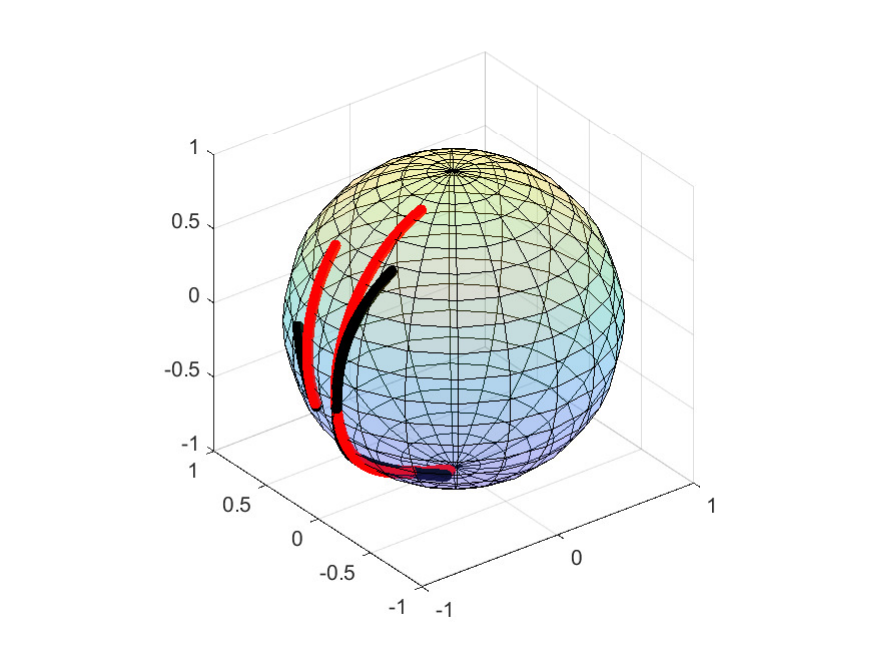}
\includegraphics[width=0.32\textwidth]{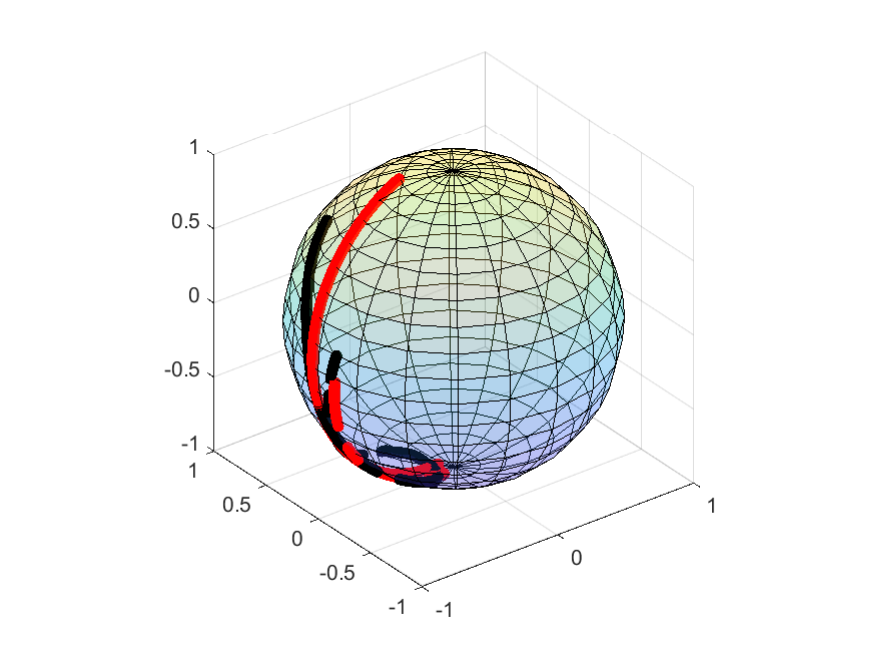}
\end{center}
\caption{Spherical bivariate curve data at times $t=1,15, 29, 43, 57,  71$ for March (first two lines) and April  (second two lines), and at times $t= 1, 15, 29, 43, 48, 52$ for May  (last two lines), representing NASA's MAGSAT spacecraft (black curve), and   magnetic field vector (red curve). The spherical coordinates are displayed at  $6000$ temporal nodes for every sampled time.}
\label{fig:1rdam2}
\end{figure}

\begin{figure}[!h]
\begin{center}
\includegraphics[width=0.32\textwidth]{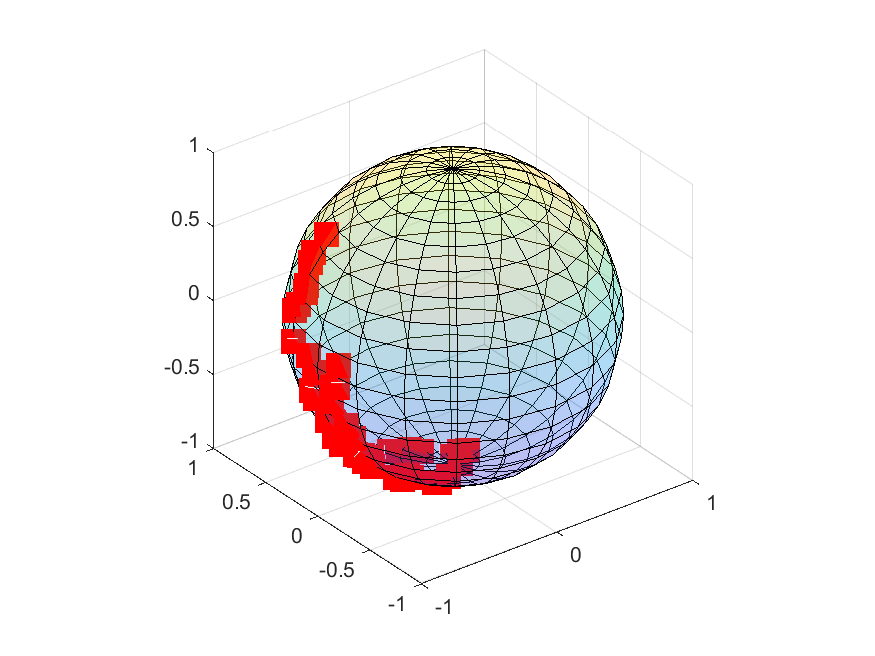}
\includegraphics[width=0.32\textwidth]{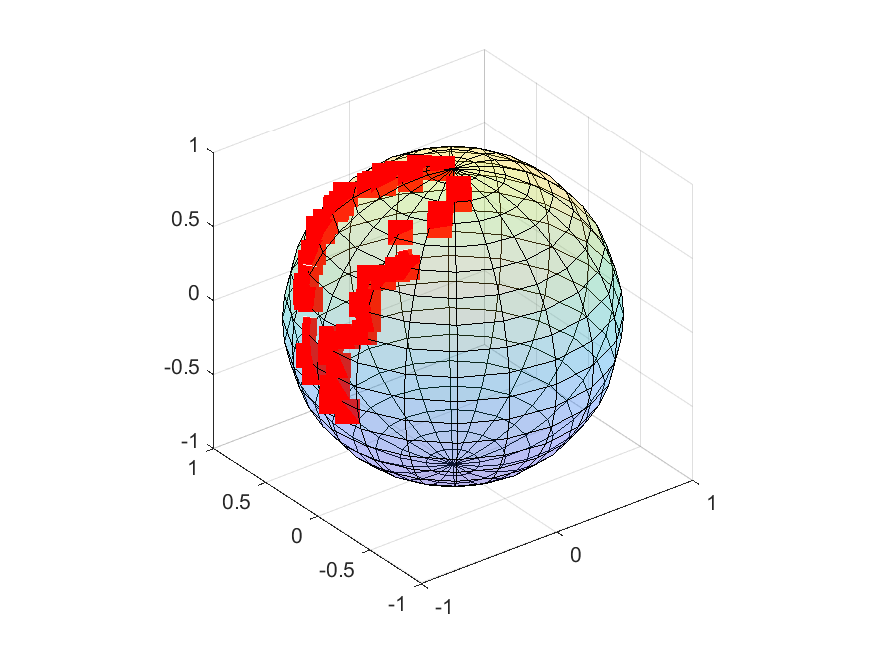}
\includegraphics[width=0.32\textwidth]{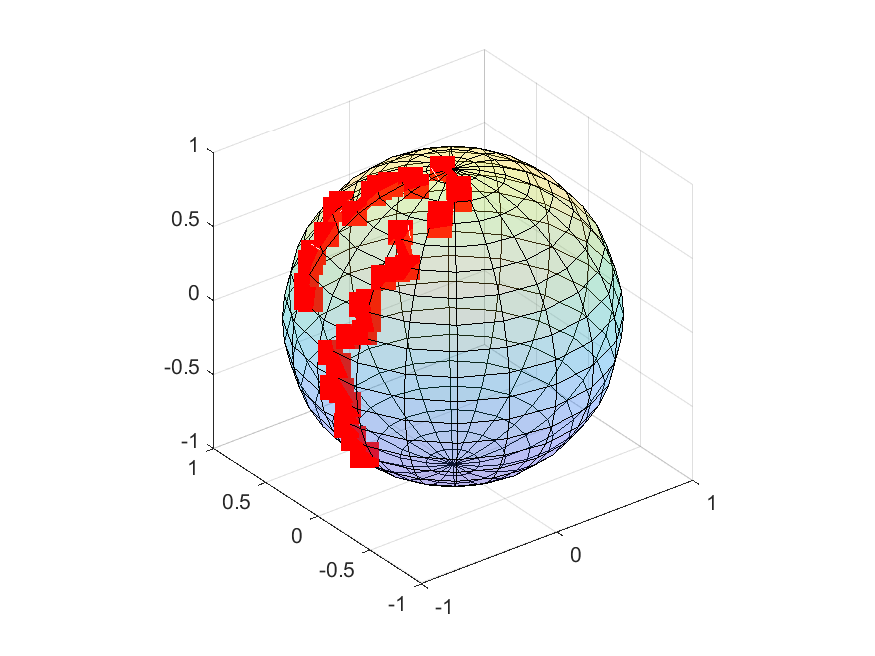}
\includegraphics[width=0.32\textwidth]{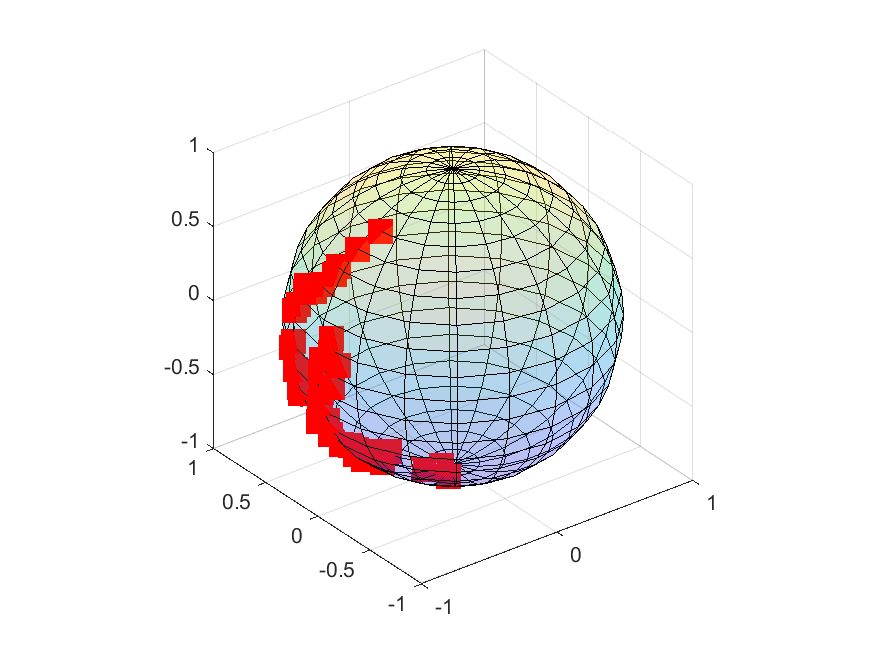}
\includegraphics[width=0.32\textwidth]{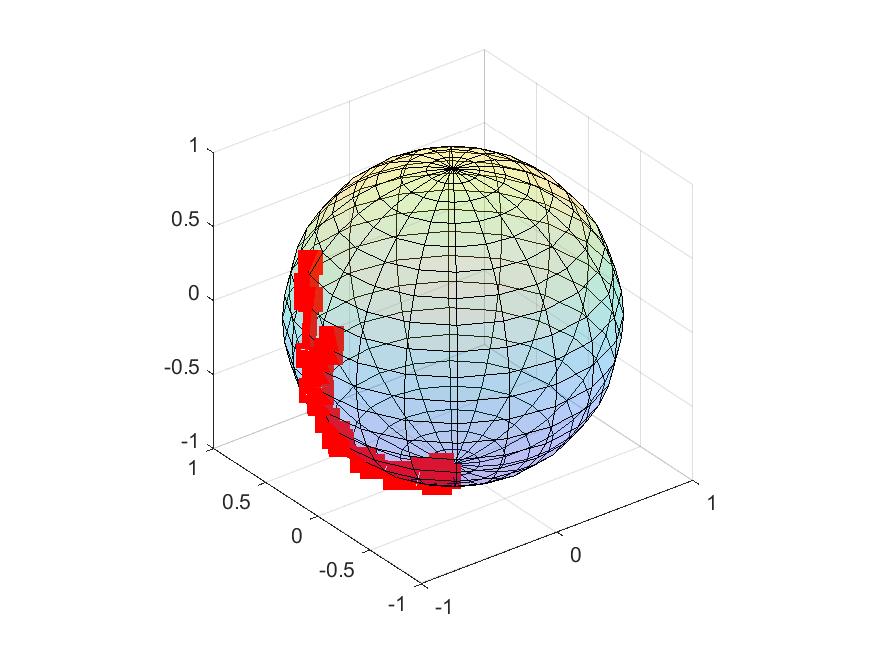}
\includegraphics[width=0.32\textwidth]{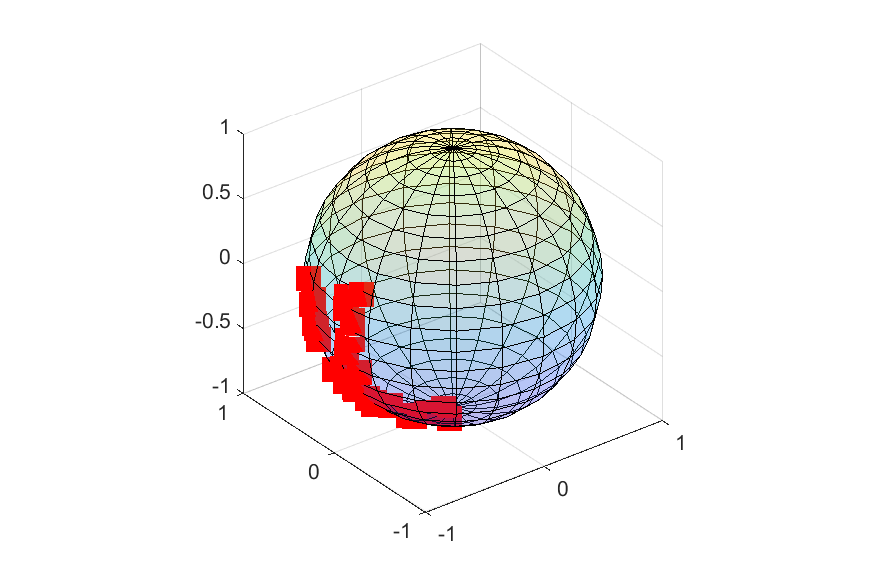}
\end{center}
\caption{Empirical intrinsic  Fr\'echet curve mean of regressors (December, 1979--February, 1980  at the top, and  March--May, 1980 at the bottom).}
\label{fig:1rdam3}
\end{figure}

\begin{figure}[!h]
\begin{center}
\includegraphics[width=0.32\textwidth]{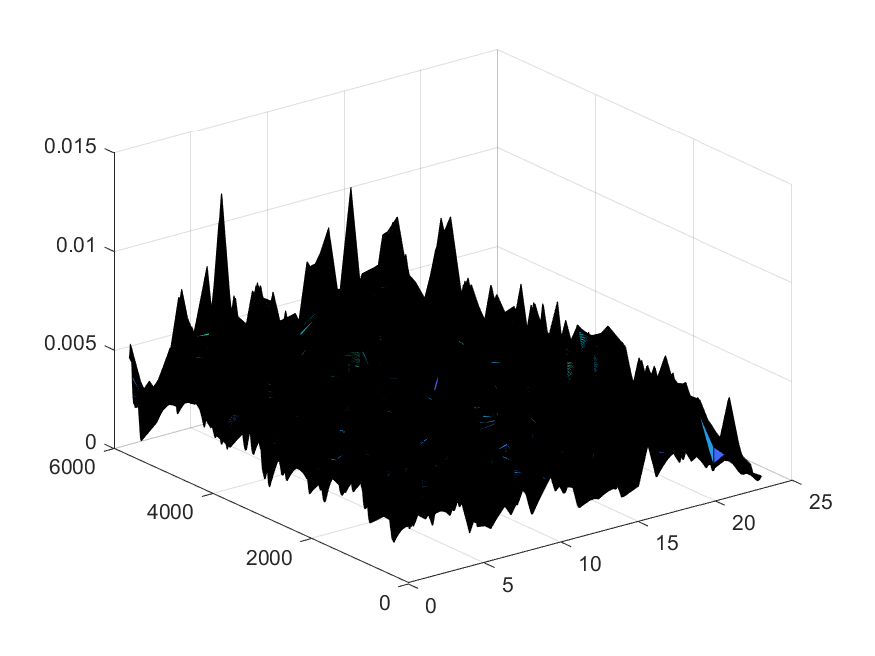}
\includegraphics[width=0.32\textwidth]{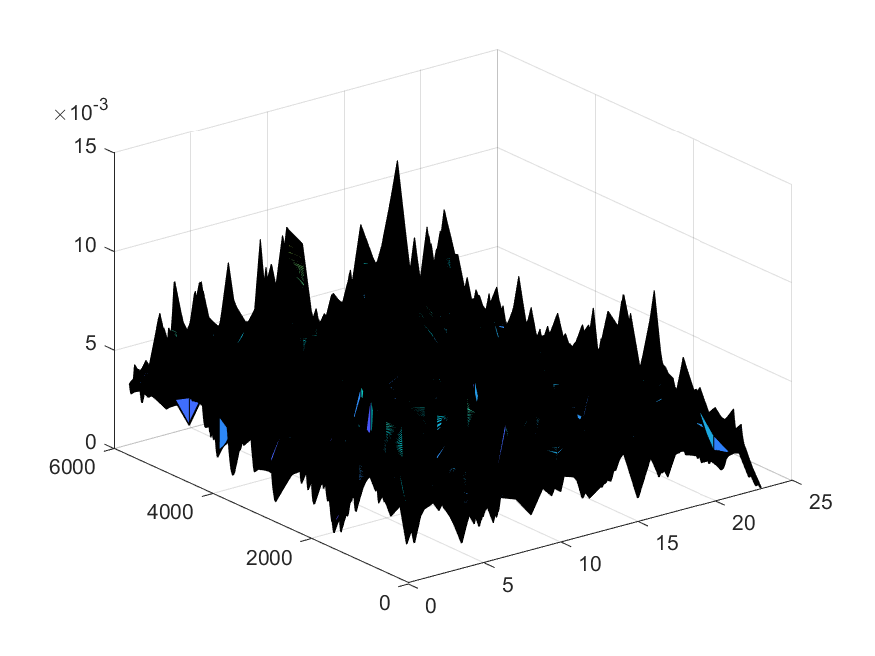}
\includegraphics[width=0.32\textwidth]{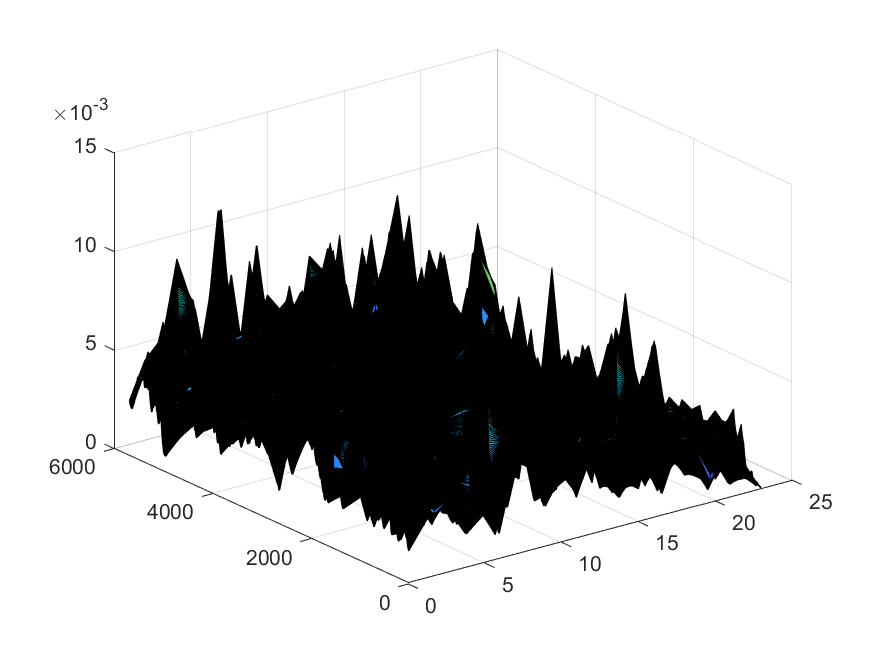}
\includegraphics[width=0.32\textwidth]{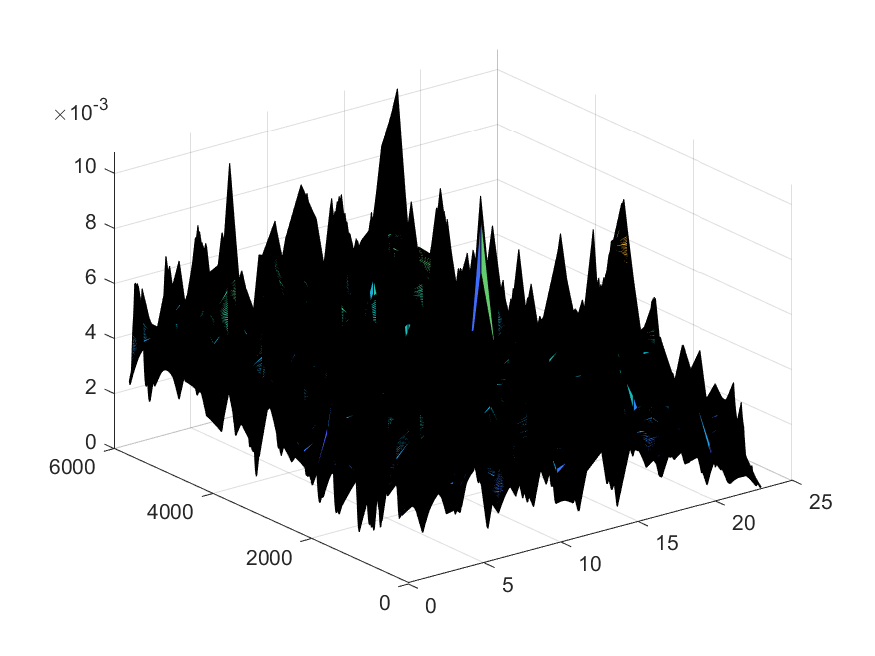}
\includegraphics[width=0.32\textwidth]{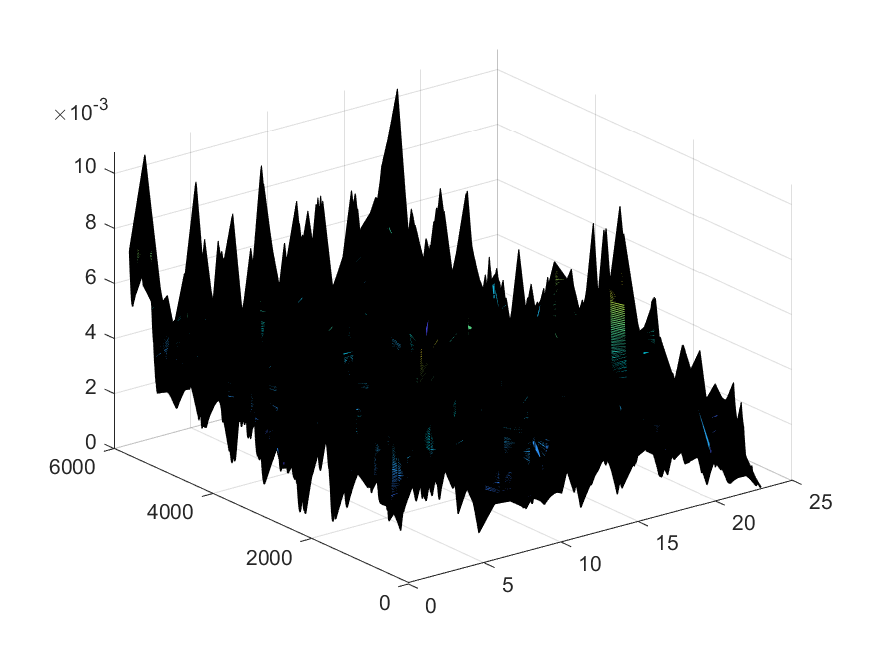}
\includegraphics[width=0.32\textwidth]{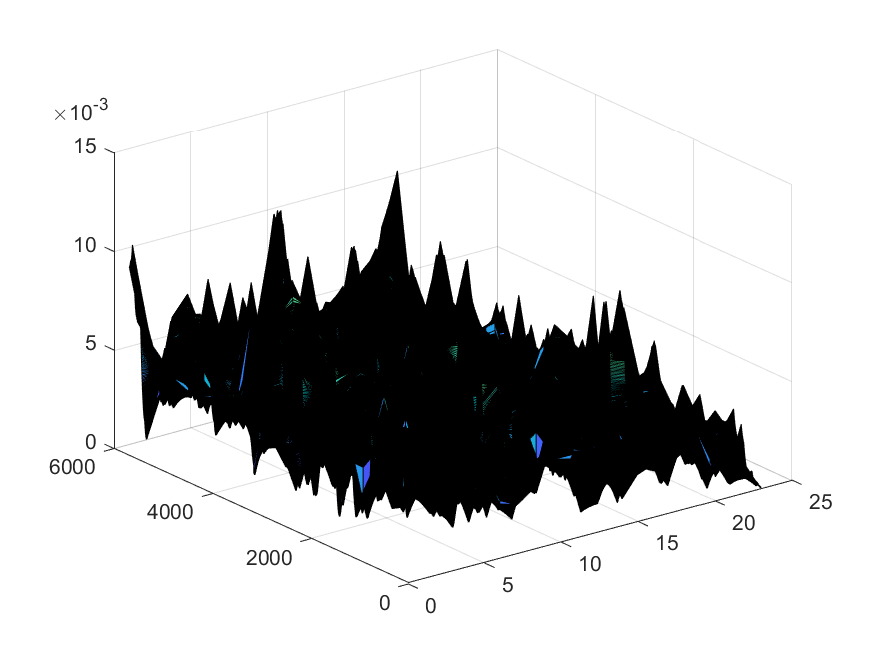}
\end{center}
\caption{$5$--fold    cross validation quadratic angular functional  errors  (December, 1979--February, 1980 at the top, and  March--May, 1980,  at the bottom).}
\label{fig:1rdam4}
\end{figure}

Similar cross validation   results are obtained for each month in the period  December, 1979--May, 1980, considering
the data set available at  NASA's National Space Science Data Center. Missing data affect
 the analysis of these months.   As commented, the sample size  for each month is  different, with  May  being the  most affected month   by missing data, although it has not prevented its inclusion in our analysis.  We have adopted the criterion that all curves should have 6000 consecutive time nodes.
For each month, the starting time  of the curves was taken randomly from the excess nodes of the corresponding multiple of 6000 nodes contained in days 3-5 of the month considered.
This, together with the existence of missing data, means that, a priori, monthly translations of curve samples  likely do not correspond to identical time intervals.
 Figures \ref{fig:1rdam1}--\ref{fig:1rdam2} display spherical bivariate  curve data   during these remaining months.  Empirical  intrinsic Fr\'echet functional  means  of regressors  are shown in Figure \ref{fig:1rdam3}, and
  5--fold  cross validation quadratic angular functional errors are provided  in Figure \ref{fig:1rdam4}.

\clearpage


\end{appendix}

\subsection*{Acknowledgements}
 \noindent This work has been supported in part by projects MCIN/ AEI/PID2022-142900NB-I00,
MCIN/ AEI/PGC2018-099549-B-I00,  and CEX2020-001105-M MCIN/AEI/ \linebreak 10.13039/501100011033).

\end{document}